\documentclass[11pt]{article}

\usepackage{geometry}
\usepackage{graphicx}
\usepackage{amssymb,amsmath,amsthm}
\usepackage{epstopdf}
\usepackage{enumerate}
\usepackage{enumitem}
\usepackage[english]{babel}
\usepackage{esint}
\usepackage{subcaption}
\usepackage{tikz}
\usepackage{xcolor}
\usepackage{pgfplots}
\pgfplotsset{compat=1.18}
\usetikzlibrary{calc}
\usetikzlibrary{calc,arrows.meta,decorations.pathreplacing}
\usepackage{hyperref}%
\hypersetup{colorlinks=true,linkcolor=blue!60!black,citecolor=blue!60!black,urlcolor=blue!60!black}

\RequirePackage{authblk}

\definecolor{ccteal}{HTML}{1D9E75}
\definecolor{ccviolet}{HTML}{7F77DD}
\definecolor{ccgreen}{HTML}{639922}

\newcommand\ccax{0.70}\newcommand\ccay{0.24}%
\newcommand\ccbx{0.30}\newcommand\ccby{0.66}%
\newcommand\ccT{7}%
\newcommand\ccS{8.1}

\newcommand\ccsetup{%
  \pgfmathsetmacro\ccdet{(\ccax)*(\ccby)-(\ccbx)*(\ccay)}%
  \pgfmathtruncatemacro\ccmax{round((\ccS)*(\ccby)/(\ccdet))}%
  \pgfmathtruncatemacro\ccmay{round(-(\ccS)*(\ccay)/(\ccdet))}%
  \pgfmathtruncatemacro\ccmbx{round(-(\ccS)*(\ccbx)/(\ccdet))}%
  \pgfmathtruncatemacro\ccmby{round((\ccS)*(\ccax)/(\ccdet))}%
  \pgfmathtruncatemacro\ccD{(\ccmax)*(\ccmby)-(\ccmay)*(\ccmbx)}%
  \ifnum\ccD<0
    \pgfmathtruncatemacro\ccmbx{-(\ccmbx)}\pgfmathtruncatemacro\ccmby{-(\ccmby)}%
    \pgfmathtruncatemacro\ccD{-(\ccD)}%
  \fi
  \pgfmathsetmacro\cclax{(\ccmax)*(\ccax)+(\ccmay)*(\ccbx)}%
  \pgfmathsetmacro\cclay{(\ccmax)*(\ccay)+(\ccmay)*(\ccby)}%
  \pgfmathsetmacro\cclbx{(\ccmbx)*(\ccax)+(\ccmby)*(\ccbx)}%
  \pgfmathsetmacro\cclby{(\ccmbx)*(\ccay)+(\ccmby)*(\ccby)}%
  \pgfmathsetmacro\ccd{max(veclen(\ccax+\ccbx,\ccay+\ccby),veclen(\ccax-\ccbx,\ccay-\ccby))}%
  \pgfmathsetmacro\ccAn{max(abs(\ccbx)+abs(\ccby),abs(\ccax)+abs(\ccay))/abs(\ccdet)}%
}
\newcommand\cctile[2]{%
  \edef\ccp{\the\numexpr(4*((#1)*(\ccmby)-(#2)*(\ccmbx))+1)/(4*\ccD)\relax}%
  \edef\ccq{\the\numexpr(4*((\ccmax)*(#2)-(\ccmay)*(#1))+1)/(4*\ccD)\relax}}
\newcommand\ccin[2]{\cctile{#1}{#2}%
  \edef\ccc{\ifnum\ccp=0 \ifnum\ccq=0 1\else 0\fi\else 0\fi}}
\newcommand\cccell[2]{(#1-.5,#2-.5) -- (#1+.5,#2-.5) -- (#1+.5,#2+.5) -- (#1-.5,#2+.5) -- cycle}
\tikzset{cclattice/.style={cm={\ccax,\ccay,\ccbx,\ccby,(0,0)}}}

\makeatletter
\newcommand{\tocnoheading}{\par\everypar{}\@endpefalse\@starttoc{toc}}
\makeatother

\newtheorem{definition}{Definition}[section]
\newtheorem{lemma}[definition]{Lemma}
\newtheorem{theorem}[definition]{Theorem}
\newtheorem{proposition}[definition]{Proposition}
\newtheorem{corollary}[definition]{Corollary}
\newtheorem{remark}[definition]{Remark}

\newcommand{\LM}[2]{\hbox{\vrule width.4pt \vbox to#1pt{\vfill
\hrule width#2pt height.4pt}}}
\newcommand{\LLL}{{\mathchoice {\>\LM{7}{5}\>}{\>\LM{7}{5}\>}{\,\LM{5}{3.5}\,}{\,\LM{3.35}{2.5}\,}}}

\def\WW{{\rm{W}}}
\def\LL{{\rm{L}}}
\def\CC{{\rm{C}}}
\def\R2{\mathbb{R}^{2}}
\def\M{\mathbb{R}^{2 \times 2}}
\def\dd{{\rm{d}}}

\title{Energies with values in networks\\ and geometric quasiconvexity}
\author[1]{Andrea Braides} 
\author[2]{Jan Kristensen}%\footnote{corresponding author} 
\author[3]{Sergio Scalabrino}\affil[1]{\small Department of Mathematics, University of Rome Tor Vergata, via della Ricerca Scientifica 1, Rome, Italy {\tt braides@mat.uniroma2.it}}
\affil[2]{\small Mathematical Institute, University of Oxford, Andrew Wiles Building, Radcliffe Observatory Quarter, Woodstock Road, Oxford, OX2 6GG, Oxford, United Kingdom
  {\tt kristens@maths.ox.ac.uk}}
\affil[3]{\small  SISSA, via Bonomea 265, Trieste, Italy {\tt sscalabr@sissa.it}}
\date{}

\def\RR{{\mathbb R}}
\def\ZZ{{\mathbb Z}}
\def\NN{{\mathbb N}}
\def\Z{{\mathcal Z}}
\def\e{\varepsilon}

\providecommand{\apmd}{\operatorname{ap\,md}}

\begin{document}
\maketitle
\begin{abstract}
We study the asymptotic behaviour, as $\varepsilon\to 0^+$, of the $p$-Dirichlet
energy of maps $u\colon\Omega\subset\R2 \to \R2$
constrained to take values in the rescaled periodic network $\varepsilon\mathcal Z$,
where $\mathcal Z$ is a standard square network in $\R2$.
The oscillations are carried by the constraint
in the target, not by the energy density, and their effect depends on $p$. For
$1<p<2$ the constraint disappears in the limit: the $\Gamma$-limit is an integral
functional on the whole of $\WW^{1,p}(\Omega , \R2 )$, whose density is given
by an asymptotic cell formula on tilted cells. This density is
equal to $(|a_1|+|a_2|)^p$ on rank-one matrices $a\otimes\nu$, $|\nu|=1$, and exceeds its
convex envelope by a quantity comparable to $|\det A|\,|A|^{p-2}$, with
constants of order $(2-p)^{-1}$. Consequently, although the original energy density is
convex, the limit density is quasiconvex but neither convex nor polyconvex:
quasiconvexity is generated by the geometry of the constraint. For $p\ge2$ the
constraint survives in the limit as $\det\nabla u=0$, and the $\Gamma$-limit is
the degenerate energy $\int_\Omega(|\nabla u_1|+|\nabla u_2|)^p\, \dd x$, that is,
the $p$-energy of $u$ considered as a map into the normed plane
$(\R2,|\cdot|_1)$, the metric limit of the networks.
\end{abstract}

\bigskip
\hfill {\em Dedicated to David Kinderlehrer on the occasion of his 85th birthday}
\bigskip

{\hypersetup{linkcolor=black}\tocnoheading}
\section{Introduction and main results}

In this paper, we study the asymptotic behavior of a prototypical family of variational integrals defined
on vector-valued Sobolev functions, and whose values are constrained to lie on a periodic network, which is seen
at a scale $\e$ that is small compared to the scale of observation. Given a bounded Lipschitz domain
$\Omega\subset\R2$ and an exponent $p>1$, we consider
\begin{equation}\label{eq1}
F_\e^p(u)=
\begin{cases}
\displaystyle\int_{\Omega} \! |\nabla u(x)|^p\, \dd x
&\hbox{if $u(x)\in\e\Z$ for almost every $x\in\Omega$}\cr
+\infty&\hbox{otherwise,}
\end{cases}
\end{equation}
where $\Z$ is a periodic network in $\R2$ and $| \cdot |$ denotes the Frobenius norm on $\M$. The prototypical example, that is
our focus here, is the union of the coordinate lines
\begin{equation}\label{eq2}
\Z=\Bigl\{(a,b)\in\R2 :\ a\in\tfrac12+\ZZ\ \hbox{ or }\ b\in\tfrac12+\ZZ\Bigr\}
=\bigcup_{i\in\ZZ^2}(i+\partial Q),\qquad Q=(-\tfrac12,\tfrac12)^2 .
\end{equation}
The effective domain of $F^p_\e$ is the class of constrained Sobolev maps 
$\WW^{1,p}(\Omega , \e\Z) = \bigl\{ u\in \WW^{1,p}(\Omega , \R2 ):\ u(x)\in\e\Z$
 for almost all  $x\in\Omega \bigr\}$.
Constraints of this type arise naturally in models of media with a fine
internal architecture, where admissible deformations, order parameters or
internal variables are forced onto a network of preferred states: lattice or
network-like microstructures, strongly anisotropic media, and discrete systems
in which the order parameter is quantized. 
Moreover, this can be regarded as a prototype for energies with values on an oscillating manifold, with the heuristic motivation
that low-energy functions will be located approximately on a network with small oscillations within the manifold (which is the case for curves \cite{BCC}). 
In all such models, one is led to ask what effective continuum theory survives the limit $\e\to 0^+$. 

We will refer to this limit as a homogenization process, although it is not a usual one.
Two features distinguish \eqref{eq1} from the energies usually encountered in
homogenization theory. First, the integrand is the $p$-Dirichlet
integrand: all oscillations are carried by the \emph{constraint}, not by the
energy density. Second, and more importantly, the oscillations take place in
the \emph{target} space rather than in the spatial $x$-variable. The
classical homogenization results deal with integrands $f(x/\e,\nabla u)$ oscillating in
$x$, and the effective energy is obtained through a cell formula on the domain
(see Braides \cite{Braides}, M\"uller \cite{Muller} and the monograph
\cite{BDF}); here, instead, the cell problem lives in the target, and the
relevant cell is the pre-image $Q_A=A^{-1}(Q)$ of the target cell under the
affine datum with gradient $A$. This is why, unlike in the classical setting, the
homogenization formula depends on the rank of the matrix $A$. More importantly,
although the original integrand is convex, the homogenized density is not: it
is quasiconvex but neither convex nor polyconvex, and the distance from its
convex envelope is controlled by the determinant. This lack of convexity is
produced by the geometry of the constraint alone and it is an important aspect of
the present work.

The main goal of the paper is to characterize the $\Gamma$-limit of $F^p_\e$
as $\e\to 0^+$ and to understand how the interplay between $p$ and the geometry
of $\Z$ affects the limit model. A dichotomy appears at $p=2$: for $p<2$ the
constraint disappears in the limit, but it leaves its trace in a genuinely
quasiconvex energy density; for $p > 2$ it survives as the differential
constraint $\det\nabla u=0$. 
For $p<2$ the result is the following.

\smallskip
\noindent{\bf Theorem A} (subcritical regime: geometric quasiconvexity).
{\em Let $1<p<2$. Then $F^p_\e$ $\Gamma$-converges, with respect to the
strong $\LL^p$-topology, to
$$
F^p_{\rm hom}(u)=\int_{\Omega} \! f^p_{\rm hom}(\nabla u(x))\, \dd x, \, u \in \WW^{1,p}( \Omega , \R2 ) ,
$$
where $f^p_{\rm hom}\colon\M \to[0,+\infty)$ is a continuous,
$p$-homogeneous, quasiconvex integrand with the following properties:

\smallskip
{\rm (i) (asymptotic homogenization formula)} on matrices of rank two
\begin{eqnarray}\label{ho-for}\nonumber
&&f^p_{\rm hom}(A)=\liminf_{T\to+\infty}\frac{|\det A|}{T^2}\inf\Bigl\{\int_{TQ_A} \! |\nabla u|^p\, \dd x: \, u\in \WW^{1,p}(TQ_A , \Z),\ \\
&& \hskip5cm
u= u_A \hbox{ \rm in a neighbourhood of } \partial TQ_A\Bigr\},
\end{eqnarray}
where  $u_A$ is the {\em projection} of the linear function $Ax$ on $\Z$ (see Definition~{\rm\ref{projlin}});

\smallskip
{\rm (ii) (qualified two-sided estimates)}
there exist $\overline c, \overline C>0$ independent of $p$ and explicitly computable,
such that
\begin{equation}\label{eq:bounds}
N(A)^p+\overline c\,(p-1)8^{\frac1{2(1-p)}}\,\frac{|\det A|}{(2-p)|A|^{2-p}}\ \le\ f^p_{\rm hom}(A)\ \le\
N(A)^p+\overline C\,\frac{|\det A|}{(2-p)|A|^{2-p}}
\end{equation}
for every $A\neq0$, where, if $A_1$ and $A_2$ denote the rows of $A$, 
\begin{equation}\label{eq:N}
N(A):=|A_1|+|A_2|.
\end{equation}
In particular, $f^p_{\rm hom}(A)=N(A)^p$ on rank-one matrices, and
the function $f^p_{\rm hom}$ is neither convex nor polyconvex.}

\smallskip
The density $f^p_{\rm hom}$ is given by an asymptotic cell formula on the
tilted cells $TQ_A=TA^{-1}(Q)$, with the projected linear map
$u_A=P_0(A\,\cdot\,)$ as boundary datum (Definitions~\ref{projlin} and
\ref{formula f_hom with A non singular}); by \eqref{eq:bounds} it has standard
$p$-growth, $|A|^p\le f^p_{\rm hom}(A)\le C|A|^p$. The restriction $p<2$ is
exactly what makes the radial projection onto the boundary of a cell an
admissible operation in $\WW^{1,p}$, the singularity $|x|^{-p}$ 
integrable in the plane precisely for $p<2$ (Lemma \ref{lemma1}): every
Sobolev map can be projected onto the network at a uniformly controlled
energy cost, which is why the limit domain is the whole of
$\WW^{1,p}(\Omega , \R2 )$. 

The lower bound in \eqref{eq:bounds} shows that $f^p_{\rm hom}$ is not convex.
Indeed, writing $I=\tfrac12\,(2e_1\otimes e_1)+\tfrac12\,(2e_2\otimes e_2)$,
by (i) and (ii) we have
\begin{equation}\label{f hom not convex counterexample}
\tfrac12 f^p_{\rm hom}(2e_1\otimes e_1)+\tfrac12 f^p_{\rm hom}(2e_2\otimes e_2)
=2^p\ <\ 2^p+\overline c\,\frac{p-1}{2-p}\,8^{\frac1{2(1-p)}}\,2^{\frac{p-2}2}\le f^p_{\rm hom}(I) .
\end{equation}
%which is compatible with rank-one convexity, the matrix
%$2e_1\otimes e_1-2e_2\otimes e_2$ having rank two.
As a consequence, $f^p_{\rm hom}$ is not polyconvex either: a polyconvex function with growth
of order $p<2$ must be convex, since an affine function of $(A,\det A)$ lying
below a function with subquadratic growth cannot depend on $\det A$.

\smallskip\noindent{\em Geometric quasiconvexity.} We consider $f^p_{\rm hom}$ as a quasiconvex function generated by the geometry
of the constraint, from which the title of this paper comes. It is worth commenting on its form.
The two terms in \eqref{eq:bounds} have a transparent meaning. The first 
is the convex envelope of $f^p_{\rm hom}$. Indeed, $N^p$ is convex and lies
below $f^p_{\rm hom}$; conversely, if $A_1,A_2\neq0$ and
$\lambda=|A_1|/N(A)$, then
$$
A=\lambda\Big(\frac1\lambda\,e_1\otimes A_1\Big)
+(1-\lambda)\Big(\frac1{1-\lambda}\,e_2\otimes A_2\Big)
$$
is a convex combination of two matrices of rank one at which, by (i),
$f^p_{\rm hom}$ takes the value $N(A)^p$. Hence $(f^p_{\rm hom})^{**}=N^p$,
and by (i) and \eqref{eq:bounds} the function $f^p_{\rm hom}$ coincides with
its convex envelope exactly on the matrices of rank at most one: it is convex
at each of them, in the sense that it has an affine minorant touching it
there, and at no matrix of rank two. Note that $N$ is a norm on $\M$ with $|A|\le N(A)\le\sqrt2\,|A|$ and
$N(a\otimes\nu)=|a_1|+|a_2|$ for every $a\in\R2$ and every unit vector
$\nu\in\R2$.
The value $N(A)^p$ is also identified in
Section \ref{section YM} as the minimal $p$-th moment of a probability
measure with barycentre $A$ supported on
\begin{equation}\label{eq:K}
K:=(e_1\otimes\R2 ) \cup( e_2\otimes\R2 )\subset \M,
\end{equation}
the set of matrices with at most one non-zero row, which contains the
gradients of all maps with values in $\e\Z$
(Lemma \ref{nabla u_e has only one row}). In other words, $N(A)^p$ is the best
lower bound that can be extracted from the pointwise constraint on the
gradients alone. For matrices of rank one, this bound is sharp, the optimal
measure being a laminate that we can generate with maps with values in $\e\Z$.
However, for $\det A\neq0$, instead, the unique optimal measure on $K$ is supported on two
matrices that are not rank-one connected, so it is not a homogeneous $\WW^{1,p}$-gradient Young
measure \cite{BJ}. The fact that the generating sequence must be gradients for maps valued in $\e\Z$ is
responsible for the nonconvexity of $f^p_{\rm hom}$. Indeed, the second term in \eqref{eq:bounds} is the energy
of the topological defects forced by the constraint. Maps with values in $\e\Z$ satisfy
$\det\nabla u_\e=0$ almost everywhere, while the affine datum has Jacobian
$\det A\neq0$; the difference is carried by \emph{vortices}, that is, by
points around which $u_\e$ winds once around a cell of the network. For 
projected linear maps, this can be computed explicitly: the distributional
Jacobian of $\e\,u_A(x/\e)$ is
$$
{\rm sign}(\det A)\,\e^2\sum_{j\in\ZZ^2}\delta_{\e A^{-1}j},
$$
with a Dirac mass, weighted by the area $\e^2$ of a
cell, at each point of a lattice of density $|\det A|\,\e^{-2}$. A vortex has
finite $p$-energy only if $p<2$, and its cost at the scale $\e/|A|$, below
which it dominates the affine part of the gradient, is of the order
$\e^2|A|^{p-2}/(2-p)$; multiplied by the density of vortices, this accounts
heuristically for the second term in \eqref{eq:bounds}, including the order
$(2-p)^{-1}$ of the constants.

The proof of the lower bound makes these heuristics rigorous without any
description of the vortices: the Jacobian of the boundary datum must be paid
for by large gradients. Precisely, a Lipschitz truncation \cite{AF} at level
$\lambda\ge4|A|$, combined with the null-Lagrangian property of the
determinant, shows that for every competitor $u$ in the cell problem the
maximal function of $|\nabla u-A|$ exceeds $\lambda$ on a set whose measure,
relative to that of the cell, is at least of order $|\det A|\,\lambda^{-2}$.
Integrating in $\lambda$ gives a lower bound of order
$|\det A|\,|A|^{p-2}/(2-p)$ for the $p$-th moment of the large values of
$|\nabla u|$ (Lemma \ref{jacobians pay gradients}), which is then added to
$N(A)^p$ by a Jensen inequality that keeps track of the tail
(Lemma \ref{jensen with tail lemma}). The upper bound follows from the
rank-one convexity of $f^p_{\rm hom}$, by laminating a rank-one matrix with a
multiple of a rotation, whose cost is computed through the projection $P_0$.
It is possible to give a more precise lower bound for the density $f^p_{\rm hom}$ using recent results \cite{burk26}.
First, we note that $N(A) \geq \| A \|$, where $\| \cdot \|$ denotes the
spectral norm, with equality for $A \in K$. Using that the Burkholder integrands defined for exponents $p \in (1,2)$ as
$$
\mathrm{B}_{p}^{\pm}(X) = \bigl( \tfrac{2}{p}-1 \bigr) \| X \|^{p} \pm \| X \|^{p-2}\det X , \, X \in \M ,
$$
are quasiconvex one may show that
\begin{equation}\label{burk-lower}
f^p_{\rm hom}(X) \geq \| X \|^{p}\,\bigg( 1 + \frac{p}{2-p} \frac{\bigl| \det X \bigr|}{\| X \|^2}\bigg)
\end{equation}
holds for all $X \in \M$. The lower bound in \eqref{burk-lower} equals
$\max \bigl\{ \mathrm{B}_{p}^{+}(X),\mathrm{B}_{p}^{-}(X) \bigr\}/(\tfrac{2}{p}-1)$ and as such is quasiconvex.

We now describe the result for $p > 2$.
By \eqref{eq:bounds}, the densities $f^p_{\rm hom}$ converge pointwise, as
$p\to 2^-$, to a function equal to $N(A)^2$ if $\det A=0$ and to
$+\infty$ otherwise, which is the limit density for $p=2$ given by Theorem B
below: the vortex energy becomes infinite and turns into a constraint.

\smallskip
\noindent{\bf Theorem B} (supercritical regime: rigidity).
{\em Let $p > 2$. Then $F^p_\e$ $\Gamma$-converges, with respect to the strong
$\LL^p$-topology, to
$$
F^p_0(u)=
\begin{cases}
\displaystyle\int_\Omega\bigl(|\nabla u_1(x)|+|\nabla u_2(x)|\bigr)^p\,\dd x
&\hbox{if $u\in \WW^{1,p}(\Omega , \R2)$ and $\det\nabla u=0$ a.e.}\cr
+\infty&\hbox{otherwise.}
\end{cases}
$$}

\smallskip
In the supercritical regime, the constraint produces a rigidity phenomenon that
survives in the limit: limit deformations have gradients of rank at most one,
and the effective energy is the degenerate energy
$N(\nabla u)^p=(|\nabla u_1|+|\nabla u_2|)^p$ induced by the network metric.
The limit density is a convex function of $(A,\det A)$; that is, it is polyconvex in
the sense of Ball \cite{Ball}, with extended values. Crossing $p=2$, the limit
density thus passes from a finite, genuinely quasiconvex function to a
polyconvex, constrained one; the two regimes agree on matrices of rank one,
where both densities equal $(|a_1|+|a_2|)^p$ for $A=a\otimes\nu$, $|\nu|=1$.
For functions $u$ with an image a curve, the proof of Theorem B 
(Section \ref{section p>2}) is short: the lower bound
follows from the identity $|\nabla u_\e|=N(\nabla u_\e)$, valid almost
everywhere for maps with values in $\e\Z$, from the convexity of $N^p$ and
from the weak continuity of the determinant, while recovery sequences are
obtained by composing $u$ with maps that replace segments in its image by
staircases in $\e\Z$.

\medskip
We want to remark how the two ways in which the exponent $2$ appears are the two sides of the
same coin.
From below, $p<2$ is the integrability condition that makes the projection of the radial projection
onto $\partial Q$ bounded in $\WW^{1,p}$ (Lemma \ref{lemma1}). Averaging over
the projection centers goes back to Federer and Fleming \cite{FF}, and was
introduced for Sobolev maps with values in manifolds by Hardt, Kinderlehrer, 
and Lin \cite{HKL}; here the network plays the role of the one-dimensional
skeleton of the cubical decomposition of the target, and projecting cell by
cell allows us to destroy the constraint in the limit.

From above, $p\ge2$ is what guarantees that the pointwise constraint
$\nabla u_\e\in K$ passes to the weak limit as the differential constraint
$\det\nabla u=0$. The determinant vanishes on $K$ and is a null Lagrangian, so
it is affine along rank-one segments, and $\{\det=0\}$ is the lamination
convex hull of $K$. For $p\ge2$ the pointwise Jacobian of a $\WW^{1,p}$ map
coincides with its distributional Jacobian, which is weakly continuous
\cite{Ball}, so that $\det\nabla u_\e=0$ passes to the limit. Below $p=2$ the
distributional Jacobian, although still weakly continuous for $p>4/3$, may
differ from the pointwise one by concentrations at vortex points, and maps
with values in $\e\Z$ can carry such vortices: the Jacobians of the maps
$\e\,u_A(x/\e)$ computed above converge weakly$^*$ to $\det A$, although
$\det\nabla(\e\,u_A(x/\e))=0$ almost everywhere. This is the phenomenon
exhibited, e.g.~by Ball and Murat \cite{BMu}, with the map $x\mapsto x/|x|$, which
shows that in dimension $n$ the determinant is not $\WW^{1,p}$-quasiconvex for
$p<n$. In our
setting, it means that for $p<2$ the inclusion $\nabla u\in K$ is compatible
with every affine boundary datum: $u_A$ satisfies $\nabla u_A\in K$ almost
everywhere and $u_A(x)=Ax$ on $\partial(TQ_A)$ for every odd integer $T$,
while for $p\ge2$ such a map exists only if $\det A=0$.

The dichotomy also has a topological reading. For $p>2$ every
$u\in \WW^{1,p}(\Omega , \e\Z)$ is continuous and, if $\Omega$ is simply
connected, it lifts to the universal covering of $\e\Z$, which is a tree, so
 no vortex can occur. For $p<2$, instead, the maps that wind around a cell, such
as $x\mapsto\e\,x/(2|x|_\infty)$ near the origin, have finite $p$-energy, and
 are precisely the singularities created by the projection of
Lemma \ref{lemma1} at the points where $u$ meets the projection centers. The
borderline case $p=2$ is the analog, for graph-valued maps, of the lifting
problem for circle-valued maps \cite{BN,BBM}, which we do not pursue here.

For $p\ge 2$ Theorem B can also be interpreted as a study of the limits of some metric spaces.
We refer to Section \ref{metric remarks} for a more detailed discussion of this aspect.

\smallskip\noindent{\bf Relation to the literature}

\noindent
{\em Homogenization under constraints on the target.}
The homogenization of integral functionals for maps constrained to a fixed smooth
manifold has been studied by Babadjian and Millot \cite{BM1}, who
introduced a notion of tangential homogenization modeled on the tangential
quasiconvexity of Dacorogna, Fonseca, Mal\'y and Trivisa \cite{DFMT}.
In those works the target is fixed and smooth and the oscillation is in the
space variable, so that our setting is in a precise sense complementary.

\smallskip\noindent{\em Oscillating targets and effective metrics.}
The closest antecedent is the analysis of Braides, Cancedda and Chiad\`o Piat
\cite{BCC}, who consider the Dirichlet integral of \emph{curves} taking values
in fast-oscillating manifolds collapsing onto a linear subspace, and show that
the limit is described by a Finsler metric obtained by optimizing the
contribution of the oscillations on level sets, with formulas mixing
homogenization and convexification. Their result is the one-dimensional
domain counterpart of our supercritical regime, in which the oscillating target
degenerates into a metric. Related constructions of effective metrics by
$\Gamma$-convergence appear in Braides, Buttazzo, and Fragal\`a \cite{BBF}.

\smallskip\noindent{\em Quasiconvexity produced by homogenization.}
Quasiconvexity in the sense of Morrey \cite{Morrey} is the natural convexity
condition for vectorial integral functionals, and polyconvexity in the sense
of Ball \cite{Ball} is its most widely used sufficient condition; see
\cite{burk26,GC26,Dacorogna,GK22} for the relations between the various notions. Explicit
quasiconvex functions that are not polyconvex are plentiful, even in the
$2\times2$ case, but explicit examples are fewer. The early examples include \cite{AD,Sverak},
and homogenization is a natural source of them: Braides \cite{Braides-poly} showed that polyconvexity is not preserved
by periodic homogenization. The mechanism we exhibit is of a different
nature: the energy density is convex and fixed, and it is the non-convex
constraint in the target that generates a quasiconvex, non-convex limit
density.

\smallskip\noindent{\em Homogeneous quasiconvex functions of low degree.}
On $\RR^{2\times2}$, functions positively homogeneous of degree $p\ge2$ that
are quasiconvex but not convex are easy to produce: $|\det A|^{p/2}$ is
polyconvex and not convex. For $p<2$ the picture changes since, as observed
after Theorem A, polyconvex functions with growth of order $p$ are convex.
That quasiconvex functions with subquadratic growth need not be convex was first observed
by \v{S}ver\'ak \cite{Sverak-sub}. Positively homogeneous functions of
degree one were studied by M\"uller \cite{Muller-hom} and by Dacorogna and
Haeberly \cite{DH}; in this class, Kirchheim and Kristensen \cite{KK} proved
that rank-one convexity forces convexity at every matrix of rank at most
one, while the examples of \cite{Sverak-sub,Muller-hom,DH} show that
convexity may fail at matrices of rank two and even quasiconvexity may fail
at rank two matrices in higher dimensions \cite{GC26}. For every $p\in(1,2)$, the
density $f^p_{\rm hom}$ is a quasiconvex function homogeneous of degree $p$
that exhibits this picture in its sharpest form, and it arises naturally from
a homogenization problem: it is convex at every matrix of rank at most one
and at no matrix of rank two, and its distance from its convex envelope is
comparable to $|\det A|\,|A|^{p-2}$, with constants of order $(2-p)^{-1}$.

\smallskip\noindent{\em Vortices.}
The energetics of point singularities is classical for maps into the circle:
the Dirichlet energy of a vortex diverges logarithmically, which is at the
origin of the Ginzburg--Landau asymptotics \cite{BBH}, while for $p<2$ the
$p$-energy of a vortex is of order $(2-p)^{-1}$, and Hardt and Lin \cite{HL95}
showed that $p$-energy minimizing maps from planar domains into the circle,
with boundary data of degree $d$, develop exactly $|d|$ singularities, whose
behaviour as $p\to 2^-$ they described. The constants in \eqref{eq:bounds}
display the same order $(2-p)^{-1}$, the Jacobian, rather than the degree,
measuring the number of vortices.

\smallskip\noindent{\em Differential inclusions and gradient Young measures.}
The set $K$ in \eqref{eq:K} is the union of two two-dimensional subspaces of
$\RR^{2\times2}$; the rank-one segments joining them lie in $\{\det=0\}$,
which is both the lamination convex hull and the quasiconvex hull of $K$.
In the moment problem of Section \ref{section YM}, whose value is $N(A)^p$,
the optimal measure is a laminate when $A$ has rank one and, by Ball and
James \cite{BJ}, it is not a gradient Young measure when $\det A\neq0$.
Which measures supported on $K$ are generated by
gradients of maps with values in $\e\Z$, a constrained version of the
characterization of gradient Young measures in Sobolev spaces by Kinderlehrer
and Pedregal \cite{KP}, is a question of independent interest. For the general
theory of differential inclusions and of hulls of sets of matrices we refer to
\cite{Gromov,MS,Kirchheim,Szekelyhidi}.

\smallskip\noindent{\em Maps into singular spaces and normed planes.}
At fixed $\e$ the target is a one-dimensional metric graph, so that the
variational theory of maps into singular targets applies: Gromov and Schoen
\cite{GS}, Korevaar and Schoen \cite{KS}, Eells and Fuglede \cite{EF},
Daskalopoulos and Mese \cite{DM}. Its estimates, however, are tied to the
curvature scale of the target. The shortest embedded circles of $\e\Z$ are the
boundaries of the cells, of length $4\e$, so that $(\e\Z,d_\e)$ is
CAT$(\kappa)$ exactly for $\kappa\ge\pi^2/(4\e^2)$ \cite[Part~II]{BH}; the
regularity theory for maps into targets with curvature bounded above by
$\kappa>0$ requires the image to lie in a ball of radius smaller than
$\pi/(2\sqrt\kappa)$, which for the optimal value $\kappa=\pi^2/(4\e^2)$ equals
$\e$ (see \cite{HKW} and, for singular spaces,
\cite{Serbinowski,GGZZ,Gennaioli}), and therefore gives no estimate uniform
in $\e$. In the limit, the convergence $d_\e\to|\cdot|_1$ is an
instance of the stable norm of a periodic metric space \cite{Pansu,Burago},
which would replace $|\cdot|_1$ for a general connected periodic network.
Finally, for $p=2$ the limit functional of Theorem B is the energy of a map
into a normed plane, a setting in which energy and area interact subtly:
several notions of area coexist, and Lytchak and Wenger \cite{LW1,LW2} showed
that energy minimization solves the Plateau problem once the notion of area
is chosen appropriately, quasiconvexity of the area being a key ingredient.
Comparing the quasiconvexity of these notions of area with that of our limit
densities seems an interesting question, which we do not pursue here.

\smallskip
{\bf Plan of the paper:}
Section \ref{section p<2} is devoted to the subcritical case $p<2$, containing the proof of Theorem A(i), although the construction of staircase paths in Section~\ref{section paths} is also used in the regime $p\ge 2$. In the subcritical regime, a fundamental issue is the estimate of projections on the networks (Section \ref{sec:proj}), together with the adaptation of nonlinear-homogenization techniques to functionals constrained in the target space; in particular, the modification of boundary values in Section \ref{sec:mo-ba}. In Section \ref{sec:en-fu} we separately define the homogenized energy densities on matrices of rank two and of rank one, respectively. In the first case, through an asymptotic homogenization formula, in the second case, by optimizing over paths. This is suggested by the necessity of two different variants of the blow-up method
applied in Section \ref{subsection gamma liminf}. Moreover, since the lower bound gives two different formulas, we prove that the resulting function satisfies the continuity properties necessary to prove the $\Gamma$-limsup inequality through a density argument in Section \ref{sec:ga-lims}.
Section \ref{section YM} develops the Young-measure counterpart of this analysis, identifying $N(A)^p$ as the value of a moment problem on $K$ and constructing, for rank-one $A$, the optimal laminate using maps with values in $\e\Z$.
Section \ref{section bounds on f hom} contains the proof of Theorem A(ii); that is, of the two-sided estimate
\eqref{eq:bounds}. Section \ref{section p>2} treats the case $p\ge2$ (Theorem B). In this case, the functions in the limit domain satisfy $\det\nabla u=0$. The necessary notions from Topology and Graph Theory are presented and used to overcome the difficulties due to the rigidity of this constraint, which forbids the usual arguments by density.
Section \ref{metric remarks} concludes the paper with some remarks on the connections with convergence issues of metric spaces.

\section{The subcritical case \texorpdfstring{$p<2$}{p<2}}\label{section p<2}
As noted in the Introduction the behaviour of functionals \eqref{eq1} is different in the two cases $p<2$ and $p\ge 2$. We start with the case $p<2$.

We set $Q=(-1/2,1/2)^2$, so that the network $\mathcal Z$ in \eqref{eq2} can be written as
$$
\Z=\bigcup_{i\in\ZZ^2} (i+\partial Q).
$$
For $a \in \RR$ we let $\lfloor a \rfloor$ denote the largest number in $\frac{1}{2}+\mathbb Z$ less than or equal to $a$.

The asymptotic part of Theorem A in the Introduction is the following result.
\begin{theorem}[subcritical homogenization theorem]\label{thm:sub:main}
Let $1<p<2$ and let $\Omega\subset\mathbb R^2$ be a Lipschitz bounded oper set. Let $\mathcal Z\subset\mathbb R^2$ be given by \eqref{eq2}, and let $F^p_\e$ be the functionals defined in 
\eqref{eq2}, with domain $\WW^{1,p}(\Omega,\e\mathcal Z)$. Then the $\Gamma$-limit as $\e\to 0$ with respect to the strong $\LL^p(\Omega)$-convergence and weak
$\WW^{1,p}(\Omega , \R2)$-convergence  of $F^p_\e$ exists and it is given by
\begin{equation}\label{gamma limit formula expression}
F^p_{\rm hom}(u)= \int_\Omega f^p_{\rm hom}(\nabla u(x))\dd x
\end{equation}
with domain $\WW^{1,p}(\Omega , \R2)$. The  quasiconvex function $f^p_{\rm hom} \colon \M \to [0,+\infty)$ is characterized on matrices $A$ with full rank by the 
{\em asymptotic homogenization formula~\eqref{ho-for}}.
\end{theorem}

This section will be devoted to the proof of Theorem \ref{thm:sub:main}. To that end, we will first study the properties of projections on $\Z$,
that ensure that the $\Gamma$-limit is finite on the whole $\WW^{1,p}(\Omega , \R2)$, and allow us to state the boundary conditions in the homogenization formula. 
A crucial step is the characterization of  $f^p_{\rm hom}$ on the subset of matrices of rank at most one, showing that $f^p_{\rm hom}(a\otimes\nu)=(|a_1|+|a_2|)^p$
for all $a,\nu\in\RR^2$ with $|\nu|=1$, or in other words that this formula coincides with the extension by continuity of the formula for full-rank matrices.
This is one the new technical parts of the arguments used. The liminf inequality is then proved by a blow-up argument, using the two different characterizations,
and the limsup inequality by using the correctors given by the homogenization formula. 

Further properties of $f^p_{\rm hom}$, which imply that $f^p_{\rm hom}$ is not convex, will be proved in Section \ref{section bounds on f hom}, completing the proof of Theorem A.

\subsection{Projections}\label{sec:proj}
In the case $p<2$ every function $u\in \WW^{1,p}(\Omega , \R2)$ can be ``projected'' to a $\WW^{1,p}$ function with values in $\e\mathcal Z$.
This is done by locally projecting on the boundary of an $\e$-cube from a suitable point in the interior.
In this section, we introduce the necessary notation.

For all $c\in\frac{1}{2}Q=(-1/4,1/4)^2$ and $i\in \ZZ^2$ we set 
$$
Q^\e_i=\e i+\e Q,\qquad x^{\e,i}_c=\e i+\e c,
$$
and consider the {\em projection 
from $x^{\e,i}_c$ to $\partial Q^\e_i$}
$
P^{\e,i}_c\colon Q^\e_i\setminus\{x^{\e,i}_c\} \to \partial Q^\e_i
$
defined by 
$$
P^{\e,i}_c(x)= y\hbox{ such that } \{y\}=\{x^{\e,i}_c+t(x-x^{\e,i}_c): t>0\}\cap \partial Q^\e_i.
$$
This projection is continuous in $Q^\e_i\setminus\{x^{\e,i}_c\}$ and smooth outside the
four segments that join $x^{\e,i}_c$ to the vertices of $Q^\e_i$. Setting $P_c:=P^{1,0}_c$,
for fixed $c$ we have
$$
P^{\e,i}_c(x)=P^{\e,0}_c(x-\e i)+\e i=\e P^{1,i}_c\Big(\frac{x}{\e}\Big)
=\e i+\e P_c\Big(\frac{x}{\e}-i\Big),
$$
so that, for $x$ outside these segments,
$$
|\nabla P^{\e,i}_c(x)|=\Big|\nabla P_c\Big(\frac{x}{\e}-i\Big)\Big|
\le \frac{5}{2|\frac{x}{\e}-i-c|}= \frac{5\e}{2\,|x-x^{\e,i}_c|},
$$
where we used $|\nabla P_c(x)|\leq \frac{5}{2|x-c|}$. To prove this estimate, 
one can divide the square $Q$ in the four triangles given by joining $c$ with the four vertices of $Q$.
In the interior of each triangle, its definition can be found explicitly; for example, in the right triangle $P_c(x)=(\frac{1}{2}, c_2+\frac{x_2-c_2}{x_1-c_1}(\frac{1}{2}-c_1))$,
hence the gradient is immediate. By reversing roles and signs for the four triangles, it holds the (sharp) bound
\begin{gather*}
  |\nabla P_c(x)| \leq \frac{K(c)}{|x-c|}, \\ K(c):=\max\bigg\{\frac{(\frac{1}{2}+|c_2|)^2}{\frac{1}{2}-|c_1|}+\frac{1}{2}-|c_1| \,,\,
  \frac{(\frac{1}{2}+|c_1|)^2}{\frac{1}{2}-|c_2|}+\frac{1}{2}-|c_2| \bigg\}\leq \frac{5}{2},  
\end{gather*}
which follows by choosing $|c_1|=|c_2|=\frac{1}{4}$, corresponding to the four corners of $\frac{1}{2}Q$. 

Alternatively, for $c\in\frac12 Q$ and $x$ out of the four segments joining $c$ to the corners, let $\nu$ be the outer normal to the side of $Q$ met by the
ray from $c$ through $x$, and let $h:=\frac12-c\cdot\nu\ge\frac{1}{4}$ be the distance of $c$ from that side.
Then $P_c(x)=c+\frac{h}{(x-c)\cdot\nu}\,(x-c)$, and a direct computation gives
$$
|\nabla P_c(x)|=\frac{h\,|x-c|}{((x-c)\cdot\nu)^2}=\frac{|P_c(x)-c|^2}{h\,|x-c|}\le\frac{9}{2}\,\frac{1}{|x-c|},
$$
since $|P_c(x)-c|\le\frac{3\sqrt2}{4}$.
\begin{lemma}[projection on $\e\mathcal Z$]\label{lemma1}
For all $u\in \WW^{1,p}(\Omega,\RR^2)$ there exists a function $P_\e(u)\in \WW^{1,p}(\Omega,\e\Z)$
such that $\|P_\e(u)-u\|_\infty\le\e$, $P_\e(u)(x)=u(x)$ if $u(x)\in\e\Z$,  and 
$$
\|\nabla P_\e(u)\|_p\le C\|\nabla u\|_p,
$$
with $C=C(p)$ depending only on $p$.
\end{lemma}
\begin{proof}
With fixed $c\in \frac{1}{2}Q$ we define 
$$
P^\e_c(u(x))=P^{\e,i}_c(u(x)) \hbox{ if }u(x)\in Q^\e_i\setminus\{x^{\e,i}_c\}.
$$
Note that there is no ambiguity if $\partial Q^\e_i\cap \partial Q^\e_j$
since in that case $P^\e_c(u(x))=u(x)$. Moreover, if $u(x)\in Q^\e_i\setminus\{x^{\e,i}_c\}$ then 
\begin{equation}
|\nabla P^\e_c(u)(x)|\le \e \frac{C}{|u(x)-x^{\e,i}_c|}|\nabla u(x)|.
\end{equation}

Let 
$\Omega^\e_i=\{x\in\Omega: u(x)\in Q^\e_i\}$.
Then, since $|\frac{1}{2}Q|=\frac{1}{4}$,
\begin{eqnarray*}
4\int_{\frac{1}{2}Q}\int_{\Omega^\e_i} |\nabla P^\e_c(u)|^p\dd x\,\dd c
&\le& C\e^p\int_{\frac{1}{2}Q}\int_{\Omega^\e_i}\frac{1}{|u(x)-x^{\e,i}_c|^p}|\nabla u|^p\dd x\,\dd c
\\
&=& C\e^p\int_{\Omega^\e_i}\int_{\frac{1}{2}Q}\frac{1}{|u(x)-x^{\e,i}_c|^p}\dd c\,|\nabla u|^p\dd x
\le C\int_{\Omega^\e_i}|\nabla u|^p\dd x.
\end{eqnarray*}
Note that here we use $p<2$ so that $c\mapsto \frac{1}{|c|^p}$ is integrable
and 
$$
\int_{\frac{1}{2}Q}\frac{1}{|u(x)-x^{\e,i}_c|^p}\dd c\le C\e^{-p}
$$
for all $x$.
We then deduce that for all $i\in\ZZ^2$ there exists $c_i\in \frac{1}{2}Q$ such that
\begin{eqnarray*}
\int_{\Omega^\e_i} |\nabla P^\e_{c_i}(u)|^p\dd x
&\le& 
4\int_{\frac{1}{2}Q}\int_{\Omega^\e_i} |\nabla P^\e_c|^p\dd x\,\dd c
\le C\int_{\Omega^\e_i}|\nabla u|^p\dd x.
\end{eqnarray*}
If we set
\begin{equation}\label{pid}
P_\e(u)(x)= P^{\e,i}_{c_i}(u(x)) \hbox{ if }u(x)\in Q^\e_i
\end{equation}
then $P_\e(u)$ satisfies the desired estimate.\end{proof}

\begin{definition}[projection of a linear function]\label{projlin} For all $A\in\RR^{2\times 2}$, we define the {\em projection} $u_A(x)= P_0(A x)$, where $P_0$ is constructed as
in \eqref{pid} with $\e=1$ and $c_i=0$ for all $i$. In particular, denoting $Q_z:=z+Q, \ z \in \mathbb Z^2$, the explicit form of $P_0$ reads
\begin{equation}\label{explicit form of P_0}
P_0(x)= z+ \frac{1}{2}\frac{x-z}{|x-z|_{\infty}}, \quad x \in \overline{Q}_z\setminus z.    
\end{equation}
We will also define the projection $\psi(x):=x/|x|_\infty, \ x \in \overline{Q}\setminus \{0\}$, so that
$$
P_0(x)=z+ \frac{1}{2}\psi(x-z), \quad x \in \overline{Q}_z\setminus z. 
$$
\end{definition}

If $x\in\partial Q_z$ then we have $|x-z|_\infty=1/2$, so that $P_0(x)=x$,  $P_0$ coincides with the identity on the net, and $P_0$ is $\mathcal{Z}$-valued.
Moreover, for every $z \in \mathbb Z^2$ the bound $|\nabla P_0(x)|\leq |x-z|^{-1}$ holds for almost every $x \in z+Q$. 
Note that $P_0(Ax) \neq P_{1}(Ax)$ as given by Lemma \ref{lemma1}, since the latter is defined with varying centers $c_i$ with respect to the squares $Q_i$.
The former is just given by the composition of $P_0$ and the linear function.

\begin{lemma}\label{rescaled P_0 satisfies  C|U| bound}
Given $u_A$ as in Definition {\rm \ref{projlin}}, and set $C_A:=\int_{Q}|(\nabla P_0)(x)A|^p \dd x$, the rescaled function $\overline u_\e(x):=\e u_A(x/\e)$ satisfies
$$
\lim_{\e \to 0} \int_U |\nabla \overline u_\e|^p \dd x = C_A|U|
$$
for all fixed bounded measurable subsets $U$.
\end{lemma}
\begin{proof}
By the change of variables 
$$
\int_U |\nabla \overline{u}_\e(x)|^p\dd x= \int_U |(\nabla P_0)(Ax/\e)A|^p\dd x=\frac{1}{|\det A|}\int_{AU} |(\nabla P_0)(x/\e)A|^p\dd x.
$$
Since $x \mapsto |(\nabla P_0)(x/\e)A|^p$ is the rescaling of a $\mathbb{Z}^2$-periodic function in  $\LL^1(Q)$ (since $p<2$), it converges  weakly in $\LL^1_{\rm loc}$ to its average $C_A$.
Hence,
\begin{eqnarray*}
\lim_{\e \to 0}\int_U |\nabla \overline u_\e|^p\, \dd x &=& \lim_{\e \to 0}\frac{1}{|\det A|} \int_{AU} |(\nabla P_0)(x/\e)A|^p\,\dd x\\
&=& \frac{|AU|}{|\det A|}\int_{Q}|(\nabla P_0)(x)A|^p\, \dd x = C_A |U|
\end{eqnarray*}
by the Riemann-Lebesgue lemma. \end{proof}

\subsection{Modification of boundary conditions}\label{sec:mo-ba}
As often is the case in variational problems, it will be useful to construct converging sequences satisfying suitable boundary conditions
and with controlled energy. Usually, such sequences match the same boundary conditions as their limit. In this case, this is not possible,
since we have varying constraints, and the result must be stated correspondingly.

\begin{lemma}\label{lemma variation of BCs}
Let $\overline u_\e\to u$ be such that $\overline u_\e\in \WW^{1,p}(\Omega , \e\Z)$ and
there exists a constant $C$ such that
$$
\limsup_{\e\to 0} \int_U|\nabla \overline u_\e|^p\dd x\le C|U|
$$
for all fixed open sets $U$ with Lipschitz boundary. Then for all $u_\e\in \WW^{1,p}(\Omega , \e\Z)$ 
with $u_\e\to u$ there exists $v_\e\in \WW^{1,p}(\Omega , \e\Z)$
with 
$$
\liminf_{\e\to0}F_\e(v_\e)\le\liminf_{\e\to 0} F_\e( u_\e),
$$
and $v_\e=\overline u_\e$ on a neighbourhood of $\partial\Omega$. 
\end{lemma}

\begin{proof}
It suffices to consider the case $\liminf_{\e\to 0} F^p_\e( u_\e)<+\infty$,
and prove that for all $\eta>0$ we may find $v_\e\in \WW^{1,p}(\Omega, \e\Z)$
with 
$$
\liminf_{\e\to0}F^p_\e(v_\e)\le(1+\eta)\liminf_{\e\to 0} F^p_\e( u_\e) +\eta
$$
and $v_\e=\overline u_\e$ on a neighbourhood of $\partial\Omega$.
For the construction of $v_\e$ we may follow a well-known procedure 
by De Giorgi. With fixed $\delta>0$ and $N\in\NN$, $N>0$, we consider
the sets
$$
\Omega^\delta_i=\Bigl\{x\in\Omega: {\rm dist}\,(x,\partial\Omega)\le\frac{i}{N}\delta
\Bigr\}
$$
for $i\in\NN$, and the functions 
$$
\varphi_i(x)=\Bigl(\frac{N}{\delta}{\rm dist}\bigl(x,\Omega^\delta_i\bigr)\Bigr)\wedge 1.
$$
We have $\varphi_i\in \WW^{1,\infty}(\Omega)$, $0\le \varphi_i\le 1$, $\varphi_i=0$ on $\Omega^\delta_i$
and $\varphi_i=1$ on $\Omega\setminus\Omega^\delta_{i+1}$; in particular,
$\nabla \varphi_i\neq 0$ only in $\Omega^\delta_{i+1}\setminus\Omega^\delta_i$.
We set $\Omega^\delta=\Omega^\delta_N$.

We consider the functions
$$
w^i_\e= P_\e\big(\varphi_i u_\e+ (1-\varphi_i)\overline u_\e\big),
$$
where $P_\e$ is defined in Lemma \ref{lemma1}. Note that 
$w^i_\e=u_\e$ on $\Omega\setminus \Omega^\delta_{i+1}$ and 
$w^i_\e=\overline u_\e$ on $\Omega^\delta_i$.
Using Lemma \ref{lemma1}, we can first estimate 
\begin{eqnarray*}
\int_\Omega|\nabla w^i_\e|^p\,\dd x
&=&\int_{\Omega^\delta_i}|\nabla \overline u_\e|^p\,\dd x
+\int_{\Omega^\delta_{i+1}\setminus\Omega^\delta_i}|\nabla w^i_\e|^p\,\dd x
+\int_{\Omega\setminus \Omega^\delta_{i+1}}|\nabla u_\e|^p\,\dd x
\\
&\le&\int_{\Omega^\delta_i}|\nabla \overline u_\e|^p\,\dd x
+C_p\int_{\Omega^\delta_{i+1}\setminus\Omega^\delta_i}|\nabla(\varphi_i u_\e+ (1-\varphi_i)\overline u_\e)|^p\,\dd x\\
&&
+\int_{\Omega\setminus \Omega^\delta_{i+1}}|\nabla u_\e|^p\,\dd x
\\
&\le&\int_{\Omega^\delta_i}|\nabla \overline u_\e|^p\,\dd x
+C_p\int_{\Omega^\delta_{i+1}\setminus\Omega^\delta_i}\bigl(|\nabla\varphi_i|^p|u_\e-
 \overline u_\e|^p+ |\nabla \overline u_\e|^p+  |\nabla u_\e|^p\bigr)\dd x
\\
&&+\int_{\Omega\setminus \Omega^\delta_{i+1}}|\nabla u_\e|^p\,\dd x
\\
&\le& C_p\int_{\Omega^\delta}|\nabla \overline u_\e|^p\,\dd x +
C_p\int_{\Omega^\delta_{i+1}\setminus\Omega^\delta_i}|\nabla u_\e|^p\,\dd x\\
&&
+C_p\int_{\Omega^\delta_{i+1}\setminus\Omega^\delta_i}\Bigl(\frac{N}{\delta}\Bigr)^p|u_\e-
 \overline u_\e|^p\dd x
+\int_{\Omega}|\nabla u_\e|^p\,\dd x,
\end{eqnarray*}
with $C_p$ a positive constant depending only on $p$.
We can now define $w_\e= w^i_\e$, where  $i\in \{1,\ldots N-1\}$ is such that 
$$
\int_{\Omega^\delta_{i+1}\setminus\Omega^\delta_i}|\nabla u_\e|^p\,\dd x\ 
\le\  \frac{1}{N-1} 
\int_{\Omega^\delta}|\nabla u_\e|^p\,\dd x
\ \le\ \frac{1}{N-1} 
\int_{\Omega}|\nabla u_\e|^p\,\dd x,
$$
so that we have 
\begin{eqnarray*}
\int_\Omega|\nabla w_\e|^p\,\dd x
&\le& C_p\int_{\Omega^\delta}|\nabla \overline u_\e|^p\,\dd x\\
&&+\Bigl(1+\frac{C_p}{ N-1}\Bigr)\int_{\Omega}|\nabla u_\e|^p\,\dd x
+C_p\int_{\Omega}\Bigl(\frac{N}{\delta}\Bigr)^p|u_\e-
 \overline u_\e|^p\dd x.
\end{eqnarray*}
By the hypotheses on $\overline u_\e$, passing to the liminf as $\e\to0$ we obtain
$$
\liminf_{\e\to0} \int_\Omega|\nabla w_\e|^p\,\dd x
\le 
C|\Omega^\delta|+\Bigl(1+\frac{C_p}{N-1}\Bigr)\liminf_{\e\to 0}\int_{\Omega}|\nabla u_\e|^p\,\dd x,
$$
which produces the desired estimate, choosing $\delta$ and $N$ appropriately.
\end{proof}
Note that the assumption on $\overline{u}_\e$ of Lemma \ref{lemma variation of BCs} is satisfied by
the rescaled function $\e\, u_A(x/\e)$ in Lemma \ref{rescaled P_0 satisfies  C|U| bound}. This fact will be used repeatedly in the following.

\subsection{Paths in the network}\label{section paths}
In this section $\e>0$ is fixed, and we analyze paths with values in $\e\Z$.
Such paths will be interpreted as sections of linear functions in the rank-one case.

\smallskip
We call {\em vertices} of $\e\Z$ the points of
$\e(\frac12+\ZZ)^2$, and for $w\in\RR^2$ we write $|w|_1:=|w_1|+|w_2|$, so that
$|w|\le|w|_1\le\sqrt2\,|w|$. The network $\e\Z$ contains the horizontal and the
vertical line through each of its vertices; every point of $\RR^2$ has a vertex
at $|\cdot|_1$-distance at most $\e$; and, for every vertex $\beta$, the
reflections $\xi\mapsto(2\beta_1-\xi_1,\xi_2)$ and
$\xi\mapsto(\xi_1,2\beta_2-\xi_2)$ map $\e\Z$ onto itself. Every path in $\e\Z$
from a vertex $\beta$ to a vertex $\beta'$ has length at least
$|\beta'-\beta|_1$. We construct paths that attain this length and stay at
distance of order $\e$ from the segment $[\beta,\beta']$, with no condition on
its direction. They are used with $\e=1$ in Section~\ref{section p<2}, for the
continuity of $f^p_{\rm hom}$ and for the construction of laminates, and in
Section~\ref{section p>2} for the recovery sequences.

\begin{lemma}\label{lem:staircase}
Let $\beta$ and $\beta'$ be vertices of $\e\Z$, let $a:=\beta'-\beta$, and let
$L>0$. There exists a piecewise-affine map $v\colon[0,L]\to\e\Z$, with monotone
components, $v(0)=\beta$ and $v(L)=\beta'$, such that:

\smallskip
{\rm(i)} for almost every $t\in(0,L)$ either
$v'(t)=\frac{|a|_1}{L}\,{\rm sign}(a_1)\,e_1$ or
$v'(t)=\frac{|a|_1}{L}\,{\rm sign}(a_2)\,e_2$, and, if $a\neq0$, the first
alternative holds on a set of measure $\frac{|a_1|}{|a|_1}L$; in particular,
$$
|v(s)-v(t)|\le|v(s)-v(t)|_1=\frac{|a|_1}{L}\,|s-t|\qquad\hbox{for all }s,t\in[0,L];
$$

\smallskip
{\rm(ii)} $\displaystyle\Big|v(t)-\beta-\frac tL\,a\Big|\le\sqrt2\,\e$ for all $t\in[0,L]$;

\smallskip
{\rm(iii)} $\displaystyle\int_0^L|v'(t)|^p\,\dd t=\frac{|a|_1^p}{L^{p-1}}$ for every $p\ge1$.
\end{lemma}

\begin{proof}
If $a_1a_2=0$, the segment $[\beta,\beta']$ lies on one of the lines of $\e\Z$
through $\beta$, and $v(t):=\beta+\frac tL\,a$ has all the properties, with $0$
in place of $\sqrt2\,\e$ in (ii). Otherwise, after one or two reflections about
$\beta$, we may suppose that $a_1,a_2>0$, and we write $a=\e(m_1,m_2)$, with
$m_1,m_2\in\NN$.

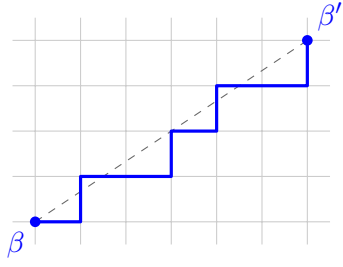
\begin{figure}[ht!]
\centering
\begin{tikzpicture}[x=0.6cm,y=0.6cm]
  \draw[gray!40,thin,step=1] (-0.5,-0.5) grid (6.5,4.5);%
  \draw[black!60,dashed] (0,0) -- (6,4);
  \draw[blue,very thick,line join=round]
    (0,0) -- (1,0) -- (1,1) -- (3,1) -- (3,2) -- (4,2) -- (4,3) -- (6,3) -- (6,4);
  \fill[blue] (0,0) circle (2pt) node[below left] {$\beta$};
  \fill[blue] (6,4) circle (2pt) node[above right] {$\beta'$};
\end{tikzpicture}
\caption{The staircase of Lemma~\ref{lem:staircase} for
$\beta'-\beta=\e(6,4)$. The grid lines, at mutual distance $\e$, belong to
$\e\Z$; the dashed line is the segment $[\beta,\beta']$.}\label{fig:staircase}
\end{figure}

For $k\in\{0,1,\dots,m_1\}$, let $x_k:=k$ and let $y_k$ be the integer closest to
$km_2/m_1$ (the smaller one in the case of a tie), so that $y_0=0$, $y_{m_1}=m_2$,
$k\mapsto y_k$ is nondecreasing, and $|m_1y_k-m_2x_k|\le m_1/2$. Let $\gamma$
be the polygonal path through the points
$$
(x_0,y_0),\ (x_1,y_0),\ (x_1,y_1),\ (x_2,y_1),\ \dots,\
(x_{m_1},y_{m_1-1}),\ (x_{m_1},y_{m_1})=(m_1,m_2),
$$
which alternates horizontal steps of length $1$ with vertical steps of length
$y_k-y_{k-1}\ge0$ (see Fig.~\ref{fig:staircase}). The path $\beta+\e\gamma$
lies on the horizontal and vertical lines through the vertices
$\beta+\e\ZZ^2$, hence in $\e\Z$. Since the components of $\gamma$ are
nondecreasing, its length is $m_1+m_2$; we parameterize it by arclength,
$\sigma\mapsto\gamma(\sigma)$, $\sigma\in[0,m_1+m_2]$, and set
$$
v(t):=\beta+\e\,\gamma\Big(\frac{m_1+m_2}{L}\,t\Big),\qquad t\in[0,L].
$$
Then $v'=\frac{|a|_1}{L}\,e_1$ on the horizontal steps and
$v'=\frac{|a|_1}{L}\,e_2$ on the vertical ones, and integrating the first
component over $[0,L]$ shows that the horizontal steps take a time of measure
$\frac{a_1}{|a|_1}L$. Since the components of $v$ are monotone and, at almost
every time, only one of them changes, $|v(t)-v(s)|_1=\int_s^t|v'|\,\dd\tau
=\frac{|a|_1}{L}(t-s)$ for $s\le t$. This proves (i), and (iii) follows from
$|v'|=\frac{|a|_1}{L}$ almost everywhere.

As for (ii), let $t\in[0,L]$ and $(x,y):=\gamma((m_1+m_2)t/L)$. Since $\gamma$
has nondecreasing components, $x+y=(m_1+m_2)t/L$, so that
$$
\beta+\frac tL\,a-v(t)=\e\Big(\frac{x+y}{m_1+m_2}\,(m_1,m_2)-(x,y)\Big)
=\frac{\e\,(m_1y-m_2x)}{m_1+m_2}\,(1,-1).
$$
On each step of $\gamma$ the quantity $m_1y-m_2x$ is an affine function of $t$,
so that it suffices to estimate it at the points listed above. At $(x_k,y_k)$
we have $|m_1y_k-m_2x_k|\le m_1/2$, and at $(x_k,y_{k-1})$, $k\ge1$, we have
$|m_1y_{k-1}-m_2x_k|=|m_1y_{k-1}-m_2x_{k-1}-m_2|\le m_1/2+m_2$. Hence
$|m_1y-m_2x|\le m_1+m_2$, and (ii) follows.
\end{proof}

\subsection{Homogenized energy functions}\label{sec:en-fu} We now introduce the homogenized energy density that will allow us to describe the $\Gamma$-limit. It is defined differently on matrices of rank $2$ and on matrices of rank $1$.

\begin{definition}\label{formula f_hom with A non singular}
For all $A\in\RR^{2\times 2}$ with $\emph{rank}(A)=2$, we define 
\begin{eqnarray}\nonumber
f^p_{\rm hom}(A)&=&\liminf_{T\to+\infty}\frac{|\det A|}{T^2}\inf\Bigl\{\int_{TQ_A}|\nabla u|^p\,\dd x: u\in \WW^{1,p}(TQ_A,
\Z),\ \\
&& \label{full rank formula}\qquad\qquad \qquad\qquad 
u= u_A \hbox{ \rm in a neighbourhood of } \partial TQ_A\Bigr\},
\end{eqnarray}
where $u_A$ is as in Definition {\rm\ref{projlin}} and $Q_A=A^{-1}(Q)$ (see Fig.~{\rm\ref{FigQA}}).

\begin{figure}[ht!]
\centering
\begin{tikzpicture}[x=0.5cm,y=0.5cm]
  \def\ccTT{5}%
   \def\ccpax{0.85}\def\ccpay{0.6}\def\ccpbx{-0.4}\def\ccpby{0.85}%
  \pgfmathtruncatemacro\cch{(\ccTT-1)/2}%
  \pgfmathsetmacro\cclow{max(\ccTT/2,\ccTT/2*(abs(\ccpay)+abs(\ccpby)))}%
  \begin{scope}[shift={(-6.4,0)}]
    \fill[ccviolet!15] (-\ccTT/2,-\ccTT/2) rectangle (\ccTT/2,\ccTT/2);
    \fill[ccviolet!45] (-.5,-.5) rectangle (.5,.5);
    \foreach \i in {-\cch,...,\cch}{%
      \draw[black!35,line width=.2pt] (\i-.5,-\ccTT/2) -- (\i-.5,\ccTT/2) (-\ccTT/2,\i-.5) -- (\ccTT/2,\i-.5);}
    \draw[ccviolet,line width=.8pt] (-\ccTT/2,-\ccTT/2) rectangle (\ccTT/2,\ccTT/2);
    \node[font=\scriptsize] at (0,0) {$Q$};
  \end{scope}
  \begin{scope}[shift={(6.4,0)},cm={\ccpax,\ccpay,\ccpbx,\ccpby,(0,0)}]
    \fill[ccteal!15] (-\ccTT/2,-\ccTT/2) rectangle (\ccTT/2,\ccTT/2);
    \fill[ccteal!45] (-.5,-.5) rectangle (.5,.5);
    \foreach \i in {-\cch,...,\cch}{%
      \draw[black!35,line width=.2pt] (\i-.5,-\ccTT/2) -- (\i-.5,\ccTT/2) (-\ccTT/2,\i-.5) -- (\ccTT/2,\i-.5);}
    \draw[ccteal,line width=.8pt] (-\ccTT/2,-\ccTT/2) rectangle (\ccTT/2,\ccTT/2);
    \node[font=\scriptsize] at (0,0) {$Q_A$};
  \end{scope}
  \draw[->,line width=.6pt] (-3.3,0) -- (3.1,0) node[midway,above] {$A^{-1}$};
  \node[below] at (-6.4,-\cclow) {$TQ$};
  \node[below] at (6.4,-\cclow) {$TQ_A=A^{-1}(TQ)$};
\end{tikzpicture}
\caption{The domain $TQ_A=A^{-1}(TQ)$ of the test functions in the definition of $f^p_{\rm hom}$.}\label{FigQA}
\end{figure}
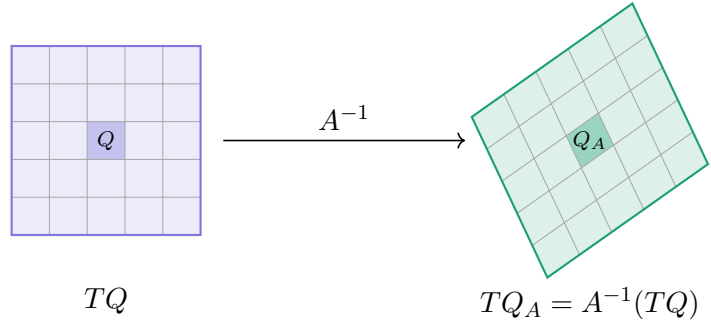

If {\rm rank}$(A)\le 1$ then we set
$$
f^p_{\rm hom}(A)=(|a_1|+|a_2|)^p,
$$
where $A=a \otimes \nu$, with $\nu\in S^1$ and $a=(a_1,a_2)$. Equivalently, for any rank-one matrix $A=a\otimes \xi$ we set
$$
f^p_{\rm hom}(A)=(|A_1|+|A_2|)^p=(|a_1|+|a_2|)^p|\xi|^p,
$$
where $A_i=a_i \xi$ are the rows of $A$. 
\end{definition}

Although this property is not immediate, $f^p_{\rm hom}$ is $p$-homogeneous, as proved in Lem\-ma~\ref{lem:homog} below.

\smallskip
In the following, when necessary to shorten the notation,  for fixed $T>0$ and $A$ nonsingular, $\mathcal{A}_T(A)$ will denote the set of competitors in the infimum problem in \eqref{full rank formula}.

\begin{remark}[Estimate of $f^p_{\rm hom}$ on the identity map $A=I$]
\rm
Take $u^T(x)=P_0(x)$ to test $f^p_{\rm hom}(I)$, for fixed $T>0$. Then
$$
\nabla u^T(x)= (\nabla P_0)(x) = \frac{1}{2}(\nabla \psi)(x-z), \quad x=  z+ \left[-\tfrac{1}{2},\tfrac{1}{2}\right]^2 \setminus \{0\}.
$$
Suppose for simplicity that $T$ is odd, so that we cover an exact number of cells without some partial cells. Then, since $A=I$ gives $TQ_A=TQ$,
\begin{eqnarray*}
f^p_{\rm hom}(I)&\leq& \frac{1}{T^2}\int_{TQ} |\nabla u_\e|^p \dd x = \frac{1}{T^2}\sum_{i=1}^{T^{2}} \int_{Q_z}\frac{1}{2^p}|\nabla \psi(x-z)|^p \dd x \\ &=&\frac{1}{T^2}\sum_{i=1}^{T^{2}} \int_{Q}\frac{1}{2^p}|\nabla \psi(y)|^p \dd y= \frac{1}{2^p} \int_{Q} |\nabla \psi(y)|^p \dd y.
\end{eqnarray*}
By symmetry, we can divide $Q$ into the four triangles each spanned by an angle of $\pi/2$, and compute the integral only for one of them. Let $\Delta :=\{ 0< y_1 < 1/2, |y_2|< y_1  \}$ be the first triangle, then 
$$
\psi(y_1,y_2)=\left(1, \tfrac{y_2}{y_1}\right), \quad \nabla \psi(y_1,y_2)=\begin{pmatrix}
 0 & 0\\  -\frac{y_2}{y_1^2}& \frac{1}{y_1}    
\end{pmatrix}, 
$$
(note $\nabla \psi \in K$ has $\det \nabla \psi=0$, see Lemma \ref{nabla u_e has only one row}), and
$$
|\nabla \psi|=\frac{|y|}{y_1^2}. 
$$
In polar coordinates $|\nabla \psi|(r,\theta)=1/(r \cos^2\theta ) $, with domain $|\theta| < \pi/4$ and $0< r < \frac{1}{2 \cos \theta}$. In this domain $\cos^2\theta > 1/2$, so that $|\nabla \psi|\leq 1/r \in \LL^p_{\rm loc}$ for $p<2$. We can now compute
$$
\frac{1}{2^p} \int_{Q} |\nabla \psi(y)|^p \dd y = \frac{4}{2^p} \int_{\Delta} |\nabla \psi(y)|^p \dd y
$$
\begin{eqnarray*}
\int_{\Delta} |\nabla \psi(y)|^p \dd y = \int_{-\pi/4}^{\pi/4}\frac{1}{\cos^{2p} \theta} \int_{0}^{\tfrac{1}{2\cos\theta}} r^{1-p} \, \dd r \, \dd\theta   = \frac{1}{(2-p)2^{2-p}}\int_{-\pi/4}^{\pi/4}  \cos^{-(p+2)}\theta  \, \dd\theta.
\end{eqnarray*} 
As $1<p<2$ and in $|\theta|< \pi/4$ it holds $1/\sqrt{2}< \cos \theta \leq 1$, the integral above is strictly increasing in $p$ and we have the bounds (attained at the endpoints $p=1$ and $p=2$) 
$$
\sqrt{2}+\ln(1+\sqrt{2})< \int_{-\pi/4}^{\pi/4}  \cos^{-(p+2)}\theta  \, \dd\theta < \frac{8}{3}, \quad 1<p<2.
$$
Therefore
\begin{equation}\label{f on the ID upper bound}
f^p_{\rm hom}(I)\leq \frac{4}{2^p} \int_{\Delta} |\nabla \psi(y)|^p \dd y \leq \frac{8}{3(2-p)}.
\end{equation}
This estimate will be useful for deriving an upper bound of $f^p_{\rm hom}$ in Proposition \ref{upper bound of f hom with N(A)}. 
\end{remark}

The following lemma states that the lower limit in \eqref{full rank formula} is indeed a limit. In the rest of the section we use the shorthand $u= u_A \sim \partial TQ$ to mean $u= u_A$ in a neighbourhood of $\partial TQ$.
\begin{lemma}
The liminf in definition \eqref{full rank formula} is actually a limit; that is,
\begin{eqnarray*}
f^p_{\rm hom}(A)=\lim_{T\to+\infty}\frac{|\det A|}{T^2}\inf\Bigl\{\int_{TQ_A}|\nabla u|^p\,\dd x: u\in \WW^{1,p}(TQ_A,
\Z),\
u= u_A \sim \partial TQ_A\Bigr\}.
\end{eqnarray*}
\end{lemma}

\begin{proof}
Fix $A\in \RR^{2\times 2}$ nonsingular and define
$$
g_T:= \frac{|\det A|}{T^2} \inf\Bigl\{\int_{TQ_A}|\nabla u|^p\,\dd x: u\in \WW^{1,p}(TQ_A,
\Z),\quad 
u= u_A \sim \partial TQ_A\Bigr\}.
$$
Let $u_T$ satisfy the boundary conditions of the infimum and be such that
$$
\frac{|\det A|}{T^2}\int_{TQ_A}|\nabla u_T|^p\,\dd x \leq g_T + \frac{1}{T}.
$$
Set $N_T:=\lfloor T\rfloor+1$, $w_T:=N_T-T\in(0,1]$, and, for $S\ge T$,
$m:=\bigl\lfloor\frac{S-T}{2N_T}\bigr\rfloor$, and $\mathbf{I}:=\{-m,\dots,m\}^2$. Then $\#\mathbf{I}=(2m+1)^2$, and the cells
$A^{-1}(N_T\mathbf{i})+TQ_A$, $\mathbf{i}\in\mathbf{I}$, are pairwise disjoint and contained in $SQ_A$. Define
$$
u_S(x)=\begin{cases}
u_T\bigl(x-A^{-1}(N_T\mathbf{i})\bigr)+N_T\mathbf{i} & \hbox{if } x\in A^{-1}(N_T\mathbf{i})+TQ_A,\ \mathbf{i}\in\mathbf{I}\\
u_A(x) & \hbox{otherwise in } SQ_A.
\end{cases}
$$
Since $P_0(y+z)=P_0(y)+z$ for $z\in\ZZ^2$, we have $u_A(x-A^{-1}(N_T\mathbf{i}))+N_T\mathbf{i}=u_A(x)$,
so $u_S=u_A$ in a neighbourhood of the boundary of each cell; as $\Z+\ZZ^2=\Z$, it follows that
$u_S\in \WW^{1,p}(SQ_A,\Z)$ and $u_S=u_A$ near $\partial SQ_A$.

We now estimate the energy of $u_A$ on $L_{S,T}:=SQ_A\setminus\bigcup_{\mathbf{i}\in\mathbf{I}}(A^{-1}(N_T\mathbf{i})+TQ_A)$.
Since $P_0(y+z)=P_0(y)+z$, the function $|(\nabla P_0)(y)A|^p$ is $\ZZ^2$-periodic, so that
$$
\int_{z+Q}|(\nabla P_0)(y)A|^p\,\dd y=C_A\qquad\hbox{for every } z\in\ZZ^2,
$$
with $C_A$ as in Lemma \ref{rescaled P_0 satisfies  C|U| bound}.
In the variable $y=Ax$ we have $A(L_{S,T})\subset(G\times I_S)\cup(I_S\times G)$, where $I_S:=(-\frac S2,\frac S2)$ and
$$
G:=I_S\setminus\bigcup_{i=-m}^{m}\Bigl(N_Ti-\frac T2,\,N_Ti+\frac T2\Bigr)
$$
is the union of $2m$ intervals of length $w_T\le1$ and of two intervals of length less than $N_T$. Since an interval of length $L$ is covered by at most $\lceil L\rceil+1$ intervals $[l-\frac12,l+\frac12]$, $l\in\ZZ$
(it may straddle their endpoints), $G$ is covered by at most $2\cdot 2m+2(N_T+1)=4m+2N_T+2$ of them and $I_S$ by
at most $\lceil S\rceil+1\le S+2$,
so that $A(L_{S,T})$ is covered by at most $2(4m+2N_T+2)(S+2)$ squares $z+Q$, $z\in\ZZ^2$. Therefore
\begin{eqnarray*}
\frac{|\det A|}{S^2}\int_{L_{S,T}}|\nabla u_A|^p\,\dd x&=&\frac{1}{S^2}\int_{A(L_{S,T})}|(\nabla P_0)(y)A|^p\,\dd y
\\ &\le& r(S,T):=2\,C_A\,\frac{(4m+2N_T+2)(S+2)}{S^2}.
\end{eqnarray*}
We can now estimate $g_S$ using $u_S$ as a test function; since the constants $N_T\mathbf{i}$ do not affect the gradient,
\begin{eqnarray*}
   g_S &\leq& \frac{|\det A|}{S^2} \int_{SQ_A} |\nabla u_S|^p\,\dd x \\
   &=& \frac{|\det A|}{S^2} \left( \sum_{\mathbf{i} \in \mathbf{I}} \int_{A^{-1}(N_T\mathbf{i}) + TQ_A} |\nabla u_T(x-A^{-1}(N_T\mathbf{i}))|^p \, \dd x + \int_{L_{S,T}} |\nabla u_A|^p\, \dd x \right) \\
   &\leq& \frac{|\det A|}{S^2}\, (2m+1)^2 \int_{TQ_A} |\nabla u_T|^p \, \dd x + r(S,T) \\
   &\leq& \frac{T^2(2m+1)^2}{S^2} \left(g_T+\frac{1}{T} \right) + r(S,T) \\
   &\leq& g_T +\frac{1}{T}  + r(S,T),
\end{eqnarray*}
where in the last inequality we used $T(2m+1)\le\frac{T}{N_T}(S-T)+T\le S$. Since $m/S\to 1/(2N_T)$ as $S\to+\infty$,
$$
\limsup_{S\to+\infty} g_S\le g_T+\frac{1}{T}+\frac{4\,C_A}{N_T}\qquad\hbox{for every } T>0 .
$$
Taking the lower limit as $T\to+\infty$, and recalling that $N_T>T$, we obtain
$\limsup\limits_{S\to+\infty} g_S\le\liminf\limits_{T\to+\infty} g_T$, which proves the claim.
\end{proof}
\begin{remark}\label{why the tilted lattice cell formula}
\emph{We use the tilted squares $TQ_A=A^{-1}(TQ)$ in the cell formula \eqref{full rank formula} because of the upper
bound in Proposition \ref{limsup for piecewise affine functions}, where we need a recovery sequence whose gradient
is periodic, in order to apply the Riemann--Lebesgue lemma. The reason is the boundary datum $u_A(x)=P_0(Ax)$,
which projects the affine function $Ax$ on the $\ZZ^2$-periodic grid $\Z$ \emph{in the target space}: since
$P_0(y+z)=P_0(y)+z$ for $z\in\ZZ^2$, we have $u_A(x+A^{-1}z)=u_A(x)+z$, so that $u_A(x)-Ax$ is periodic with
respect to the tilted lattice $A^{-1}\ZZ^2$ \emph{in the domain}, whose cells $A^{-1}(z+Q)$ are mapped by $A$ onto
the squares $z+Q$ of the grid. Hence, for $T\in\mathbb{N}$, a competitor on $TQ_A$ that equals $u_A$ near
$\partial TQ_A$ can be translated by the vectors $A^{-1}(T\mathbf{k})$, $\mathbf{k}\in\ZZ^2$, adding $T\mathbf{k}$ to its
values, and the copies glue into a function whose gradient is periodic with periodicity cell $TQ_A$. By
Lemma \ref{cubes and losanghe lemma}, the squares $TQ$ give the same value $f^p_{\rm hom}(A)$, so the tilted
squares are only a matter of convenience.}
\end{remark}
\begin{lemma}\label{cubes and losanghe lemma}
Let \emph{rank}$(A)=2$ and define
$$
f^{\#,p}_{\rm{hom}}(A):=\liminf_{T\to+\infty}\frac{1}{T^2}\inf\Bigl\{\int_{TQ}|\nabla u|^p\,\dd x: u\in \WW^{1,p}(TQ,
\Z),\
u= u_A \sim \partial TQ\Bigr\}.
$$
Then
$f^{\#,p}_{\rm{hom}}(A) = f^{p}_{\rm{hom}}(A)$.
\end{lemma}

\begin{proof}
Denote
$$
g^{\#}_T(A)=\frac{1}{T^2}\inf\Bigl\{\int_{TQ}|\nabla u|^p\,\dd x: u\in \WW^{1,p}(TQ,
\Z),\
u= u_A \sim \partial TQ\Bigr\},
$$
$$
g_T(A)=\frac{|\det A|}{T^2}\inf\Bigl\{\int_{TQ_A}|\nabla u|^p\,\dd x: u\in \WW^{1,p}(TQ_A,
\Z),\
u= u_A \sim \partial TQ_A\Bigr\}.
$$
We now prove that for $T\in\mathbb{N}$ and $S\gg T$
$$
g^{\#}_S(A) \leq g_T(A) +o_T(1) + C \frac{T}{S},
$$
where $C$ depends on $A$ and $p$. Since the liminf defining $f^{p}_{\rm{hom}}(A)$ is a limit, this will give, taking
the liminf first for $S\to \infty$ and then for $T\to \infty$ along the integers, the first inequality
$f^{\#,p}_{\rm{hom}}(A) \leq f^{p}_{\rm{hom}}(A)$.

\begin{figure}[ht!]
\centering
\begin{tikzpicture}[scale=0.9]
  \ccsetup
  \pgfmathsetmacro\ccsqX{\ccT/2+(abs(\ccax)+abs(\ccbx))/2}%
  \pgfmathsetmacro\ccsqY{\ccT/2+(abs(\ccay)+abs(\ccby))/2}%
  \pgfmathsetmacro\ccsqdet{abs(\ccax*\ccby-\ccay*\ccbx)}%
  \pgfmathsetmacro\ccsqI{0.5+\ccT*(abs(\ccbx)+abs(\ccby))/(2*\ccsqdet)}%
  \pgfmathsetmacro\ccsqJ{0.5+\ccT*(abs(\ccax)+abs(\ccay))/(2*\ccsqdet)}%
  \def\ccmeets#1#2{%
    \edef\ccsqi{\the\numexpr#1\relax}\edef\ccsqj{\the\numexpr#2\relax}%
    \edef\ccsqx{\the\dimexpr\ccax pt*(\ccsqi)+\ccbx pt*(\ccsqj)\relax}%
    \edef\ccsqy{\the\dimexpr\ccay pt*(\ccsqi)+\ccby pt*(\ccsqj)\relax}%
    \def\ccc{1}%
    \ifdim\ccsqx<\ccsqX pt\else\def\ccc{0}\fi
    \ifdim\ccsqx>-\ccsqX pt\else\def\ccc{0}\fi
    \ifdim\ccsqy<\ccsqY pt\else\def\ccc{0}\fi
    \ifdim\ccsqy>-\ccsqY pt\else\def\ccc{0}\fi
    \ifdim\ccsqi pt<\ccsqI pt\else\def\ccc{0}\fi
    \ifdim\ccsqi pt>-\ccsqI pt\else\def\ccc{0}\fi
    \ifdim\ccsqj pt<\ccsqJ pt\else\def\ccc{0}\fi
    \ifdim\ccsqj pt>-\ccsqJ pt\else\def\ccc{0}\fi}%
  \pgfmathsetmacro\ccW{max(max(abs(\cclax)+abs(\cclbx),abs(\cclay)+abs(\cclby))/2+\ccd/2,\ccT/2+\ccd)+0.3}%
  \pgfmathtruncatemacro\ccK{ceil(\ccAn*\ccW)+1}\pgfmathtruncatemacro\ccKm{-\ccK}%
  \pgfmathsetmacro\ccu{7.5/(2*\ccW)}%
  \pgfmathsetmacro\cclim{\ccT/2+0.12}%
  \pgfmathsetmacro\cclimx{(\ccT-abs(\ccax)-abs(\ccbx))/2}%
  \pgfmathsetmacro\cclimy{(\ccT-abs(\ccay)-abs(\ccby))/2}%
  \begin{scope}[x=\ccu cm,y=\ccu cm]
    \clip (-\ccW,-\ccW) rectangle (\ccW,\ccW);
    \fill[ccviolet!15] (-\ccT/2,-\ccT/2) rectangle (\ccT/2,\ccT/2);
    \begin{scope}[cclattice]
      \foreach \i in {\ccKm,...,\ccK}{\foreach \j in {\ccKm,...,\ccK}{%
        \edef\ccx{\the\dimexpr\ccax pt*(\i)+\ccbx pt*(\j)\relax}%
        \edef\ccy{\the\dimexpr\ccay pt*(\i)+\ccby pt*(\j)\relax}%
        \edef\cco{\ifdim\ccx>\cclimx pt 0\else\ifdim\ccx<-\cclimx pt 0\else
          \ifdim\ccy>\cclimy pt 0\else\ifdim\ccy<-\cclimy pt 0\else 1\fi\fi\fi\fi}%
        \ifnum\cco=1 \fill[ccteal!35,draw=ccteal!35,line width=.3pt] \cccell{\i}{\j};\fi}}
      \foreach \j in {\ccKm,...,\ccK}{%
        \draw[black!35,line width=.2pt] (\j+.5,\ccKm-1) -- (\j+.5,\ccK+1)
                                        (\ccKm-1,\j+.5) -- (\ccK+1,\j+.5);}
      \foreach \i in {\ccKm,...,\ccK}{\foreach \j in {\ccKm,...,\ccK}{%
        \edef\ccx{\the\dimexpr\ccax pt*(\i)+\ccbx pt*(\j)\relax}%
        \edef\ccy{\the\dimexpr\ccay pt*(\i)+\ccby pt*(\j)\relax}%
        \edef\cco{\ifdim\ccx>\cclim pt 1\else\ifdim\ccx<-\cclim pt 1\else
          \ifdim\ccy>\cclim pt 1\else\ifdim\ccy<-\cclim pt 1\else 0\fi\fi\fi\fi}%
        \ifnum\cco=1
          \ccmeets{\i}{\j}%
          \ifnum\ccc=1 \def\cccol{ccteal!70!black}\else\def\cccol{black!45}\fi
          \node[circle,inner sep=0pt,minimum size=1.5pt,fill=\cccol] at (\i,\j) {};
        \fi}}

      \foreach \i in {\ccKm,...,\ccK}{\foreach \j in {\ccKm,...,\ccK}{%
        \ccmeets{\i}{\j}\let\cca\ccc
        \ccmeets{\i+1}{\j}\let\ccr\ccc
        \ccmeets{\i}{\j+1}\let\cct\ccc
        \ifnum\cca=\ccr\space\else\draw[ccteal,line width=.8pt,line cap=round] (\i+.5,\j-.5) -- (\i+.5,\j+.5);\fi
        \ifnum\cca=\cct\space\else\draw[ccteal,line width=.8pt,line cap=round] (\i-.5,\j+.5) -- (\i+.5,\j+.5);\fi}}
    \end{scope}
    \draw[ccviolet,line width=.8pt] (-\ccT/2,-\ccT/2) rectangle (\ccT/2,\ccT/2);
    \node[fill=white,inner sep=1.5pt] at (0,0) {$SQ$};
  \end{scope}
  \begin{scope}[shift={(-3.2,-4.25)}]
    \filldraw[fill=ccviolet!15,draw=ccviolet] (0,-.12) rectangle ++(.4,.24);
    \node[right] at (.45,0) {$SQ$};
    \begin{scope}[shift={(1.75,0)},x=.36cm,y=.36cm,cclattice]
      \filldraw[fill=ccteal!35,draw=ccteal] \cccell{0}{0};
    \end{scope}
    \node[right] at (2.0,0) {$TA^{-1}(Q+k)\subset SQ$};
    \draw[ccteal,line width=.8pt,line cap=round] (0,-.55) -- ++(.4,0);
    \node[right] at (.45,-.55) {union of the cells meeting $SQ$};
  \end{scope}
\end{tikzpicture}
\caption{From lozenges to squares}\label{FigA}
\end{figure}
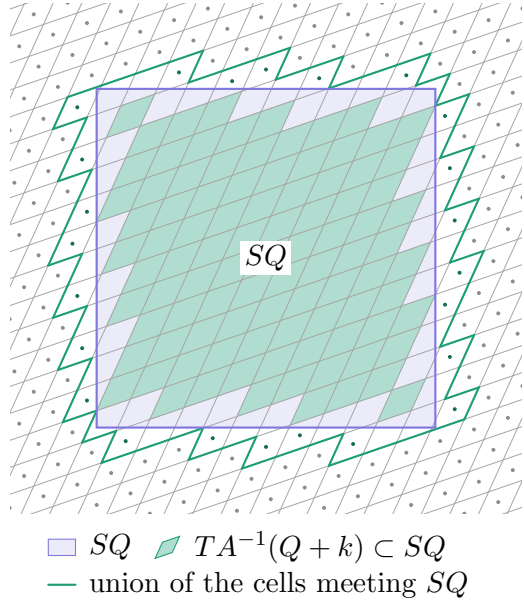

Let $u_T$ be a quasi-minimum for $g_T$
$$
\frac{|\det A|}{T^2}\int_{TQ_A}|\nabla u_T|^p\,\dd x \leq g_T + \frac{1}{T}.
$$
Referring to Fig.~\ref{FigA}, we want to fill in the big square of size $S$ with parallelograms $TQ_A$ and estimate
the energy of $u_A$ on the remaining set, where some parallelograms do not fit completely in the square; the thickness
of this boundary layer is at most $\delta_T:=T\cdot \textrm{diam}(Q_A)\leq T \sqrt{2}|A^{-1}|$, since $\sqrt{2}=\textrm{diam}(Q)$.

Let $\mathbf{I}$ be the set of indices $\mathbf{i}=(i_1,i_2) \in T\ZZ^2$ for which $A^{-1}(\mathbf{i}) + TQ_A \subset SQ$
(that is, the $T$-parallelograms are well contained). We define $u_S$ as
$$
u_S(x)=\begin{cases}\displaystyle u_T\big(x-A^{-1}(\mathbf{i})\big)+\mathbf{i} &\hbox{ if $x \in A^{-1}(\mathbf{i}) + TQ_A$ for some $\mathbf{i}\in\mathbf{I}$}\cr
u_A(x)=P_0(Ax) &\hbox{ if $x \in L_{S,T}:=SQ \setminus \bigcup_{\mathbf{i} \in \mathbf{I}}(A^{-1}(\mathbf{i}) + TQ_A)$.}
\end{cases}
$$
Since $\mathbf{i}\in T\ZZ^2\subset\ZZ^2$ and $P_0(y+z)=P_0(y)+z$ for $z\in\ZZ^2$, we have
$u_A(x-A^{-1}(\mathbf{i}))+\mathbf{i}=u_A(x)$, so that $u_S=u_A$ in a neighbourhood of the boundary of each
parallelogram; as $\Z+\ZZ^2=\Z$, it follows that $u_S\in \WW^{1,p}(SQ,\Z)$ and $u_S=u_A$ near $\partial SQ$.

Since $\nabla u_A$ is not bounded, we estimate its energy on $L_{S,T}$ by covering $L_{S,T}$ with whole parallelograms. The parallelograms
$A^{-1}(\mathbf{i})+TQ_A$, $\mathbf{i}\in T\ZZ^2$, are pairwise disjoint, so that
$L_{S,T}$ is contained, up to a null set, in the union of those that meet $SQ$ without being contained in it. Each of
these meets $\partial SQ$; hence it lies in $\{x:{\rm dist}(x,\partial SQ)\le\delta_T\}$, whose measure is at most
$8S\delta_T$ for $S\ge2\delta_T$ (as it is made by two concentric squares); since $|TQ_A|=T^2/|\det A|$, their number is at most
$$
\frac{8S\delta_T}{|TQ_A|}\leq 8\sqrt2\,|\det A|\,|A^{-1}|\,\frac{S}{T}.
$$
Moreover, since $T\in\mathbb{N}$, each square $\mathbf{i}+TQ$ is, up to a null set, the union of $T^2$ unit squares,
on each of which the $\ZZ^2$-periodic function $|(\nabla P_0)(y)A|^p$ has integral $C_A$, with $C_A$ as in
Lemma \ref{rescaled P_0 satisfies  C|U| bound}. Hence
$$
\int_{A^{-1}(\mathbf{i})+TQ_A}|\nabla u_A|^p\,\dd x=\frac{1}{|\det A|}\int_{\mathbf{i}+TQ}|(\nabla P_0)(y)A|^p\,\dd y
=\frac{T^2}{|\det A|}\,C_A ,
$$
and therefore
$$
\int_{L_{S,T}}|\nabla u_A|^p\,\dd x\leq 8\sqrt2\,C_A\,|A^{-1}|\,S\,T .
$$
We can also estimate the cardinality of $\mathbf{I}$ as
$\#\mathbf{I}\leq \frac{S^2}{|TQ_A|}= |\det A| \frac{S^2}{T^2}$.
Then
\begin{eqnarray*}
   g^{\#}_S &\leq& \frac{1}{S^2} \int_{SQ} |\nabla u_S|^p\,\dd x \\
   &=& \frac{1}{S^2} \left( \sum_{\mathbf{i} \in \mathbf{I}} \int_{A^{-1}(\mathbf{i}) + TQ_A} |\nabla u_T(x-A^{-1}(\mathbf{i}))|^p \, \dd x + \int_{L_{S,T}} |\nabla u_A|^p\, \dd x \right) \\
   &\leq& \frac{1}{S^2} \left( \#\mathbf{I} \int_{TQ_A} |\nabla u_T|^p \, \dd x + 8\sqrt2\,C_A\,|A^{-1}|\,S\,T  \right) \\
   &\leq & \frac{|\det A|}{T^2} \int_{TQ_A} |\nabla u_T|^p \, \dd x + 8\sqrt2\,C_A\,|A^{-1}|\,\frac{T}{S}\\
   &\leq & g_T+\frac{1}{T}  + 8\sqrt2\,C_A\,|A^{-1}|\,\frac{T}{S},
\end{eqnarray*}
which gives the claim with $C=8\sqrt2\,C_A|A^{-1}|$.

For the opposite inequality we reverse the roles of squares and parallelograms, with one difference: a translated
copy $u^\#_T(x-\tau)+A\tau$ of a quasi-minimum $u^\#_T$ for $g^\#_T$ agrees with $u_A$ near the boundary of
$\tau+TQ$ only if $\tau\in A^{-1}\ZZ^2$, and the translates of $TQ$ by the lattice $A^{-1}\ZZ^2$ do not tile the plane.
We therefore set $d:=\textrm{diam}(Q_A)$ and, for every $\mathbf{j}\in\ZZ^2$, choose $\mathbf{k}_{\mathbf{j}}\in\ZZ^2$ with
$|A^{-1}\mathbf{k}_{\mathbf{j}}-(T+d)\mathbf{j}|\le d/2$, which is possible since the cells $A^{-1}(\mathbf{k}+Q)$ cover the
plane and have diameter $d$. The squares $A^{-1}\mathbf{k}_{\mathbf{j}}+TQ$ are then contained in $(T+d)(\mathbf{j}+Q)$,
hence pairwise disjoint, and we place in them the copies $u^\#_T(x-A^{-1}\mathbf{k}_{\mathbf{j}})+\mathbf{k}_{\mathbf{j}}$ for all
$\mathbf{j}$ with $(T+d)(\mathbf{j}+Q)\subset SQ_A$, setting $u_S=u_A$ elsewhere in $SQ_A$. The frames
$(T+d)(\mathbf{j}+Q)\setminus(A^{-1}\mathbf{k}_{\mathbf{j}}+TQ)$ have width at most $d$, so covering them by cells
$A^{-1}(z+Q)$, $z\in\ZZ^2$, and using periodicity as above, their total energy is at most $C\,d\,S^2/T$, while the
boundary layer of $SQ_A$ is treated as before. This gives
$$
g_S\le g^\#_T+\frac1T+C\Bigl(\frac dT+\frac TS\Bigr),
$$
with $C$ depending on $A$ and $p$; letting first $S\to+\infty$ and then taking the lower limit as $T\to+\infty$, we obtain
$f^{p}_{\rm hom}(A)\le f^{\#,p}_{\rm hom}(A)$.
\end{proof}
Denote for simplicity
$\mathcal R:=\{A\in\RR^{2\times2}:\det A\neq0\}$. For $A\in\mathcal R$, we recall the definition of $f^p_{\rm hom}$ in \eqref{formula f_hom with A non singular} in a shortened notation,
$$
f^p_{\rm hom}(A)=\lim_{T\to\infty}g_T(A),\qquad
g_T(A):=\inf\Big\{\frac{1}{|TQ_A|}\int_{TQ_A}|\nabla v|^p\,\dd x:\ v\in\mathcal A_T(A)\Big\},
$$
where $\mathcal A_T(A)$ is the already defined set of competitors
$v\in \WW^{1,p}(TQ_A,\RR^2)$ with $v(x) \in\mathcal Z$ for almost every $x$ and $v(x)=P_0(Ax)$ on
$\partial TQ_A$. For $\det A=0$, $f^p_{\rm hom}(A)=N(A)^p$, where $N(A)=|A_1|+|A_2|$ with $A_i$ the rows of $A$.

\begin{lemma}[$p$-homogeneity]\label{lem:homog}
For every $A\in\RR^{2\times2}$ and $s\in\RR\setminus\{0\}$,
$$f^p_{\rm hom}(sA)=|s|^p f^p_{\rm hom}(A).
$$
\end{lemma}

\begin{proof}
If $\det A=0$, then
$f^p_{\rm hom}(sA)=N(sA)^p=|s|^pN(A)^p=|s|^pf^p_{\rm hom}(A)$.

Let $\det A\neq0$ and $s>0$. Since $TQ_{sA}=s^{-1}TQ_A$ and dilations
preserve traces, the map $v\mapsto w$, $w(y):=v(y/s)$, is a bijection from
$\mathcal A_T(sA)$ onto $\mathcal A_T(A)$: indeed, $w$ takes values in
$\mathcal Z$, and on $\partial TQ_A$ we have $w(y)=P_0(sA\,y/s)=P_0(Ay)$.
Moreover $\nabla w(y)=s^{-1}\nabla v(y/s)$ and $|TQ_A|=s^2|TQ_{sA}|$, so that
$$
\frac{1}{|TQ_A|}\int_{TQ_A}|\nabla w|^p\,\dd y=s^{-p}\frac{1}{|TQ_{sA}|}\int_{TQ_{sA}}|\nabla v|^p\,\dd x .
$$
Hence $g_T(sA)=s^p\,g_T(A)$ for every $T$, and letting $T\to\infty$ gives the
claim for $s>0$.

For $s=-1$: since $Q=-Q$, we have $TQ_{-A}=TQ_A$. Since
$\mathcal Z=-\mathcal Z$ and $P_0(-y)=-P_0(y)$, the map $v\mapsto-v$ is an
energy-preserving bijection from $\mathcal A_T(A)$ onto $\mathcal A_T(-A)$.
Thus $g_T(-A)=g_T(A)$, and the case $s<0$ follows by combining the two.
\end{proof}

\subsection{Proof of the lower bound in the homogenization theorem
}\label{subsection gamma liminf}
The proof of Theorem \ref{thm:sub:main} is subdivided into this section and the following one for the lower and the upper limits, respectively.

\smallskip
The proof of the lower bound can be achieved by the blow-up argument 
of Fonseca and M\"uller \cite{FM}.  Its formalization adapted to homogenization 
problems can be found in \cite{BMS}, to which we refer for details. 
In this particular case, we note that the domain of the $\Gamma$-limit is trivially contained
in $\WW^{1,p}(\Omega,\RR^2)$. In order to prove the lower bound, we fix $u_\e\to u$ weakly in $\WW^{1,p}(\Omega,\RR^2)$
with $\liminf_{\e\to0}F^p_\e(u_\e)<+\infty$. Upon extracting subsequences, we 
can suppose that this liminf is actually a limit, and that the measures
$$
\mu_\e(B)=\int_B|\nabla u_\e|^p\dd x
$$
weakly$^*$-converge to some $\mu$. It suffices then to prove that for $\mathcal{L}^2$-almost every $x_0$ we have
\begin{equation}\label{loba}
\frac{\dd\mu}{\dd x}(x_0)\ge f^p_{\rm hom}(\nabla u(x_0)).
\end{equation}
Select a point $x_0 \in \Omega$ which is a Lebesgue point of $u$ and $\nabla u$, a point of approximate differentiability
of $u$, a point of almost everywhere convergence of $u_\e \to u$, and such that the Radon-Nikodým derivative of $\mu$ with
respect to $\mathcal{L}^2$ exists and is finite. Note that $\mathcal{L}^2$-almost every $x_0 \in \Omega$ satisfy these
properties.
For simplicity of notation, upon a translation, we can suppose that $x_0=0$. Denote $A=\nabla u(0)$.

 We first treat the case of $\text{rank}(A)=2$.
 Write
\begin{eqnarray*}
\frac{\dd\mu}{\dd x}(0)=\lim_{\rho\to 0} \frac{\mu(\rho Q)}{\rho^2}=
\lim_{\rho\to 0}\liminf_{\e\to 0} \frac{\mu_\e(\rho Q)}{\rho^2}= 
\lim_{\rho\to 0}\liminf_{\e\to 0} \frac{1}{\rho^2}\int_{\rho Q}|\nabla u_\e|^p\dd x,
\end{eqnarray*}
noting that for this liminf is actually a limit except for a countable number of $\rho$.
We can then find $\e=\e(\rho)$ such that 
\begin{eqnarray*}
\frac{\dd\mu}{\dd x}(0)=
\lim_{\rho\to 0}\frac{1}{\rho^2}\int_{\rho Q}|\nabla u_\e|^p\dd x.
\end{eqnarray*}
We can then consider the scaled functions
$$
v_\e(y)=\frac{1}{\rho} u_\e (\rho y) -\frac{\e}{\rho} \left \lfloor \frac{u(0)}{\e} \right \rfloor, \quad y \in Q
$$
so that $v_\e$ takes values in $\frac{\e}{\rho} \mathcal{Z}$ (and the integer translation is irrelevant) and
\begin{eqnarray*}
\frac{\dd\mu}{\dd x}(0)=
\lim_{\rho\to 0}\int_{Q}|\nabla v_\e|^p\dd y.
\end{eqnarray*}
We have that $v_\e$ converges to $v(y)=Ay$ in $\LL^p(Q,\RR^2)$, since $0$ is a point of approximate differentiability so that
$$
\tilde v_\e(y)=\frac{u_\e(\rho y)-u(0)}{\rho}
$$
tends to $Ay$ in $\LL^p(Q,\RR^2)$, see e.g.~Remark 3.72 in \cite{AFP},  and $v_\e-\tilde v_\e = \frac{\e}{\rho} \lambda_\e\to 0$ uniformly, where $0\leq \lambda_\e \leq 1$. It is then not restrictive to suppose also that $\nabla v_\e$ weakly converges to the constant value $A=\nabla u(0)$, since $\nabla v_\e (y)=\nabla u_\e(\rho y) $ giving that the gradients $|\nabla v_\e|$ are equi-bounded in $\LL^p(Q,\RR^{2 \times 2})$.
If we set $\eta=\frac{\e}{\rho}$ we can apply Lemma \ref{lemma variation of BCs} 
with $\Omega=Q$, $\eta$ in place of $\e$,
$v_\e$ in place of $u_\e$, $v$ in place of $u$ and the functions $\overline u_\e$
given by $\eta \,u_A(y/\eta)$, where $u_A$ is as in Definition \ref{formula f_hom with A non singular}. We check that the assumption on such $\overline u_\e$ holds. Since $\overline u_\e (y) = \eta u_A (\frac{y}{\eta})$, then by Lemma \ref{rescaled P_0 satisfies  C|U| bound} we have 
\begin{equation}\label{UU}
    \lim_{\e \to 0} \int_U |\nabla \overline u_\e|^p \dd y =C_A|U|
\end{equation}
for all fixed open sets $U\subset Q$ with Lipschitz boundary.
Also, $\overline u_\e \to v(y)=Ay$ in $\LL^p(Q,\RR^2)$ since
$$
\int_{Q}| Ay-\eta P_0(A \tfrac{y}{\eta})|^p \dd y= \int_{Q} \eta^p \underbrace{\left| A\tfrac{y}{\eta}-P_0(A \tfrac{y}{\eta})\right|^p}_{\leq 1} \dd y \leq \eta^p \to 0,
$$
by the definition of $P_0$. By \eqref{UU} with $Q$ in place of $U$, the gradients are clearly equi-bounded, which gives the weak convergence of $\overline u_\e $ to $v$ in $\WW^{1,p}(Q,\RR^2)$.
We then deduce that we may suppose that $v_\e(y)=\eta \,u_A(y/\eta)$ in a neighbourhood
of $\partial Q$, and that 
\begin{gather*}
\frac{\dd\mu}{\dd x}(0)\ge
\liminf_{\eta \to 0}\inf\Bigl\{ \int_{Q}|\nabla w|^p\dd y:  
w\in \WW^{1,p}(Q,\eta\Z),\\
\hskip15em w(y)=\eta \,u_A(y/\eta)\hbox{ in a neighbourhood
of }\partial Q\Bigr\}.
\end{gather*}

By taking $T=\frac{1}{\eta}$ and rewriting the right-hand side of this inequality as a minimization 
problem on $TQ$; that is, by using the transformation
$$
u(x)= T w\left(\frac{x}{T}\right), \quad \nabla u(x)=\nabla w\left(\frac{x}{T}\right),
$$
 changing variables in the integral, we obtain the lower bound (\ref{loba}). Note that we chose coordinate cubes $Q$ rather than the tilted cubes $Q_A$ since the lower bound is indifferent to the choice; the two formulations agree by Lemma \ref{cubes and losanghe lemma}.  \smallskip

We now deal with the singular case of $\text{rank}(A)=1$. Denote $A=\nabla u(0)$, with $A=a \otimes \nu$. Let $Q_\nu:=(0,1)\nu + (0,1)\nu^\perp$ be the unit cube with sides parallel to $\nu$ and $\nu^\perp$, and write
\begin{eqnarray*}
\frac{\dd\mu}{\dd x}(0)&=& \lim_{\rho\to 0} \frac{\mu(\rho Q_\nu)}{|\rho Q_\nu|} =\lim_{\rho\to 0} \frac{1}{\rho^2} \mu(\rho Q_\nu)\\
&=&
\lim_{\rho\to 0}\liminf_{\e\to 0} \frac{1}{\rho^2} \mu_\e(\rho Q_\nu)= 
\lim_{\rho\to 0}\liminf_{\e\to 0} \frac{1}{\rho^2}\int_{\rho Q_\nu}|\nabla u_\e|^p\dd x,
\end{eqnarray*}
noting that this liminf is actually a limit except for a countable number of $\rho$.
We can then find $\e=\e(\rho)$ such that 
\begin{eqnarray*}
\frac{\dd\mu}{\dd x}(0)=
\lim_{\rho\to 0}\frac{1}{\rho^2}\int_{\rho Q_\nu}|\nabla u_\e|^p\dd x.
\end{eqnarray*}
Define now the one-dimensional restriction of $u_\e$ on $Q_\nu$
$$
u_\e^{y}(t)=u_\e(y\nu^\perp+t\nu), \quad t \in (0,\rho), 
$$
which is absolutely continuous for $\mathcal{L}^1$-almost every $y \in (0,\rho)$. For fixed $\rho$ it holds
\begin{eqnarray}\label{first slicing}
\frac{1}{\rho^2}\int_{\rho Q_\nu}|\nabla u_\e|^p\dd x&
\ge& \frac{1}{\rho^2}\int_{\rho Q_\nu}|\nabla u_\e\cdot \nu|^p\dd x=
\frac{1}{\rho^2}\int_0^\rho\int_0^\rho
\Big|\frac{\dd u_\e^{y}(t)}{\dd t}\Big|^p\dd t \, \dd y.
\end{eqnarray}
Since $u_\e^{y}(t)\in \e \mathcal Z$ for $\mathcal{L}^1$-almost all $y$, taking into account that the image of $u_\e^{y}$ is a polygonal in $\e \mathcal Z$, we now prove that for fixed $\rho$ and $\mathcal{L}^1$-almost all $y \in (0,\rho)$
\begin{equation}\label{inequality step 1 liminf}
\int_0^\rho
\Big|\frac{\dd u_\e^{y}(t)}{\dd t}\Big|^p\dd t\ge \rho \Big(\Big|\frac{a^{\e,\rho,y}_1}\rho\Big|
+\Big|\frac{a^{\e,\rho,y}_2}\rho\Big|\Big)^p,
\end{equation}
where $a^{\e,\rho,y}:=(a^{\e,\rho,y}_1,a^{\e,\rho,y}_2)=u_\e^{y}(\rho)-u_\e^{y}(0)$.

\begin{figure}[ht!]
\begin{center}
\begin{tikzpicture}[scale=1]
\def\eps{0.5}

\draw[step=\eps, gray!40, thin] (0,0) grid (4,3);

\draw[->] (0,0) -- (4.5,0);
\draw[->] (0,0) -- (0,3.5);

\draw[thick, blue]
(0,0) -- (1,0) -- (1,1) -- (2,1) -- (2,2) -- (2.7,2);

\fill (0,0) circle (2pt);
\fill (2.7,2) circle (2pt) node[above right]{$(a_1^{\varepsilon,\rho,y}, a_2^{\varepsilon,\rho,y})$};

\node[blue] at (2.5,0.5) {\textcolor{black}{$u_\varepsilon^y((0,\rho)) \subset \e \mathcal{Z}$}};

\end{tikzpicture}
\caption{Path of $u_\e^y$ in $\e \mathcal{Z}$}\label{piecewise path}
\end{center}
\end{figure}
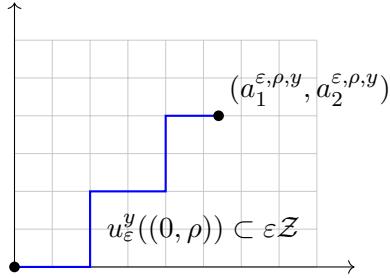

To this end, fix $\rho >0$ and $y \in (0,\rho)$ a point of absolute continuity of $t \mapsto u^y_\e(t)$. Given $\dfrac{\dd u_\e^{y}(t)}{\dd t}= (u_\e^{y})'(t)$ and $|(u_\e^{y})'(t)|^p= \big(|(u_\e^{y})'_1(t)|^2+(u_\e^{y})'_2(t)|^2\big)^{\frac{p}{2}} $, suppose first $u_\e^{y}$ has piecewise-constant speed, either horizontally with (signed) speed $v^\e_1$ or vertically with speed $v^\e_2$, say, in both cases in a finite union of intervals in $(0,\rho)$. There exists $\lambda\in(0,1)$ such that $|(u_\e^{y})'(t)|^p=|v^\e_1|^p$ in a set of $t$ of measure $\lambda\rho$, and  $|(u_\e^{y})'(t)|^p=|v^\e_2|^p$ in a  set of $t$ of measure $(1-\lambda)\rho$. Hence, we get
$$
\int_0^\rho
\Big|\frac{\dd u_\e^{y}(t)}{\dd t}\Big|^p\dd t= \lambda \rho |v^\e_1|^p + (1-\lambda)\rho |v^\e_2|^p
$$
with the relations
$$
\lambda \rho v^\e_1= a_1^{\e,\rho,y}, \quad (1-\lambda)\rho v^\e_2= a_2^{\e,\rho,y}  
$$
(see Fig.~\ref{piecewise path}). As depicted in the figure, the signs of $v_1^\e$ and $v_2^\e$ are given by the signs of the final point components $a_1^{\e,\rho,y}$ and $a_2^{\e,\rho,y}$. 
We consider the parameter $\lambda$ as free, with the speeds adjust accordingly. Thus, the integral rewrites as
\begin{eqnarray*}
\int_0^\rho
\Big|\frac{\dd u_\e^{y}(t)}{\dd t}\Big|^p\dd t&= &\lambda \rho |v^\e_1|^p + (1-\lambda)\rho |v^\e_2|^p= \frac{1}{(\lambda \rho)^{p-1}}|a_1^{\e,\rho,y}|^p + \frac{1}{((1-\lambda) \rho)^{p-1}}|a_2^{\e,\rho,y}|^p \\
&=& \rho\left(\frac{1}{\lambda^{p-1}} \left|\frac{a_1^{\e,\rho,y}}{\rho}\right|^p + \frac{1}{(1-\lambda)^{p-1}} \left|\frac{a_2^{\e,\rho,y}}{\rho}\right|^p \right) .
\end{eqnarray*}
Minimizing the last quantity in $\lambda$; that is, for $$\lambda_{\min}= \dfrac{|a_1^{\e,\rho,y}|}{|a_1^{\e,\rho,y}|+ |a_2^{\e,\rho,y}|},
$$ this expression simplifies to \eqref{inequality step 1 liminf}.

In the general case, where $(u_\e^y)'(t)$ is not piecewise constant, we can use Jensen's inequality on each subinterval of $(0,\rho)$ where $u_\e^y$
lies in one of the segments of $\e\mathcal Z\setminus
\e \left(\frac{1}{2}+\mathbb Z\right)^2$
\begin{eqnarray*}
\int_0^\rho
\Big|\frac{\dd u_\e^{y}(t)}{\dd t}\Big|^p\dd t &= &\int_0^{t_1}|(u^y_\e)'_1(t)|^p\dd t + \int_{t_1}^{t_2}|(u^y_\e)'_2(t)|^p\dd t+...+\int_{t_n}^\rho |(u^y_\e)'_1(t)|^p \dd t \\
&\geq& \frac{1}{t_1^{p-1}}[(u^y_\e)_1(t_1)-(u^y_\e)_1(0)]^p + \frac{1}{(t_2-t_1)^{p-1}}[(u^y_\e)_2(t_2)-(u^y_\e)_1(t_1)]^p\\&&
\hskip2cm\dots + \frac{1}{(\rho-t_n)^{p-1}}[(u^y_\e)_1(\rho)-(u^y_\e)_1(t_n)]^p\\
&=& \sum_{i=1}^{n+1} \frac{c_i^p}{\lambda_i^{p-1}},
\end{eqnarray*}
where $\lambda_i:=t_i-t_{i-1}$, (with $t_0=0, t_{n+1}=\rho$), $\sum_1^{n+1} \lambda_i= \rho $ and $c_i:= (u_\e^y)_k(t_{i})-(u_\e^y)_k(t_{i-1}) > 0$, for $k=1$ or $k=2$, with $\sum_1^{n+1}c_i= |a^{\rho,\e,y}_1|+|a^{\rho,\e,y}_2|$. This expression is minimized by
$$
\lambda_i= \rho \ \frac{c_i}{\sum_{j=1}^{n+1}c_j},
$$
and its minimum value is $\dfrac{1}{\rho^{p-1}}(\sum_{i=1}^{n+1}c_i)^p= \dfrac{1}{\rho^{p-1}}( |a^{\e,\rho,y}_1|+|a^{\e,\rho,y}_2|)^p$, again obtaining \eqref{inequality step 1 liminf}.

Finally, since $a^{\e,\rho,y}/\rho$ tends to $a$ for almost all $y$ as $\rho\to 0$, we can suppose it holds the limit
\begin{equation}\label{inequality step 2 liminf}
\lim_{\rho \to 0} \frac{1}{\rho} \int_0^\rho \left(\left|\frac{a_1^{\e,\rho,y}}{\rho}\right|+\left|\frac{a_2^{\e,\rho,y}}{\rho}\right|\right)^p \dd y =(|a_1|+|a_2|)^p,
\end{equation}
which entails 
\begin{equation}\label{lower bounded in rank 1 case}
\liminf_{\rho\to 0}\frac{1}{\rho^2}\int_{\rho Q_\nu}|\nabla u_\e|^p\dd x\ge (|a_1|+|a_2|)^p=f^p_{\rm hom}(A),
\end{equation}
as desired, by combining \eqref{first slicing} and \eqref{inequality step 1 liminf}. 

\subsection{Proof of the upper bound in the homogenization theorem
}\label{sec:ga-lims}
In order to prove the upper bound, we can proceed by a density argument. First, we exhibit
recovery sequences for a linear target function of gradient of maximum rank, and consequently by a gluing procedure, which makes use of Lemma \ref{lemma variation of BCs}, for piecewise-affine functions with gradients of maximum rank. After that, we prove in Proposition \ref{f hom is continuous} that the function $f_{\rm{hom}}^p$ is continuous, so that the functional in (\ref{gamma limit formula expression}) is strongly continuous in $\WW^{1,p}(\Omega,\RR^2)$. Due to this fact, the upper bound is proved for a general $u \in \WW^{1,p}(\Omega,\RR^2)$ by the density in norm of piecewise-affine functions.
\begin{proposition}\label{limsup for piecewise affine functions}
The upper bound of the $\Gamma$-limit of $F^p_\e$ with $F^p_{\rm hom}$ in \emph{(\ref{gamma limit formula expression})} holds for piecewise-affine functions with non-singular gradients. 
\end{proposition}
\begin{proof}
Let $u$ be the linear target function $u(x)=A x$ with $A \in \RR^{2\times 2}$ nonsingular, fix $\eta>0$, and let
$u_A(x)= P_0(A x)$ with $P_0$ as in Definition \ref{projlin}. Since the liminf in \eqref{full rank formula} is a limit,
there exist $T\in\mathbb{N}$ and $u^T\in \WW^{1,p}(TQ_A,\Z)$ with $u^T= u_A$ in a neighbourhood of $\partial TQ_A$ such that
$$
\frac{|\det A|}{T^2} \int_{TQ_A}|\nabla u^T|^p\,\dd x\le f^p_{\rm hom}(A)+\eta.
$$
Note that we chose the tilted cubes $TQ_A$: since $T\in\mathbb{N}$, the parallelograms $A^{-1}(T\mathbf{k})+TQ_A$,
$\mathbf{k}\in\ZZ^2$, tile the plane, and the translations $A^{-1}(T\mathbf{k})$ are compatible with $u_A$, in the sense
that $u_A(x-A^{-1}(T\mathbf{k}))+T\mathbf{k}=u_A(x)$, since $P_0(y+z)=P_0(y)+z$ for $z\in\ZZ^2$. Hence the function
$$
\tilde u(x):=u^T\big(x-A^{-1}(T\mathbf{k})\big)+T\mathbf{k}\qquad\hbox{if } x\in A^{-1}(T\mathbf{k})+TQ_A,\ \mathbf{k}\in\ZZ^2,
$$
coincides with $u_A$ in a neighbourhood of the boundary of each parallelogram, so that, as $\Z+\ZZ^2=\Z$,
$\tilde u\in \WW^{1,p}_{\rm loc}(\RR^2,\Z)$. Moreover, $\tilde u(x+A^{-1}(T\mathbf{k}))=\tilde u(x)+T\mathbf{k}$, so that
$\psi(x):=\tilde u(x)-Ax$ is periodic with respect to the lattice $A^{-1}(T\ZZ^2)$, with periodicity cell $TQ_A$,
and so are $\nabla\psi=\nabla\tilde u-A$ and $|\nabla\tilde u|^p$, as $A$ is constant.

We claim that the functions $u_\e(x):= \e\tilde u(x/\e)$ weakly converge to $u$ in $\WW^{1,p}(\Omega,\RR^2)$ and
$$
\lim_{\e\to 0} F^p_\e(u_\e)\le (f^p_{\rm hom}(A)+\eta)|\Omega|,
$$
which yields the desired inequality. Indeed, $u_\e(x)=Ax+\e\psi(x/\e)$, and by the Riemann--Lebesgue lemma
$\psi(\cdot/\e)$ is bounded in $\LL^p(\Omega,\RR^2)$, so that $u_\e\to u$ strongly in $\LL^p(\Omega,\RR^2)$; moreover
$$
\nabla u_\e(x)=A+\nabla\psi\Big(\frac x\e\Big)\rightharpoonup A+\frac{|\det A|}{T^2}\int_{TQ_A}\nabla\psi\,\dd x=A
\qquad\hbox{weakly in } \LL^p(\Omega,\RR^{2\times2}),
$$
since the gradient of a periodic function has zero average on the periodicity cell. Finally, since $u_\e(x)\in\e\Z$
almost everywhere and $|\nabla u_\e(x)|^p=|\nabla\tilde u(x/\e)|^p$, again by the Riemann--Lebesgue lemma
$$
\lim_{\e\to 0} F^p_\e(u_\e)=\lim_{\e\to0}\int_\Omega \Big|\nabla \tilde u
\Big(\frac x\e\Big)\Big|^p\, \dd x
=|\Omega|\,\frac{|\det A|}{T^2}\int_{TQ_A}|\nabla u^T|^p\,\dd x\le (f^p_{\rm hom}(A)+\eta)|\Omega|.
$$
This proves the upper bound in the case of a linear target function.

Let now $u \in \WW^{1,p}(\Omega,\mathbb R^2)$ be piecewise affine; that is,
$$
\Omega= \bigcup_{i=1}^m \Omega_i \cup N,
$$
where $\Omega_i$ are pairwise disjoint open subsets of $\Omega$, $|N|=0$ and
$$
u(x)=A_ix+b_i \quad \text{for } x \in \Omega_i, \quad i=1,\dots,m,
$$
with $A_i \in \mathbb R^{2 \times 2}$ non-singular and $b_i \in \mathbb R^2$.
Since the liminf in \eqref{full rank formula} is a limit, there exist $T\in\mathbb{N}$ and
$u_i^T\in \WW^{1,p}(TQ_{A_i},\Z)$ with $u_i^T= u_{A_i}$ in a neighbourhood of $\partial TQ_{A_i}$ such that
$$
\frac{|\det A_i|}{T^2} \int_{TQ_{A_i}}|\nabla u_i^T|^p\,\dd x\le f^p_{\rm hom}(A_i)+\eta,
$$
for all $i=1,\dots,m$. As in the case of one affine function, we extend $u_i^T$ to
$$
\tilde u_i(x):=u_i^T\big(x-A_i^{-1}(T\mathbf{k})\big)+T\mathbf{k}\qquad\hbox{if } x\in A_i^{-1}(T\mathbf{k})+TQ_{A_i},\ \mathbf{k}\in\ZZ^2,
$$
so that $\tilde u_i\in \WW^{1,p}_{\rm loc}(\RR^2,\Z)$ and $\tilde u_i(x)-A_ix$ is periodic with respect to the lattice
$A_i^{-1}(T\ZZ^2)$. The functions
$$
u^i_\e(x):= \e \tilde u_i\left( \frac{x}{\e}\right)+ \e \left \lfloor \frac{b_i}{\e} \right \rfloor  \in \WW^{1,p}_{\rm loc}(\mathbb R^2, \e \mathcal{Z}),
$$
where the integer part is taken componentwise, are good candidates when restricted on each $\Omega_i$, since
$\e \lfloor \frac{b_i}{\e} \rfloor - b_i= \e (  \lfloor \frac{b_i}{\e}  \rfloor  -\frac{b_i}{\e} )\to 0$, so that $\e \lfloor \frac{b_i}{\e} \rfloor \to b_i$ also in $\LL^p(\Omega_i,\mathbb R^2)$ and trivially in $\WW^{1,p}(\Omega_i,\mathbb R^2)$, and
$$
\e \tilde u_i\left( \frac{x}{\e}\right) \to A_ix \quad \text{weakly in }\WW^{1,p}(\Omega_i,\mathbb R^2)
$$
by the same proof as before, performed for one affine function $Ax$, now applied to the functions $A_ix$ and $\e \tilde u_i\left( \frac{x}{\e}\right)$, on $\Omega_i$, recalling that the Riemann-Lebesgue lemma works on every bounded open set. Thus, we have
$$
u_\e^i(x)= \e \tilde u_i\left( \frac{x}{\e}\right)+ \e \Big \lfloor \frac{b_i}{\e} \Big \rfloor \to A_ix+b_i \quad \text{weakly in }\WW^{1,p}(\Omega_i,\mathbb R^2)
$$
for all $i=1,\dots,m$. Moreover, since $\nabla u_\e^i(x)= \nabla \tilde u_i\left(\frac{x}{\e}\right)$ and
$|\nabla\tilde u_i|^p$ is periodic with periodicity cell $TQ_{A_i}$, by the Riemann--Lebesgue lemma
\begin{equation}\label{limsup inequality for each Omega_i}
\lim_{\e \to 0} \int_{\Omega_i} |\nabla u_\e^i(x)|^p \dd x = |\Omega_i|\,\frac{|\det A_i|}{T^2}\int_{TQ_{A_i}}|\nabla u_i^T|^p\,\dd x \leq (f^p_{\rm hom}(A_i)+\eta)|\Omega_i| .
\end{equation}

The issue now is that we cannot simply merge the $u_\e^i$ to a global $u_\e$ in all $\Omega$ because we cannot guarantee that they coincide $\mathcal{H}^1$ almost everywhere on the interfaces across the different $\Omega_i$. We need to modify $u_\e^i$ suitably while also keeping the inequality (\ref{limsup inequality for each Omega_i}) true. To this end, define
$$
\overline{u}_\e(x):= \e P_0\left(\frac{u(x)}{\e}\right).
$$
We note that $\overline{u}_\e \in \WW^{1,p}(\Omega, \e \mathcal{Z})$. Also, by
$$
\overline{u}_\e(x)=\e P_0\left(\frac{A_ix+b_i}{\e} \right)=A_ix+b_i+\e\,\theta\left(\frac{A_ix+b_i}{\e} \right) \quad \text{for } x \in \Omega_i, \quad i=1,\dots,m,
$$
where $\theta(y):=P_0(y)-y$ is bounded and $\ZZ^2$-periodic, we find that
$$
\overline{u}_\e(x) \to A_ix+b_i \quad \text{weakly in }\WW^{1,p}(\Omega_i,\mathbb R^2), \quad i=1,\dots,m,
$$
using the same argument performed for the case of one affine function.

We now apply the variation of boundary conditions Lemma \ref{lemma variation of BCs} on each $\Omega_i$ in place of $\Omega$, to each $u^i_\e$  and $\overline{u}_\e$ in common for all $i=1,\dots,m$, to find a sequence $v^i_\e \in \WW^{1,p}(\Omega_i, \e\Z)$ with
$$
\limsup_{\e\to0} \int_{\Omega_i}|\nabla v^i_\e|^p \dd x \le\limsup_{\e\to 0} \int_{\Omega_i}|\nabla u^i_\e|^p \dd x
$$
and $v^i_\e=\overline u_\e$ on a neighbourhood of $\partial\Omega_i$ (at the very end of the proof of Lemma \ref{lemma variation of BCs} take the $\limsup$ instead of $\liminf$ on both sides). By the inequality above, the gradients of $v_\e^i$ are equi-bounded in $\LL^p(\Omega_i,\RR^{2 \times 2})$, as the gradients of $u_\e^i$ are, thus we check only the convergence in $\LL^p(\Omega_i,\RR^{2})$ of the functions $v_\e^i$ to $A_ix+b_i$. From the explicit form of $v_\e^i$ in the proof of Lemma \ref{lemma variation of BCs},
$$
v_\e^i=P_\e(\varphi u_\e^i +(1-\varphi)\overline{u}_\e),
$$
where $\varphi$ is a smooth cutoff function. Now we know that $\varphi u_\e^i +(1-\varphi)\overline{u}_\e \to A_ix+b_i$ and $v_\e^i$ converges to its same limit since
$$
\lim_{\e\to 0}
\int_{\Omega_i}|v_\e^i-(\varphi u_\e^i +(1-\varphi)\overline{u}_\e )|^p \dd x \leq \lim_{\e\to 0}(\sqrt2\,\e)^p|\Omega_i| =0
$$
by definition of $P_\e$.

Now we can merge $v_\e^i$, since they coincide with the common function $\overline{u}_\e$ close to the boundaries of $\Omega_i$; that is, we define in $\Omega$
$$
u_\e(x):= v_\e^i(x), \quad 
\hbox{if }x \in \Omega_i, \ i\in\{1,\dots,m\},
$$
which gives a function $u_\e \in \WW^{1,p}(\Omega, \e \mathcal{Z})$. We check that
$u_\e \to u$ weakly in $\WW^{1,p}(\Omega,\mathbb R^2)$,
by the linearity of the integral and the convergence of each $v_\e^i$. Finally,
\begin{eqnarray*}
\limsup_{\e \to 0} F_\e^p(u_\e)&=& \limsup_{\e \to 0} \int_\Omega |\nabla u_\e|^p \dd x = \limsup_{\e \to 0} \left( \sum_{i=1}^m \int_{\Omega_i} |\nabla v_\e^i|^p \dd x \right) \\ &\leq& \sum_{i=1}^m \limsup_{\e \to 0}\int_{\Omega_i} |\nabla v_\e^i|^p \dd x \leq \sum_{i=1}^m \limsup_{\e \to 0}\int_{\Omega_i} |\nabla u_\e^i|^p \dd x \\ &\leq& \sum_{i=1}^m (f^p_{\rm hom}(A_i)+\eta)|\Omega_i|
=\int_{\Omega} f^p_{\rm hom}(\nabla u)\,\dd x + \eta|\Omega|,
\end{eqnarray*}
which concludes the proof by sending $\eta \to 0$.\end{proof}

\begin{proposition}\label{f hom is continuous}
The function $f^p_{\rm hom}$ as in Definition \rm\ref{formula f_hom with A non singular} is continuous. 
\end{proposition}
\begin{proof} Since the formula defining $f^p_{\rm hom}(A)$ for $A$ of full rank is continuous with respect to the boundary conditions and the set of full-rank matrices is an open set, 
it suffices to prove the continuity for the limit matrices $A$ that are singular. Due to the explicit form of $f^p_{\rm hom}(A)$ when $A$ is of rank $1$, we may also limit the analysis to approximating matrices of rank $2$.

\smallskip\noindent
\textbf{Proof of lower semicontinuity.}
Let $A=a \otimes \nu$ and $A_\eta=A + B_\eta$ with $B_\eta \to 0$ as $\eta \to 0$. We want to show 
$$
\liminf_{\eta \to 0}f_{\rm hom}^p(A_\eta) \geq (|a_1|+|a_2|)^p=f_{\rm hom}^p(A).
$$

We preliminarily remark that in the definition of the homogenized energy density, we can use $Q_\nu$ with sides parallel to $\nu, \nu^\perp$ instead of $Q$ since the proof of Lemma \ref{cubes and losanghe lemma} works the same when an orthogonal transformation is applied; namely, we can write
$$
f^p_{\rm hom}(A_\eta)= \liminf_{T \to \infty}\frac{1}{T^2}\inf\left\{ \int_{TQ_\nu}|\nabla u|^p \, \dd x: u \in \WW^{1,p}(TQ_\nu,\mathcal{Z}), \ u = u_{A_\eta} \sim \partial TQ_\nu      \right\}.
$$
Using this definition of $f^p_{\rm hom}(A_\eta)$, 
from any competitor of the infimum problem $u$ define the rescaled function
$$
v_T(x)=\frac{1}{T}u(T x), \quad x \in Q_\nu,
$$
which takes values in $\frac{1}{T}\mathcal{Z}$, and record
$$
\frac{1}{T^2}\int_{TQ_\nu}|\nabla u|^p \, \dd x= \int_{Q_\nu}|\nabla v_T|^p \, \dd x, \quad v_T(x)=\frac{1}{T}u_{A_\eta}(Tx) \sim \partial Q.
$$
Following the same slicing procedure used for the proof of the lower bound in the rank-1 case, we define the one-dimensional restriction of $v_T$ on $Q_\nu$
$$
v_T^{y}(t)=v_T(y\nu^\perp+t\nu), \quad t \in (0,1), 
$$
which is absolutely continuous for $\mathcal{L}^1$-almost every $y \in (0,1)$. It holds
\begin{eqnarray*}
\int_{ Q_\nu}|\nabla v_T|^p \, \dd x&
\ge& \int_{ Q_\nu}|\nabla v_T \cdot \nu|^p\, \dd x=\int_0^1\int_0^1
\Big|\frac{\dd v_T^{y}(t)}{\dd t}\Big|^p\, \dd t \, \dd y,
\end{eqnarray*}
and by the same minimization performed to prove inequality (\ref{inequality step 1 liminf}), here with $\rho=1$, we get
$$
\int_0^1\int_0^1
\Big|\frac{\dd v_T^{y}(t)}{\dd t}\Big|^p\, \dd t \, \dd y \geq \int_0^1 \big(|a_{\eta,1}^{T,y}|+|a_{\eta,2}^{T,y}|\big)^p \, \dd y,
$$
where
$$
a_\eta^{T,y}:=(a_{\eta,1}^{T,y},a_{\eta,2}^{T,y})=v_T^y(1)-v_T^y(0). 
$$
Since we have the boundary constraint $v_T(x)=\frac{1}{T}u_{A_\eta}(Tx)$ close to $\partial Q$, and the points $(0,y)$ and $(1,y)$ are in $\partial Q$ for all $y \in(0,1)$, we have that
$$
a_\eta^{T,y}=\frac{1}{T}u_{A_\eta}(T(y\nu^\perp+\nu))-\frac{1}{T}u_{A_\eta}(T(y\nu^\perp)), \quad y \in(0,1)
$$
and, moreover, the uniform bound
\begin{eqnarray*}
\left| a_\eta^{T,y}-A_\eta \nu \right|&=& \Big|\frac{1}{T}u_{A_\eta}(T(y\nu^\perp+\nu))-A_\eta(y\nu^\perp+\nu) - \frac{1}{T}u_{A_\eta}(T(y\nu^\perp)) +A_\eta(y\nu^\perp) \Big|\\&=&\frac{1}{T}\Big|u_{A_\eta}(T(y\nu^\perp+\nu))-A_\eta(T(y\nu^\perp+\nu))-u_{A_\eta}(T(y\nu^\perp)) +A_\eta(T(y\nu^\perp))  \Big|\\
&\leq& \frac{2}{T},
\end{eqnarray*}
since $u_A(x)=P_0(Ax)$ and $|P_0(Ax)-Ax|\leq 1$ by the definition of $P_0$. Therefore,
$$
\lim_{T \to \infty} a_\eta^{T,y}=A_\eta \nu= a+B_\eta \nu,
$$
uniformly in $y \in (0,1)$, and hence also in $\LL^p_y(0,1)$, yielding
$$
\lim_{T \to \infty} \int_0^1 \big(|a_{\eta,1}^{T,y}|+|a_{\eta,2}^{T,y}|\big)^p \, \dd y = (|a_1+(B_\eta \nu)_1|+|a_2+(B_\eta \nu)_2|)^p.
$$
By the arbitrariness of $u$, since the right-hand side does not depend on $u$, we get that
$$
f_{\rm hom}^p(A_\eta) \geq (|a_1+(B_\eta \nu)_1|+|a_2+(B_\eta \nu)_2|)^p,
$$
and finally
$$
\liminf_{\eta \to 0}f_{\rm hom}^p(A_\eta) \geq (|a_1|+|a_2|)^p.
$$
More specifically, take for each $T$ an almost minimizer $\overline{u}$ such that
$$
\frac{1}{T^2}\int_{TQ_\nu}|\nabla \overline{u}|^p \, \dd x -\frac{1}{T} \leq \inf\left\{ \frac{1}{T^2}\int_{TQ_\nu}|\nabla u|^p \, \dd x: u \in \WW^{1,p}(TQ,\mathcal{Z}), \ u = u_{A_\eta} \sim \partial TQ_\nu      \right\}.
$$
Then
\begin{eqnarray*}
f_{\rm hom}^p(A_\eta) &\geq& \liminf_{T \to \infty}\left(\frac{1}{T^2}\int_{TQ_\nu}|\nabla \overline{u}|^p \, \dd x -\frac{1}{T} \right)  =\liminf_{T \to \infty} \int_{Q_\nu}|\nabla \overline{v}_T|^p \, \dd x \\ &\geq&  \liminf_{T \to \infty}\int_0^1\int_0^1
\Big|\frac{\dd v_T^{y}(t)}{\dd t}\Big|^p\, \dd t \, \dd y \geq \lim_{T \to \infty} \int_0^1 \big(|a_{\eta,1}^{T,y}|+|a_{\eta,2}^{T,y}|\big)^p \, \dd y \\&=& (|a_1+(B_\eta \nu)_1|+|a_2+(B_\eta \nu)_2|)^p,  
\end{eqnarray*}
and we conclude by taking the $\liminf$ as $\eta \to 0$ on both sides.

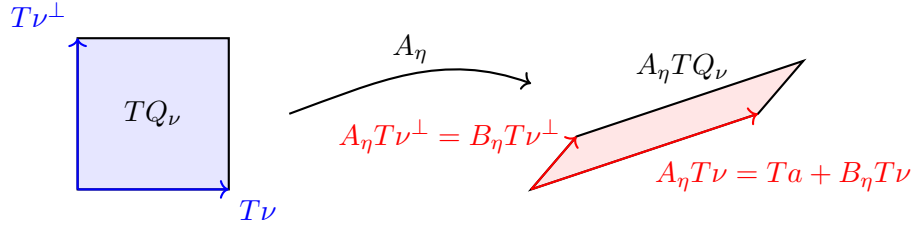
\begin{figure}[ht!]
\centering
\begin{tikzpicture}[scale=2]

\coordinate (nu) at (1,0);
\coordinate (nuperp) at (0,1);

\draw[thick, fill=blue!10]
    (0,0) -- (nu) -- ($(nu)+(nuperp)$) -- (nuperp) -- cycle;
\node at (0.5,0.5) {$TQ_\nu$};

\draw[->, thick, blue] (0,0) -- (nu) node[below right] {$T\nu$};
\draw[->, thick, blue] (0,0) -- (nuperp) node[above left] {$T\nu^\perp$};

\begin{scope}[xshift=3cm]
\coordinate (Tnu) at (1.5,0.5);
\coordinate (Tnuperp) at (0.3,0.35);

\draw[thick, fill=red!10]
    (0,0) -- (Tnu) -- ($(Tnu)+(Tnuperp)$) -- (Tnuperp) -- cycle;
\node at (1.0,0.8) {$A_\eta TQ_\nu$};

\draw[->, thick, red] (0,0) -- (Tnu)
    node[pos=0.5, below right] {$A_\eta T\nu=Ta+ B_\eta T\nu$};
\draw[->, thick, red] (0,0) -- (Tnuperp)
    node[left=2pt] {$A_\eta T\nu^\perp =B_\eta T\nu^\perp$};
\end{scope}

\draw[->, thick]
    (1.4,0.5) to[out=20,in=160] node[above] {$A_\eta$} (3.0,0.7);

\end{tikzpicture}
\caption{Stretched image of the cube $TQ_\nu$ under $A_\eta$.}
\label{fig:stretched-image}
\end{figure}

\smallskip\noindent
\textbf{Proof of upper semicontinuity. }Let $A=a \otimes \nu$ and let $A_\eta=A + B_\eta$ with $B_\eta \to 0$ as $\eta \to 0$. The behaviour of $A_\eta$ in terms of the deformation of the reference cube is pictured in Fig.~\ref{fig:stretched-image}: since $A_\eta$
becomes singular as $\eta\to0$, the parallelogram $A_\eta TQ_\nu$ is
increasingly stretched and flattened.
We want to show that
$$
\limsup_{\eta \to 0}f_{\rm hom}^p(A_\eta) \leq f_{\rm hom}^p(A)=(|a_1|+|a_2|)^p.
$$

It is sufficient to prove that for each $\eta>0$ there exists $u^T_\eta \in \WW^{1,p}(TQ_\nu,\Z)$,
with $T \gg 1$ such that $\dfrac{1}{T}u_\eta^T(Tx) \to A_\eta x \text{ weakly in }\WW^{1,p}(Q_\nu,\mathbb R^2)$ as $T\to \infty$ and
$$
\frac{1}{T^2}\int_{TQ_\nu}|\nabla u^T_\eta|^p\,\dd x \leq (|a_1|+|a_2|)^p + o_{\eta \to 0}(1).
$$
The final step, after proving the inequality above, in order to conclude the upper semicontinuity will follow by using Lemma \ref{lemma variation of BCs} to match the boundary conditions.  
Note that we are using $Q_\nu$ with sides parallel to $\nu, \nu^\perp$ instead of $Q$, but the proof of Lemma \ref{cubes and losanghe lemma} works the same when an orthogonal transformation is applied.

\begin{figure}[ht!]
\centering
\begin{tikzpicture}[scale=2.7]

\def\nstripes{7}
\def\jstripe{4}
\def\zx{1.3}
\def\zw{1.45}
\def\zh{0.62}
\def\kh{0.17}
\colorlet{Sjfill}{orange!30}
\colorlet{Sjline}{orange!75!black}

\pgfmathtruncatemacro{\nlast}{\nstripes-1}
\pgfmathsetmacro{\ylo}{\jstripe/\nstripes}
\pgfmathsetmacro{\yhi}{(\jstripe+1)/\nstripes}
\pgfmathsetmacro{\zy}{(\ylo+\yhi-\zh)/2}
\pgfmathsetmacro{\zt}{\zy+\zh}
\pgfmathsetmacro{\zs}{\zt-\kh}
\pgfmathsetmacro{\zr}{\zx+\zw}

\foreach \i in {0,...,\nlast} {
    \pgfmathsetmacro{\ybot}{\i/\nstripes}
    \pgfmathsetmacro{\ytop}{(\i+1)/\nstripes}
    \ifodd\i
        \fill[blue!20] (0,\ybot) rectangle (1,\ytop);
    \else
        \fill[blue!10] (0,\ybot) rectangle (1,\ytop);
    \fi
}
\foreach \i in {1,...,\nlast} {
    \pgfmathsetmacro{\y}{\i/\nstripes}
    \draw[blue!50, thin] (0,\y) -- (1,\y);
}

\fill[Sjfill] (0,\ylo) rectangle (1,\yhi);
\draw[thick, blue!70!black] (0,0) rectangle (1,1);
\draw[thick, Sjline] (0,\ylo) rectangle (1,\yhi);
\node at (0.5,{(\ylo+\yhi)/2}) {$S_j$};

\node at (0.5,{2.5/\nstripes}) {$TQ_\nu$};

\draw[->, thick, blue!80!black] (0,0) -- (1,0) node[below right] {$T\nu$};
\draw[->, thick, blue!80!black] (0,0) -- (0,1) node[above left] {$T\nu^\perp$};

\draw[decorate, decoration={brace, amplitude=4pt}, thick]
    (-0.05,\ylo) -- (-0.05,\yhi)
    node[midway, left=5pt] {$\dfrac{1}{|B_\eta|}$};

\fill[orange!8] (1,\ylo) -- (\zx,\zy) -- (\zx,\zt) -- (1,\yhi) -- cycle;
\draw[densely dashed, Sjline] (1,\ylo) -- (\zx,\zy)  (1,\yhi) -- (\zx,\zt);

\fill[Sjfill] (\zx,\zy) rectangle (\zr,\zs);
\shade[bottom color=Sjfill, top color=orange!65] (\zx,\zs) rectangle (\zr,\zt);
\draw[Sjline] (\zx,\zs) -- (\zr,\zs);
\draw[thick, Sjline] (\zx,\zy) rectangle (\zr,\zt);

\node[font=\small] at ({(\zx+\zr)/2},{(\zy+\zs)/2}) {$u$ constant in $y$};
\node[font=\small] at ({(\zx+\zr)/2},{(\zs+\zt)/2}) {transition in $y$};
\node[above] at ({(\zx+\zr)/2},\zt) {$S_j$};

\draw[decorate, decoration={brace, amplitude=4pt}, thick]
    ({\zr+0.04},{\zs-0.005}) -- ({\zr+0.04},\zy)
    node[midway, right=5pt] {$\dfrac{1}{|B_\eta|}-k$};
\draw[decorate, decoration={brace, amplitude=3pt}, thick]
    ({\zr+0.04},\zt) -- ({\zr+0.04},{\zs+0.005})
    node[midway, right=4pt] {$k$};

\draw[->, thick, blue!80!black]
    ({\zr+0.85},{\zy+0.05}) -- ({\zr+0.85},{\zt-0.05})
    node[midway, right=3pt] {$\nu^\perp$};

\end{tikzpicture}
\caption{Stripe-by-stripe construction on the cube $TQ_\nu$. Left: $TQ_\nu$ is
sliced into stripes parallel to $\nu$, of thickness $1/|B_\eta|$; $S_j$ is one
of them. Right: enlargement of $S_j$, split into a part of thickness
$1/|B_\eta|-k$, where $u_\eta^T$ is constant in $y$, and a thin layer of thickness $k$,
where the transition in $y$ takes place.}
\label{fig:stripes}
\end{figure}
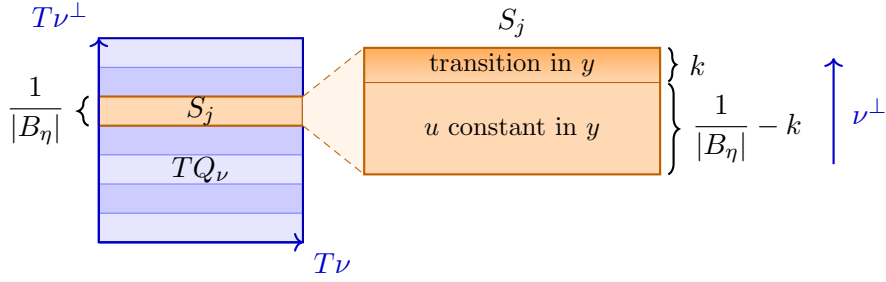

We now describe the construction of $u^T_\eta$. This will be done by subdividing the reference cube in stripes, and repeating a construction stripe by stripe. We parameterize the cube $TQ_\nu$ with $(t,y)$; that is,
$$
x = t\nu+y\nu^\perp, \quad (t,y)\in(0,T)\times(0,T).
$$
We want to construct the function
$$
u_\eta^T(x)=u_\eta^T(t\nu+y\nu^\perp), \quad (t,y)\in(0,T)\times(0,T).
$$
As shown in Figure~\ref{fig:stripes}, we divide the cube $TQ_\nu$ into a number of $|B_\eta| T$ stripes of thickness $1/|B_\eta|$ each, and without loss of generality we can assume $|B_\eta| T \in \mathbb N$ and $T \gg 1/|B_\eta|$, up to a negligible error. Then each stripe is subdivided in two regions of height $1/|B_\eta|-k$ and $k$, with $k \ll 1/|B_\eta|$. We build the function $ u_\eta^T$ on each stripe so that the derivative in $y$ is concentrated only on the $k$-stripe, whose energy will be negligible at the limit.

Fix the stripe $S_j:=(0,T) \times (y_j,y_{j+1})$, with $y_j=\frac{j}{|B_\eta|}$ and $j=0,\dots, |B_\eta| T -1$. In the lower and upper sides of the stripe we define 
$$t \mapsto u_\eta^T(t,y_j) \text{ and } t \mapsto u_\eta^T(t,y_{j+1})
$$ in such a way that it follows a ``zigzag" path (on $\mathcal{Z}$) almost along the image of $(0,T) \times \{y_j\}$ and of $(0,T) \times \{y_{j+1}\}$ through the linear map $A_\eta$, up to a uniform error. 

Define
\begin{eqnarray}
\tilde a^{j,0}&:=& \text{lattice point in }  \Big(\frac{1}{2}+\mathbb{Z}\Big)^2 \text{ closest  to }y_j B_\eta \nu^\perp  , \label{definition of initial point of vector a^j}
\\
\tilde a^{j,T}&:=& \text{lattice point in }  \left(\frac{1}{2}+\mathbb{Z}\right)^2 \text{ closest to }Ta+B_\eta T \nu + y_j B_\eta \nu^\perp . \label{definition of final point of vector a^j}
\end{eqnarray}
Both points belong to $\mathcal{Z}$.
These are (the projections of) the initial and final point of the image  through the map $A_\eta$ of the lower side of $S_j$.
Note that we can estimate the length of the vector from $\tilde a^{j,0}$ to $\tilde a^{j,T}$, $a^{j}:=\tilde a^{j,T} -\tilde a^{j,0}$ with
\begin{equation}\label{a eta T estimate vector part 2}
|a^j-Ta|=|\tilde a^{j,T}-(Ta+B_\eta T \nu+y_j B_\eta \nu^\perp)+ y_j B_\eta \nu^\perp - \tilde a^{j,0} + B_\eta T \nu| \leq 2+T|B_\eta|.
\end{equation}

The above estimate holds component-wise, 
\begin{equation}\label{a eta T estimate component wise part 2}
|a_i^{j}-Ta_i|\leq |a^{j}-Ta| \leq 2+T|B_\eta|, \quad i=1,2.    
\end{equation}

Let $v^j\colon[0,T]\to\Z$ be the path given by Lemma~\ref{lem:staircase} (with
$\e=1$ and $L=T$) from $\tilde a^{j,0}$ to $\tilde a^{j,T}$, see Fig.~\ref{piecewise path algorithm}, and set
\begin{equation}\label{v_j definition construction zig zag path}
u_\eta^T(t,y_j):=v^j(t),\qquad t\in(0,T).
\end{equation}

\begin{figure}[ht!]
\begin{center}
\begin{tikzpicture}[scale=1]
\def\eps{0.5}

\draw[step=\eps, gray!40, thin] (0,0) grid (4,3);

\draw[->] (0,0) -- (4.5,0);
\draw[->] (0,0) -- (0,3.5);

\draw[thick, blue]
(0,0) -- (0.5,0) -- (0.5,0.5) -- (1.5,0.5) -- (1.5,1) -- (2,1) -- (2,1.5) -- (2.5,1.5) -- (3,1.5) -- (3,2);

\fill (0,0) circle (2pt) node [below]{$\tilde a^{j,0}$};
\fill (3,2) circle (2pt) node[right]{$(\tilde a_1^{j,T}, \tilde a_2^{j,T})$};

\fill (0.1,0.2) circle (2pt);
\node at (-0.3,0.2) [left]{$ y_j B_\eta \nu^\perp$};
\draw[->, decorate, decoration={snake, post length=1mm}] (-0.4,0.2) -- (0.0,0.2);

\fill (2.8,2.2) circle (2pt);
\node at (2.2,2.2) [left]{$Ta+B_\eta T \nu +y_j B_\eta \nu^\perp$};
\draw[->, decorate, decoration={snake, post length=1mm}] (2.2,2.2) -- (2.7,2.2);

\draw[-, thick, red] (0,0) -- (3,2) node[pos=0.5, above left] {$a^j$};

\end{tikzpicture}
\caption{Construction of the path $v^j$ at the step $j$}\label{piecewise path algorithm}
\end{center}
\end{figure}
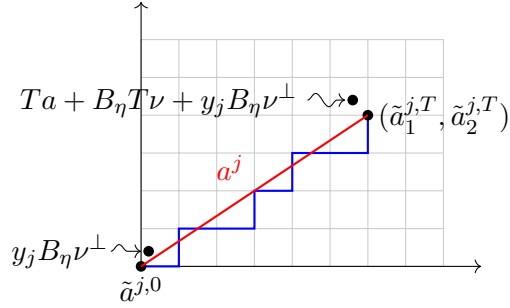

By Lemma~\ref{lem:staircase}(i), $v^j$ moves horizontally, with velocity
$v^j_1e_1$, on a set of measure $\lambda T$, and vertically, with velocity
$v^j_2e_2$, on the complementary set, where
$$
\lambda:=\frac{|a^j_1|}{|a^j_1|+|a^j_2|},\qquad
v^j_i:={\rm sign}(a^j_i)\,\frac{|a^j_1|+|a^j_2|}{T}\quad(i=1,2),
$$
so that $\lambda Tv^j_1=a^j_1$ and $(1-\lambda)Tv^j_2=a^j_2$. Moreover, by
Lemma~\ref{lem:staircase}(iii) and (ii), with
\begin{equation}\label{definition of vector line function z^j from a^j}
z^j(t):=\tilde a^{j,0}+\frac tT\,a^j,\qquad t\in(0,T),
\end{equation}
we have
\begin{equation}\label{lambda min and a^j integral}
\int_0^T\Big|\frac{\dd v^j}{\dd t}\Big|^p\,\dd t=\frac{(|a^j_1|+|a^j_2|)^p}{T^{p-1}}
\end{equation}
and
\begin{equation}\label{L infty bound of zigzag path and true vector line}
|z^j(t)-v^j(t)|\le\sqrt2\qquad\hbox{for all }t\in(0,T).
\end{equation}

The definition of 
$$
t \mapsto u_\eta^T(t,y_{j+1})=:v^{j+1}(t)
$$
follows exactly the same steps as $v^j(t)$. 
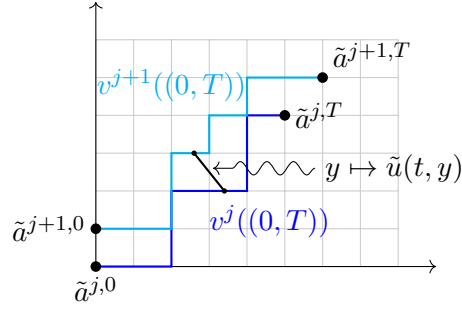
\begin{figure}[ht!]
\begin{center}
\begin{tikzpicture}[scale=1]
\def\eps{0.5}

\draw[step=\eps, gray!40, thin] (0,0) grid (4,3);

\draw[->] (0,0) -- (4.5,0);
\draw[->] (0,0) -- (0,3.5);

\draw[thick, blue]
(0,0) -- (1,0) -- (1,1) -- (2,1) -- (2,2) -- (2.5,2);
\draw[thick, cyan]
(0,0.5) -- (0,0.5) -- (1,0.5) -- (1,1.5) -- (1.5,1.5) -- (1.5,2) -- (2,2) -- (2,2.5) --(3,2.5);
\draw[thick, black]
(1.7,1) -- (1.3,1.5);
\fill (1.7,1) circle (1pt);
\fill (1.3,1.5) circle (1pt);
\node at (2.9,1.3) [right]{$ y \mapsto \tilde u(t,y)$};
\draw[->, decorate, decoration={snake, post length=2mm}] (2.9,1.3) -- (1.55,1.3);
\fill (0,0) circle (2pt) node [below]{$\tilde a^{j,0}$};
\fill (2.5,2) circle (2pt) node[right]{$\tilde a^{j,T}$};

\fill (0,0.5) circle (2pt) node [left]{$\tilde a^{j+1,0}$};
\fill (3,2.5) circle (2pt) node[above right]{$\tilde a^{j+1,T}$};

\node at (1,2.4) {$\textcolor{cyan}{v^{j+1}((0,T))}$};
\node[blue] at (2.3,0.6) {\textcolor{blue}{$v^j((0,T))$}};
\end{tikzpicture}
\caption{Schematic paths of $v^j$ and $v^{j+1}$}\label{piecewise path previous second}
\end{center}
\end{figure}

Now, in the portion of the stripe $(0,T) \times (y_j,y_j+\frac{1}{|B_\eta|}-k)$, we let $u_\eta^T(t,y)$ be constant in $y$,
$$
u_\eta^T(t,y)=u_\eta^T(t,y_j), \quad (t,y) \in (0,T) \times \Big(y_j,y_j+\frac{1}{|B_\eta|}-k\Big).
$$
In the thin stripe $y\in (y_j + \frac{1}{|B_\eta|}-k,y_{j+1})$
where $y_{j+1}=y_{j}+\frac{1}{|B_\eta|}$, the function $u_\eta^T$ will also depend on $y$.
For $y \in (y_j+\frac{1}{|B_\eta|}-k, y_{j}+\frac{1}{|B_\eta|})$ and $t \in (0,T)$ define the affine interpolation
$$
\tilde u(t,y)=u_\eta^T(t,y_j) + \big(u_\eta^T(t,y_{j+1})-u_\eta^T(t,y_j) \big) \frac{y-(y_j + \frac{1}{|B_\eta|}-k)}{k}
$$
 (see Fig.~\ref{piecewise path previous second}). 
Note that $\tilde u \notin \mathcal{Z}$ in general, thus we set
$$
u_\eta^T(t,y)= P_1(\tilde u)(t,y), \quad (t,y) \in (0,T) \times \Big(y_j+\dfrac{1}{|B_\eta|}-k, y_{j}+\dfrac{1}{|B_\eta|}\Big), 
$$
where $P_1(\tilde u)$ is constructed as in Lemma \ref{lemma1}, and we point out two facts: 

\smallskip
    (i) for all $t \in(0,T)$ we have $P_1(\tilde u)(t,y_j)=u_\eta^T(t,y_j) \in \mathcal{Z}$ and $P_1(\tilde u)(t,y_{j+1})=u_\eta^T(t,y_{j+1}) \in \mathcal{Z}$, so $u_\eta^T$ it is well defined in the full stripe $S_j$ when passing from the lower side to the upper side of it;
    
\smallskip
    (ii) the proof and the estimate (with same constant $C(p)$) in Lemma \ref{lemma1} remain the same even if we change the domain $\Omega$, in this case with the cube $TQ_\nu$ or its subdomains. Note that we are using the operator $P_1$ and not the simple composition with $P_0$ as in Definition \ref{projlin}.

\smallskip
    Thanks to (ii),  we need only to estimate the gradient of $\tilde u$ as defined on its domain after using the inequality of the type
$$
\|\nabla u_\eta^T\|_p=\|\nabla P_1(\tilde u)\|_p\leq C(p) \| \nabla \tilde u\|_p,
$$
where the $\LL^p$ norm is tacitly understood in the domain where $u_\eta^T=P_1(\tilde u)$.  Note that
$$
\frac{\partial \tilde u}{\partial y}(t,y)= \frac{u_\eta^T(t,y_{j+1})-u_\eta^T(t,y_j)}{k}.
$$
We want to give a bound to $|u_\eta^T(t,y_{j+1})-u_\eta^T(t,y_j)|$ uniform in $t$ and $j$, in order to estimate the derivative above.
By definition (\ref{definition of vector line function z^j from a^j}) and the bound (\ref{L infty bound of zigzag path and true vector line}), 
\begin{eqnarray*}
&&\hskip-1cm|u_\eta^T(t,y_{j+1})-u_\eta^T(t,y_j)|\\&\leq& |z^j(t)-u_\eta^T(t,y_j)|+ |z^{j+1}(t)-z^j(t)|+|u_\eta^T(t,y_{j+1})-z^{j+1}(t)|\\
&\leq& 2\sqrt{2}+|z^{j+1}(t)-z^j(t)|\ =\  2\sqrt{2}+\Big|(a^{j+1}_1-a^{j}_1) \frac{t}{T}, (a^{j+1}_2-a^{j}_2) \frac{t}{T} \Big| 
\\ &\leq& 2\sqrt{2}+|a^{j+1}-a^j|\ =\  2\sqrt{2} + |(\tilde a^{j+1,T}-\tilde a^{j+1,0}) - (\tilde a^{j,T}- \tilde a^{j,0})| \\
&\leq& 2\sqrt{2}+|\tilde a^{j+1,T}-\tilde a^{j,T}| + |\tilde a^{j+1,0}- \tilde a^{j,0}| \\& =&
2\sqrt{2}+ \Big|\tilde a^{j+1,T}- (Ta+B_\eta T \nu + y_{j+1} B_\eta \nu^\perp)\\
&&\hskip3cm+ (Ta+B_\eta T \nu + y_j B_\eta \nu^\perp - \tilde a^{j,T} )+\frac{B_\eta}{|B_\eta|} \nu^\perp\Big| \\
&&\hskip1cm+ \Big|\tilde a^{j+1,0}-y_{j+1}B_\eta \nu^\perp + (y_jB_\eta \nu^\perp-\tilde a^{j,0}) + \frac{B_\eta}{|B_\eta|} \nu^\perp\Big|\\
&\leq& 2\sqrt{2}+6,
\end{eqnarray*}
where we used (\ref{definition of initial point of vector a^j}) and (\ref{definition of final point of vector a^j}) for the last triangle inequality. Therefore, we find
\begin{equation}\label{estimate "12" of y derivative of tilde u}
\left|\frac{\partial \tilde u}{\partial y}(t,y) \right|= \frac{|u_\eta^T(t,y_{j+1})-u_\eta^T(t,y_j)|}{k}\leq \frac{2\sqrt{2}+6}{k}\leq \frac{9}{k}.
\end{equation}
By $0 \leq \dfrac{y-(y_j + \frac{1}{|B_\eta|}-k)}{k} \leq 1$ we also see that
$$
\left|\frac{\partial \tilde u}{\partial t}(t,y) \right|\leq 2\left(\left| \frac{\dd v^j}{\dd t}\right|+\left|\frac{\dd v^{j+1}}{\dd t}\right|\right).
$$
Finally, we compute the averaged integral: 
\begin{eqnarray*}
\frac{1}{T^2}\int_{TQ_\nu} |\nabla u_\eta^T|^p \dd x& =& \frac{1}{T^2}\int_0^T \int_0^T  |\nabla u_\eta^T(t,y)|^p \dd t \, \dd y
\\&=& \frac{1}{T^2} \sum_{j=0}^{|B_\eta|T-1} \int_{y_j}^{y_{j+1}}\int_0^T  |\nabla u_\eta^T(t,y)|^p \dd t \, \dd y.
\end{eqnarray*}
For each $j$ the summand is divided as 
\begin{eqnarray*}
&&\hskip-1cm\int_{y_j}^{y_{j+1}}\int_0^T  |\nabla u_\eta^T(t,y)|^p \dd t \, \dd y\\&=&\int_{y_j}^{y_j+\frac{1}{|B_\eta|}-k}\int_0^T \left|\frac{\dd v^j(t)}{\dd t}\right|^p \dd t \, \dd y + \int_{y_j+\frac{1}{|B_\eta|}-k}^{y_{j+1}}\int_0^T |\nabla P_1 (\tilde u)(t,y)|^p \dd t \, \dd y \\
&=& \Big(\frac{1}{|B_\eta|}-k\Big)\left(
\lambda T |v^{j}_1|^p + (1-\lambda)T |v^{j}_2|^p \right) + \int_{y_j+\frac{1}{|B_\eta|}-k}^{y_{j+1}}\int_0^T |\nabla P_1 (\tilde u)(t,y)|^p \dd t \, \dd y \\
&=&\Big(\frac{1}{|B_\eta|}-k\Big)\frac{(|a_1^j|+|a_2^j|)^p}{T^{p-1}} + \int_{y_j+\frac{1}{|B_\eta|}-k}^{y_{j+1}}\int_0^T |\nabla P_1 (\tilde u)(t,y)|^p \dd t \, \dd y \\ &\leq&  \Big(\frac{1}{|B_\eta|}-k\Big)\frac{(|a_1^j|+|a_2^j|)^p}{T^{p-1}} + C(p) \int_{y_j+\frac{1}{|B_\eta|}-k}^{y_{j+1}}\int_0^T |\nabla \tilde u(t,y)|^p \dd t \, \dd y \\ &\leq&  \Big(\frac{1}{|B_\eta|}-k\Big)\frac{(|a_1^j|+|a_2^j|)^p}{T^{p-1}} + C(p) \int_{y_j+\frac{1}{|B_\eta|}-k}^{y_{j+1}}\int_0^T \left|\frac{\partial \tilde u}{\partial y}(t,y) \right|^p \dd t \, \dd y \\ &&+ C(p) \int_{y_j+\frac{1}{|B_\eta|}-k}^{y_{j+1}}\int_0^T \left|\frac{\partial \tilde u}{\partial t}(t,y) \right|^p \dd t \, \dd y\\
&=& I^j_1+I^j_2+I^j_3,
\end{eqnarray*}
where we used $p/2<1$ to split the full gradient norm $|\nabla \tilde u|^p$ and \eqref{lambda min and a^j integral} in the third equality. We handle the three terms separately.

\smallskip
\noindent
    \emph{First term $I_1^j$.} We recall estimate (\ref{a eta T estimate component wise part 2}) which gives
    \begin{eqnarray*}
\left( |a_1^{j}|+ |a_2^{j}| \right)^p&=& \left( |a_1^{j}|-T|a_1|+ |a_2^{j}|-T|a_2|+T|a_1|+T|a_2| \right)^p  \\
&\leq& (1+\delta)(T|a_1|+T|a_2|)^p+C_\delta\left||a_1^{j}|-T|a_1|+ |a_2^{j}|-T|a_2|\right|^p \\
&\leq& (1+\delta)T^p(|a_1|+|a_2|)^p +C_\delta 2^{p-1}\left(\left||a_1^{j}|-T|a_1|\right|^p+\left||a_2^{j}|-T|a_2| \right|^p\right) \\
&\leq&(1+\delta)T^p(|a_1|+|a_2|)^p + C_\delta 2^p (2+T|B_\eta|)^p.
\end{eqnarray*}
and therefore
\begin{eqnarray*}
    &&\hskip-1cm\frac{1}{T^2}\sum_{j=0}^{|B_\eta| T-1} \left(\frac{1}{|B_\eta|}-k\right)\frac{(|a_1^j|+|a_2^j|)^p}{T^{p-1}} \\ &\leq& \frac{1}{T^2}\left(T-|B_\eta|T k \right)\left((1+\delta)T(|a_1|+|a_2|)^p +  \frac{C_\delta 2^p}{T^{p-1}} (2+T|B_\eta|)^p \right)\\
    &=& (1-|B_\eta|k)(1+\delta)(|a_1|+|a_2|)^p + (1-|B_\eta|k)C_\delta 2^p \left(\frac{2}{T}+|B_\eta| \right)^p \\
    &=&(1+\delta)(|a_1|+|a_2|)^p + C_\delta o_{\eta \to 0}(1).
\end{eqnarray*}
\emph{Second term $I_2^j$.} 
$$
I^j_2=C(p) \int_{y_j+\frac{1}{|B_\eta|}-k}^{y_{j+1}}\int_0^T \left|\frac{\partial \tilde u}{\partial y}(t,y) \right|^p \dd t \, \dd y \leq C(p) kT \left(\frac{9}{k}\right)^p = C(p)9^p\frac{T}{k^{p-1}} ,
$$
and thus 
\begin{eqnarray*}
\frac{1}{T^2}\sum_{j=0}^{|B_\eta| T-1}C(p) \int_{y_j+\frac{1}{\eta}-k}^{y_{j+1}}\int_0^T \left|\frac{\partial \tilde u}{\partial y}(t,y) \right|^p \dd t \, \dd y \leq \frac{1}{T^2}\,|B_\eta| T \,  C(p)9^p\frac{T}{k^{p-1}}=C(p) 9^p \frac{|B_\eta|}{k^{p-1}}, 
\end{eqnarray*}
which goes to $0$ as $\eta \to 0$. 

\smallskip\noindent
\emph{Third term $I_3^j$.} Since we have
$$
\left|\frac{\partial \tilde u}{\partial t}(t,y)\right|^p\leq 2^p\left(\left| \frac{\dd v^j}{\dd t}\right|+\left|\frac{\dd v^{j+1}}{\dd t}\right|\right)^p\leq 2^{2p-1} \left| \frac{\dd v^j}{\dd t}\right|^p+2^{2p-1}\left|\frac{\dd v^{j+1}}{\dd t}\right|^p,
$$
we can follow the same computations performed for the first term, equally for $j$ and $j+1$, 
\begin{eqnarray*}
I_3^j&=& C(p) \int_{y_j+\frac{1}{|B_\eta|}-k}^{y_{j+1}}\int_0^T \left|\frac{\partial \tilde u}{\partial t}(t,y) \right|^p \dd t \, \dd y\\  &\leq& C(p)2^{2p-1} \int_{y_j+\frac{1}{|B_\eta|}-k}^{y_{j+1}}\int_0^T \left( \left| \frac{\dd v^j}{\dd t}\right|^p+\left|\frac{\dd v^{j+1}}{\dd t}\right|^p \right)\dd t \, \dd y \\& =&C(p)2^{2p-1} k \Big(\frac{(|a_1^j|+|a_2^j|)^p}{T^{p-1}}  + \frac{(|a_1^{j+1}|+|a_2^{j+1}|)^p}{T^{p-1}}  \Big)\\ &\leq& C(p)2^{2p-1}  k \Big(2(1+\delta)T(|a_1|+|a_2|)^p +  2\frac{C_\delta 2^p}{T^{p-1}} (2+T|B_\eta|)^p  \Big),  
\end{eqnarray*}
applying the summation
\begin{eqnarray*}
\frac{1}{T^2}\sum_{j=0}^{|B_\eta| T-1} I_3^j &\leq& \frac{|B_\eta|C(p)2^{2p-1}k}{T}\left(2(1+\delta)T(|a_1|+|a_2|)^p +  2\frac{C_\delta 2^p}{T^{p-1}}(2+T|B_\eta|)^p\right)\\ &=& C_1|B_\eta|k(1+\delta)(|a_1|+|a_2|)^p+ C_2C_\delta |B_\eta| k \Big(\frac{2}{T}+|B_\eta| \Big)^p  ,  
\end{eqnarray*}
and everything tends to $0$ as $\eta \to 0$ since $|B_\eta| \to 0$. In conclusion,
$$
\frac{1}{T^2} \int_{TQ_\nu} |\nabla u_\eta^T|^p \dd x \leq (1+\delta)(|a_1|+|a_2|)^p +C_\delta  o_{\eta \to 0}(1),
$$
which yields the claim sending first $\eta \to 0$ and then $\delta \to 0$.

We check that the weak-convergence statement holds. Define
$$
u_\e(x)=\e u_\eta^T\left(\frac{x}{\e}\right), \quad x \in Q_\nu,
$$
where $\e=1/T$. Note that $u_\e \in \WW^{1,p}(Q_\nu, \e \mathcal{Z})$. Then, for $\eta$ fixed,  it holds
$$
\lim_{\e \to 0} u_\e(x) = A_\eta x \quad \text{weakly in }\WW^{1,p}(Q_\nu,\mathbb R^2). 
$$
We need to prove that the gradients are equi-bounded and that the functions  converge in $\LL^p(Q_\nu,\mathbb R^2)$. We note that
$$
\int_{Q_\nu} |\nabla u_\e(x)|^p \dd x=\int_{Q_\nu} |\nabla u_\eta^T\left(\frac{x}{\e}\right)|^p \dd x = \frac{1}{T^2} \int_{TQ_\nu} |\nabla u_\eta^T|^p \dd x \leq (1+\delta)(|a_1|+|a_2|)^p +C_\delta  o_{\eta}(1).
$$
Thus, the gradients are equi-bounded in $\LL^p$. To prove the convergence of $u_\e$, we note that
\begin{eqnarray*}
\int_{Q_\nu} |u_\e(x)-A_\eta x|^p \dd x&=& \int_{Q_\nu} \left|\e u_\eta^T\left(\frac{x}{\e}\right)-\e A_\eta \frac{x}{\e}\right|^p \dd x\\ &=& \e^{2+p}\int_{\frac{1}{\e}Q_\nu}|u_\eta^T(x)-A_\eta x|^p \dd x= \frac{1}{T^{2+p}}\int_{TQ_\nu}|u_\eta^T(x)-A_\eta x|^p \dd x.    
\end{eqnarray*}
Generally speaking, the integrand $|u_\eta^T(x)-A_\eta x|$ is uniformly (in $x$) close to a constant, since $u_\eta^T$ is a zigzag almost running along the image lines through $A_\eta$, and therefore the integral is of order $O(T^{-p}) \to 0$ as $T \to \infty$. Indeed, for $x=t\nu+y\nu^\perp$, $(t,y) \in (0,T) \times (0,T)$, and $A_\eta=a \otimes \nu +B_\eta$,
$$
|u_\eta^T(x)-A_\eta x|= |u_\eta^T(t,y)-at-t B_\eta \nu- y B_\eta \nu^\perp |,    
$$
if $y \in (y_j,y_j+\frac{1}{|B_\eta|}-k)$ for some $j=0,\dots,|B_\eta|T-1$; that is, where $u_\eta^T(t,y)=u_\eta^T(t,y_{j})$; then, recalling (\ref{definition of vector line function z^j from a^j}),
\begin{eqnarray*}
&&\hskip-2cm|u_\eta^T(t,y_j)-at-t B_\eta \nu- y B_\eta \nu^\perp |\\
&\leq& |u_\eta^T(t,y_j)-z^j(t)|+ |z^j(t)-at-t B_\eta \nu- y B_\eta \nu^\perp | \\ &\leq&
\sqrt{2}+ |a^j \frac{t}{T} -at-tB_\eta \nu| + |\tilde a^{j,0} -yB_\eta \nu^\perp| \\ &\leq& \sqrt{2} + |\frac{t}{T}(a^j-Ta-B_\eta T \nu)| +|\tilde a^{j,0}-y_j B_\eta \nu^\perp| + |y_j B_\eta \nu^\perp- y B_\eta \nu^\perp| \\ &\leq& \sqrt{2} +\frac{t}{T}2+1+|y-y_j||B_\eta| \leq C,
\end{eqnarray*}
where we used the (rearranged) inequality (\ref{a eta T estimate vector part 2}), the bound \eqref{L infty bound of zigzag path and true vector line} and definition (\ref{definition of initial point of vector a^j}). If instead $y \in(y_j+\frac{1}{|B_\eta|}-k,y_{j+1})$, where we have the transition, we recall that
$$
u_\eta^T(t,y)=P_1(\tilde u)(t,y), \quad
\tilde u(t,y)=u_\eta^T(t,y_j) + \big(u_\eta^T(t,y_{j+1})-u_\eta^T(t,y_j) \big) \frac{y-(y_j + \frac{1}{|B_\eta|}-k)}{k}.
$$
By the estimates
\begin{gather*}
|u_\eta^T(t,y)-\tilde u(t,y)|= |P_1(\tilde u)(t,y)-\tilde u(t,y)|\leq 1, \\
|\tilde u(t,y)-u_\eta^T(t,y_j)|=\left|\big(u_\eta^T(t,y_{j+1})-u_\eta^T(t,y_j) \big) \frac{y-(y_j + \frac{1}{|B_\eta|}-k)}{k} \right|\leq 9, 
\end{gather*}
the second coming from \eqref{estimate "12" of y derivative of tilde u}, through a concatenated triangle inequality, we get 
\begin{eqnarray*}
&&\hskip-1.5cm|u^T_\eta(t,y) -at-t B_\eta \nu- y B_\eta \nu^\perp|  \\ &\leq& |u_\eta^T(t,y)-\tilde u(t,y)|+ |\tilde u(t,y)-u_\eta^T(t,y_j)| + |u_\eta^T(t,y_j)-at-t B_\eta \nu- y B_\eta \nu^\perp | \\
&\leq& 1+9+C.
\end{eqnarray*}
Therefore,
$$
|u_\eta^T(x)-A_\eta x|= |u_\eta^T(t,y)-at-t B_\eta \nu- y B_\eta \nu^\perp | \leq C 
$$
for all $x \in TQ_\nu$, and
$$
\int_{Q_\nu} |u_\e(x)-A_\eta x|^p \dd x= \frac{1}{T^{2+p}}\int_{TQ_\nu}|u_\eta^T(x)-A_\eta x|^p \dd x \leq C^p \frac{1}{T^p}=C^p\e^{p},
$$
which tends to $0$.
We can apply Lemma \ref{lemma variation of BCs} 
with $\Omega=Q_\nu$, the functions $u_\e\to u(x) =A_\eta (x)$, and the function $\overline u_\e$
given by $\e \,u_{A_\eta}(x/\e)$, where $u_{A_\eta}$ is as in Definition \ref{formula f_hom with A non singular}. We check that the assumption on such $\overline u_\e$ holds: since $\overline u_\e (x) = \e  u_{A_\eta} (\frac{x}{\e})$,  then for any fixed $U\subset Q_\nu$ with Lipschitz boundary, we have
\begin{equation}\label{above-2}
\lim_{\e \to 0}\int_U |\nabla \overline u_\e|^p \dd x =C_{A_\eta} |U|
\end{equation}
 by Lemma \ref{rescaled P_0 satisfies  C|U| bound}.
Also, $\overline u_\e \to u(x)=A_\eta x$ in $\LL^p(\Omega,\RR^2)$ since
$$
\lim_{\e\to 0} \int_{Q_\nu}\left| A_\eta x-\e P_0(A_\eta \tfrac{x}{\e})\right|^p \dd x= \lim_{\e\to 0} \int_{Q_\nu} \e^p \left| A_\eta\tfrac{x}{\e}-P_0(A_\eta \tfrac{x}{\e})\right|^p \dd y =0,
$$
as $\left| A_\eta\tfrac{x}{\e}-P_0(A_\eta \tfrac{x}{\e})\right|\le 1$
by the definition of $P_0$. By \eqref{above-2} above with $\Omega=Q_\nu$ in place of $U$, we deduce that the gradients are equi-bounded, which gives the weak convergence of $\overline u_\e $ to $u$ in $\WW^{1,p}(Q_\nu,\mathbb R^2)$. Thus, thanks to Lemma \ref{lemma variation of BCs} we find $v_\e \in \WW^{1,p}(Q_\nu,\e \mathcal{Z})$ with 
$$
\liminf_{\e \to 0} \int_{Q_\nu}|\nabla v_\e|^p \dd x \leq \liminf_{\e \to 0} \int_{Q_\nu}|\nabla u_\e|^p \dd x
$$
and $v_\e= \overline{u}_\e= \e u_{A_\eta}(\frac{x}{\e})$ in a neighbourhood of $Q_\nu$. Scaling back to $TQ_\nu$, with $T=1/\e$,
$$
w_T(x):=T v_\e\left(\frac{x}{T}\right) \in \WW^{1,p}(TQ_\nu,\mathcal{Z}), 
$$
and $w_T=u_{A_\eta}$ in a neighbourhood of $\partial TQ_\nu$. Therefore
\begin{eqnarray*}
f^{p}_{\rm{hom}}(A_\eta)&=&\liminf_{T\to+\infty}\frac{1}{T^2}\inf\Bigl\{\int_{TQ_\nu}|\nabla u|^p\,\dd x: u\in \WW^{1,p}(TQ_\nu,
\Z),\
u= u_{A_\eta} \sim \partial TQ_\nu\Bigr\} \\
&\leq& \liminf_{T\to+\infty}\frac{1}{T^2} \int_{TQ_\nu} |\nabla w_T|^p\dd x= \liminf_{\e \to 0} \int_{Q_\nu}|\nabla v_\e|^p \dd x\\ &\leq& \liminf_{\e \to 0} \int_{Q_\nu}|\nabla u_\e|^p \dd x = \liminf_{T \to +\infty}\frac{1}{T^2}\int_{TQ_\nu} |\nabla u_\eta^T|^p \dd x\\
&\leq& (1+\delta)(|a_1|+|a_2|)^p +C_\delta  o_{\eta}(1),
\end{eqnarray*}
and  we find
$$
\limsup_{\eta \to 0}f^{p}_{\rm{hom}}(A_\eta) \leq (|a_1|+|a_2|)^p,
$$
by first sending $\eta \to 0$ and then $\delta \to 0$.
\end{proof}
\begin{proposition}
The $\Gamma$-limsup inequality for $F^p_\e$ with \emph{(\ref{gamma limit formula expression})} holds for all functions $u\in \WW^{1,p}(\Omega,\RR^2)$.    
\end{proposition}
\begin{proof}
By Proposition \ref{f hom is continuous} the function $A \mapsto f_{\rm hom}^p(A)$ is continuous, so that 
$$
u\mapsto \int_{\Omega}f_{\rm hom}^p(\nabla u(x)) \, \dd x
$$
is continuous in $\WW^{1,p}(\Omega,\RR^2)$. Fix $u \in \WW^{1,p}(\Omega,\RR^2)$; then there exists a sequence $u_k \in \WW^{1,p}(\Omega,\RR^2) $ piecewise affine, with non-singular gradients, such that $u_k \to u$ strongly as $k \to \infty$. By Proposition \ref{limsup for piecewise affine functions}, for each $k$ the upper bound holds
$$
\Gamma\hbox{-}\limsup_{\e\to 0} F^p_\e(u_k) \leq \int_{\Omega}f_{\rm hom}^p(\nabla u_k) \, \dd x.
$$
By the lower semicontinuity of the $\Gamma \hbox{-}\limsup$ we get
\begin{eqnarray*}
\Gamma\hbox{-}\limsup_{\e\to 0} F^p_\e(u)&\leq& \liminf_{k \to \infty}\Big(\Gamma\hbox{-}\limsup_{\e\to 0} F^p_\e(u_k) \Big)\\&\leq& \liminf_{k \to \infty}
\int_{\Omega}f_{\rm hom}^p(\nabla u_k) \, \dd x= \int_{\Omega}f_{\rm hom}^p(\nabla u) \, \dd x,   
\end{eqnarray*}
which gives the upper bound for the prescribed $u$.
\end{proof}

\subsection{The Young measure setting}\label{section YM}
Although the $\Gamma$-limit in Theorem A has already been proved without their use, we now introduce Young measures, as they make the geometric constraint of the problem transparent. Indeed, since $u_\e$ takes values in $\e\mathcal{Z}$, its gradient lies almost everywhere in
$$
K:=(e_1\otimes\mathbb{R}^2)\cup(e_2\otimes\mathbb{R}^2),
$$
every Young measure generated by such gradients is supported in $K$. This support condition is the starting point for the bounds on $f^p_{\mathrm{hom}}$ in Section~\ref{section bounds on f hom}, which in particular show that $f^p_{\mathrm{hom}}$ is not convex. 
In addition, Young measure give an alternative derivation of the rank-one formula $f^p_{\mathrm{hom}}(a\otimes\nu)=(|a_1|+|a_2|)^p$ to the construction of optimized taxicab paths in Section~\ref{subsection gamma liminf}. \bigskip

We write $\mathcal{P}(\mathbb{R}^N)$ for the set of Borel probability measures on $\mathbb{R}^N$.
The space $\mathbf{Y}^p(\Omega,\mathbb{R}^N)$ of $p$-\emph{Young measures} consists of all parametrized families $\nu=(\nu_x)_{x\in\Omega}\subset\mathcal{P}(\mathbb{R}^N)$ that are weakly$^*$ measurable: that is, $x\mapsto\int_{\RR^N} h\, \dd\nu_x$ is Lebesgue measurable for every $h\in \CC_0(\mathbb{R}^N)$, and satisfy
$$
  \int_\Omega\int_{\mathbb{R}^N}|z|^p\,\dd\nu_x(z)\,\dd x<\infty .
$$
A bounded sequence $(v_j)\subset \LL^p(\Omega,\mathbb{R}^N)$ \emph{generates} $\nu\in\mathbf{Y}^p(\Omega,\mathbb{R}^N)$, in symbols $v_j\overset{\mathbf{Y}}{\to}\nu$ as $j\to\infty$, if
$$
h(v_j)\rightharpoonup\int_{\RR^N} h\, \dd\nu_x\quad\text{weakly in }\LL^1(\Omega)\qquad\text{for every }h\in \CC_0(\mathbb{R}^N),
$$
that is, $\int_\Omega\varphi(x)\,h(v_j(x))\,\dd x\to\int_\Omega\varphi(x)\,\int_{\RR^N} h\,\dd\nu_x\,\dd x$ for all $\varphi\in \LL^\infty(\Omega)$ and $h\in \CC_0(\mathbb{R}^N)$.
We say that the Young measure $\nu$ is \emph{homogeneous} when there exists $\mu \in \mathcal{P}(\RR^N)$ such that $\nu_x=\mu$ for almost every $x \in \Omega$. A \emph{gradient} Young measure is generated by sequences of $\LL^p$ gradients $v_j=\nabla u_j$, for some sequence of Sobolev functions $(u_j)\subset \WW^{1,p}$. 

We record here two results from the general theory of Young measures, see \cite[Lem\-ma~4.12]{Rindlerbook} and \cite[Theorem~8.6]{FL07}; for a general Young-measure approach to lower semicontinuity, covering normal integrands and arbitrary tight (e.g.\ $\LL^p$-bounded) sequences, we refer the reader to Balder~\cite{Balder84}. 
\begin{lemma}\label{supp nu x in K lemma}
Let $\nu=(\nu_x)_{x \in \Omega} \in \mathbf{Y}^p(\Omega,\mathbb R^N)$ be the Young measure generated by a $\LL^p$-bounded sequence $(v_j) \subset \LL^p(\Omega,\mathbb R^N)$, and let $C \subset \mathbb R^N$ be a closed set. Then
$$
\emph{dist}(v_j,C) \to 0 \text{ in measure} \iff \emph{supp }\nu_x \subset C \emph{ for almost every }x \in \Omega.
$$
\end{lemma}
 We actually need the following weaker statement of the lemma above: if $v_j(x) \in C$ for almost every $x \in \Omega$ then $\text{supp }\nu_x \subset C$ for almost every $x \in \Omega$.
\begin{proposition}\label{YM LSC, normal integrands}
Let $1<p<\infty$ and $\nu=(\nu_x)_{x \in \Omega} \in \mathbf{Y}^p(\Omega,\mathbb R^N)$ be the Young measure generated by a $\LL^p$-bounded sequence $(v_j) \subset \LL^p(\Omega,\mathbb R^N)$. Let $F \colon \Omega \times \mathbb R^N \to [0,\infty)$ be a Carathéodory integrand of $p$-growth. Then  $$
(i) \quad \liminf_{j \to \infty} \int_\Omega F(x,v_j(x)) \, \dd x \geq \int_\Omega \int_{\mathbb R^N} F(x,M) \, \dd\nu_x(M) \, \dd x. 
$$
Moreover, $F(\cdot, v_j(\cdot))$ is equi--integrable if and only if
$$
(ii) \quad \lim_{j \to \infty} \int_\Omega F(x,v_j(x)) \, \dd x = \int_\Omega \int_{\mathbb R^N} F(x,M) \, \dd\nu_x(M) \, \dd x.
$$
\end{proposition} 
\bigskip 
Now we point out a pointwise differential constraint satisfied by the sequences $u_\e$ which are $\e \Z$-valued. 
Define the closed set $K \subset \mathbb R^{2\times 2}$ as 
\begin{equation}\label{definition of matrix set K}
K:= (e_1 \otimes \mathbb R^2) \cup (e_2 \otimes \mathbb R^2).
\end{equation}
So $M \in K$ takes the row-wise form
$
M=\begin{pmatrix}
M_1\\
0
\end{pmatrix} $ or $ M=\begin{pmatrix}
0\\
M_2
\end{pmatrix}$.

\begin{lemma}\label{nabla u_e has only one row}
If $u_\e \in \WW^{1,p}(\Omega,\e \mathcal{Z})$, then
$
\nabla u_\e(x) \in K$ for almost every $x \in \Omega$.
\end{lemma}
\begin{proof}
From $u_\e(x)=(u_{\e,1}(x),u_{\e,2}(x)) \in \e \mathcal{Z}$ for almost every $x \in \Omega$, 
we have that for almost every $x \in \Omega$
$$
u_{\e,1}(x) \in \frac{\e}{2}+\e\mathbb{Z} \ \text{ or } \ u_{\e,2}(x) \in \frac{\e}{2}+\e\mathbb{Z}.
$$
Since
$$
\frac{\e}{2}+\e\mathbb{Z}= \bigcup_{q \in \mathbb Z} \left\{ \frac{\e}{2}+\e q \right\}
$$
we can use the well-known property for Sobolev functions $f$, $1\leq p < \infty$, and  $c \in \mathbb R$, that
$
\nabla f=0 $ almost everywhere on $\{x: f(x)=c\}$, 
yielding that on the sets
$$
E^q_1:=\Big\{x \in \Omega: u_{\e,1}(x)= \frac{\e}{2}+\e q \Big\}, \quad E^q_2:=\Big\{x \in \Omega: u_{\e,2}(x)= \frac{\e}{2}+\e q \Big\} ,
$$
we have
$\nabla u_{\e,1}=0$ for almost every $x \in E_1^q$, and $\nabla u_{\e,2}=0$ for almost every $x \in E_2^q$.
Therefore, also $\nabla u_{\e,1}=0$ almost everywhere in $E_1=\bigcup_{q \in \mathbb Z} E_1^q$ and the same for $\nabla u_{\e,2}=0$ almost everywhere in $E_2=\bigcup_{q \in \mathbb Z} E_2^q$. Since $N:=\Omega \setminus (E_1 \cup E_2)$ is Lebesgue negligible by $u_\e(x) \in \e \mathcal{Z}$ for almost every $x$, we have that for almost every $x \in \Omega$ either
$$
\nabla u_{\e,1}=0 \implies \nabla u_\e(x)= \begin{pmatrix}
0\\
\nabla u_{\e,2}(x)
\end{pmatrix}, \quad \text{or}\quad
\nabla u_{\e,2}=0 \implies \nabla u_\e(x)= \begin{pmatrix}
\nabla u_{\e,1}(x)\\
0
\end{pmatrix},
$$
which gives $\nabla u_\e(x) \in K$.
\end{proof}

\noindent{\bf A lower bound 
in terms of Young measures.}
We fix $u \in \WW^{1,p}(\Omega,\mathbb R^2)$ and  $u_\e\rightharpoonup u$ weakly in $\WW^{1,p}(\Omega,\RR^2)$
with $\liminf_{\e\to0}F^p_\e(u_\e)<+\infty$. Upon extracting subsequences, we 
can suppose that this liminf is actually a limit, and that
$u_\e(x) \in \e \mathcal{Z}$  almost everywhere in $\Omega$.
Since in particular $F_\e^p(u_\e)=\int_\Omega |\nabla u_\e|^p\, \dd x$,  we can also assume that $\nabla u_\e$, as a $\LL^p(\Omega, \RR^{2 \times 2})$-bounded sequence, generates the Young measure $\nu=(\nu_x)_{x \in \Omega} \in \mathbf{Y}^p(\Omega,\RR^{2 \times 2})$. By condition $u_\e(x) \in \e \mathcal{Z}$ and Lemma \ref{nabla u_e has only one row} we infer $\nabla u_\e(x) \in K$ defined in \eqref{definition of matrix set K}  for almost every $x \in \Omega$. By Lemma \ref{supp nu x in K lemma} we find that
\begin{equation}\label{support of YM}
\text{supp }\nu_x \subset K \text{ for almost every }x \in \Omega.
\end{equation}
Since $p>1$,the center of mass $\overline{\nu}_x:=\int_{\RR^{2 \times 2}} M \ \dd\nu_x(M)$ equals $\nabla u(x)$ for almost every $x \in \Omega$. Indeed, for a given matrix $M$, the linear-growth function $\Phi(M)=M_{i,j}$ is such that $\Phi(\nabla u_\e)$ is $p$-uniformly integrable, since $p>1$, and therefore
$$
\int_\Omega \varphi(x) \int_{\mathbb R^{2 \times 2} } M_{i,j}\, \dd\nu_x(M) \, \dd x=  \lim_{\e \to 0}\int_{\Omega} \varphi(x) (\nabla u_\e)_{i,j}(x) \dd x = \int_{\Omega} \varphi(x) (\nabla u)_{i,j}(x) \dd x,
$$
for all $\varphi \in \CC(\overline{\Omega})$, where we have used the weak convergence $\nabla u_\e \rightharpoonup \nabla u$. Thanks to the properties above satisfied by $\nu$, we can employ Proposition \ref{YM LSC, normal integrands} (i) to find
\begin{multline*}
\liminf_{\e \to 0} F_\e^p(u_\e)= \liminf_{\e \to 0}\int_\Omega |\nabla u_\e|^p \, \dd x    \geq \int_\Omega \int_{\mathbb R^{2 \times 2}} |M|^p \, \dd\nu_x(M) \, \dd x\\
\geq \int_\Omega \inf \Big\{ \int_{\mathbb R^{2 \times 2}} |M|^p \, \dd\nu(M) \ : \ \nu \in \mathcal{P}(\RR^{2 \times 2}),\,  \text{supp }\nu \subset K, \, \overline{\nu}=\nabla u(x)   \Big\}\dd x.
\end{multline*}
Thus, a lower bound $f^p_{\mathcal{Y},\rm hom}\colon \RR^{2 \times 2} \to \RR$ ($\mathcal{Y}$ stands for Young measures) is defined as
$$
f^p_{\mathcal{Y},\rm hom}(A):= \inf \Big\{ \int_{\mathbb R^{2 \times 2}} |M|^p \, \dd\nu(M) \ : \ \nu \in \mathcal{P}(\RR^{2 \times 2}),\,  \text{supp }\nu \subset K, \, \overline{\nu}=A  \Big\}.
$$
Since $\text{supp } \nu \subset K= K_1 \cup K_2$, with $K_1= e_1 \otimes \RR^2$ and $K_2=e_2 \otimes \RR^2$, we can write $\nu$ as a convex combination of two probability measures
$$
\nu= \lambda \nu_1+ (1-\lambda)\nu_2, 
$$
where $\nu_1$ is supported on $K_1$ and $\nu_2$ is supported on $K_2$. Explicitly, $\nu_1$ is defined as $\nu_1=\frac{\nu \LLL K_1 }{\nu(K_1)}= \frac{\nu
\LLL K_1}{\lambda}$, and $\nu_2$ accordingly.
The moment condition $\overline{\nu}=A$ yields
$$
A=\lambda \overline{ \nu_1}+(1-\lambda)\overline{\nu}_2.
$$ Since $A= \begin{pmatrix}
A_1 \\ A_2    
\end{pmatrix}$ and $\nu_i$ are supported in $K_i$, which are convex subspaces of $\RR^{2\times 2}$, then also their center of masses $\overline{\nu}_i\in K_i$; hence, the moment condition row-wise reads
$$
\overline{\nu}_1=\frac{1}{\lambda}\begin{pmatrix}
A_1 \\ 0    
\end{pmatrix}= \frac{1}{\lambda} e_1 \otimes A_1 \quad \text{and} \quad \overline{\nu}_2=\frac{1}{1-\lambda}\begin{pmatrix}
0 \\ A_2    
\end{pmatrix}=\frac{1}{1-\lambda} e_2 \otimes A_2.
$$
By Jensen's inequality, applied to each $\nu_i \in \mathcal{P}(\RR^{2\times 2})$,
\begin{multline}
\int_{\mathbb R^{2 \times 2}} |M|^p \, \dd\nu(M)= \lambda \int_{\mathbb R^{2\times 2}}|M|^p \, \dd\nu_1(M)+(1-\lambda)\int_{\mathbb R^{2\times 2}}|M|^p \, \dd\nu_2(M) \\
\geq \lambda |\overline{\nu}_1|^p +(1-\lambda)|\overline{\nu}_2|^p = \frac{|A_1|^p}{\lambda^{p-1}}+ \frac{|A_2|^p}{(1-\lambda)^{p-1}}\geq (|A_1|+|A_2|)^p. \label{double jensen in YM section}
\end{multline}
The minimum $\lambda$ in the above expression is given by
$$
\lambda_\text{min}= \frac{|A_1|}{|A_1|+|A_2|}.
$$
The probability measure minimizing the problem defining $f^p_{\mathcal{Y},\rm hom}(A)$ is given by the convex combination of Dirac deltas centered on the center of masses $\overline{\nu}_1$ and $\overline{\nu}_2$, namely
$$
\nu_{\min}= \lambda \delta_{\frac{1}{\lambda} e_1 \otimes A_1}+(1-\lambda) \delta_{\frac{1}{1-\lambda}e_2 \otimes A_2},
$$
with $\lambda=\lambda_\text{min}$.
Indeed $\text{supp }\nu_{\min}= \{\frac{1}{\lambda} \, e_1 \otimes A_1\} \cup \{\frac{1}{1-\lambda} e_2 \otimes A_2\} \subset K$ and it holds
$$
\overline{\nu}_{\min}= \lambda \frac{1}{\lambda}e_1\otimes A_1+(1-\lambda)\frac{1}{1-\lambda}e_2\otimes A_2= A.
$$
The integral gives
$$
\int_{\RR^{2\times 2}} |M|^p \, \dd\nu_{\min}(M)=\lambda  \frac{|A_1|^p}{\lambda^{p}}+ (1-\lambda)\frac{|A_2|^p}{(1-\lambda)^{p}}=(|A_1|+|A_2|)^p.
$$
\begin{remark}
\emph{In the rank-one case $A= a \otimes n$, since $A_i=a_in$,
$$
(|A_1|+|A_2|)^p= (|a_1|+|a_2|)^p|n|^p=(|a_1|+|a_2|)^p,
$$
if $n$ is a unit vector. We thus recover the previous formula of $f^p_\text{hom}$.}
\end{remark}

The following proposition allows us to produce a recovery sequence for rank-1 matrices (see Corollary \ref{corollary Young measure limit rank one}).

\begin{proposition}\label{proposition recovery sequence rank one with YM}
Let $A=a\otimes n \in \RR^{2\times 2}$ with $a \in \RR^2$ and $n \in S^1$. Set $\lambda=\frac{|a_1|}{|a_1|+|a_2|}$  and $$B_1:= \frac{1}{\lambda}\begin{pmatrix}
a_1 n \\ 0    
\end{pmatrix}=\frac{a_1}{\lambda}e_1 \otimes n,
\qquad \qquad
B_2:= \frac{1}{1-\lambda}\begin{pmatrix}
0 \\ a_2 n    
\end{pmatrix}=\frac{a_2}{1-\lambda}e_2 \otimes n.
$$The homogeneous Young measure
$$
\nu_x= \lambda \delta_{B_1}+(1-\lambda) \delta_{B_2} \quad \text{for almost every }x \in \Omega
$$
is a gradient Young measure, with generating sequence $\nabla u_\e$ such that $u_\e(x) \in \e \mathcal{Z}$ almost everywhere and $\nabla u_\e$ is $p$--equi--integrable.
\end{proposition}

\begin{proof}
Given $a\in\RR^2$, let $a_k\in\mathbb Q^2$ with $a_k\to a$ as $k\to\infty$, and
for each $k$ let $T_k\in\NN$ be such that $T_ka_k\in\ZZ^2$. Set
$\beta_0:=(\frac12,\frac12)\in\RR^2$,
$z_k(t):=\beta_0+t\,a_k$, for $t\in\RR$,
and let $v_k\colon[0,T_k]\to\Z$ be the path given by Lemma~\ref{lem:staircase}
(with $\e=1$ and $L=T_k$) from $\beta_0$ to $\beta_0+T_ka_k$, which are vertices
of $\Z$. With
$$
\lambda_k:=\frac{|a_{k,1}|}{|a_{k,1}|+|a_{k,2}|},\qquad
v_{k,i}:={\rm sign}(a_{k,i})\,(|a_{k,1}|+|a_{k,2}|)\quad(i=1,2),
$$
we have $\lambda_kv_{k,1}=a_{k,1}$ and $(1-\lambda_k)v_{k,2}=a_{k,2}$, and, by
Lemma~\ref{lem:staircase}(i), there exist disjoint sets
$E^k_1,E^k_2\subset(0,T_k)$, with $|E^k_1|=\lambda_kT_k$ and
$|E^k_2|=(1-\lambda_k)T_k$, such that
$$
v'_k(t)=\begin{cases} v_{k,1}e_1, & t\in E^k_1\\ v_{k,2}e_2, & t\in E^k_2.\end{cases}
$$
The speeds $|v_{k,i}|=|a_{k,1}|+|a_{k,2}|$ converge to $|a_1|+|a_2|$, and
$|z_k(t)-v_k(t)|\le\sqrt2$ for all $t\in[0,T_k]$ by
Lemma~\ref{lem:staircase}(ii). Since $v_k(0)=z_k(0)$ and $v_k(T_k)=z_k(T_k)$,
the function $w_k:=v_k-z_k$ belongs to $\WW^{1,\infty}_0((0,T_k),\RR^2)$; we extend
it $T_k$-periodically to $\RR$ and set $v_k:=z_k+w_k$ on $\RR$. Then
$v_k(t+T_k)=v_k(t)+T_ka_k$, and, since $T_ka_k\in\ZZ^2$ and $\Z+\ZZ^2=\Z$, the
extended $v_k$ still takes values in $\Z$.
We claim that, for a suitable sequence $\e(k)\to0$ as $k\to\infty$, the maps
$u_{\e(k)}(x)=\e(k)\,v_k\big(n\cdot x/\e(k)\big)$ are the desired sequence.
For fixed $k$ we consider
$$
u_\e(x)= \e v_k\Big(\frac{n \cdot x}{\e}\Big), \quad x \in \Omega,
$$
and check the following properties.

\smallskip
1. Since  $v_k$ takes values in $\mathcal{Z}$, $u_\e(x) \in \e \mathcal{Z}$ for every $x \in \Omega$.

\smallskip 2.
    By the chain rule,
    $$
    \nabla u_\e (x) = v_k'\left( \frac{n\cdot x}{\e}\right)\otimes n
    $$
    and since $v_k$ is piecewise affine with two alternating slope values $v_{k,1}$ and $v_{k,2}$ which are uniformly close (up to a sign) to $|a_1|+|a_2|$, its derivative is uniformly bounded and therefore $\nabla u_\e$ is bounded as well; this fact gives the $p$--equi--integrability.

\smallskip 3. It remains to check the generation of the Young measure $\nu$. We note that    $v'_k(t)=w'_k(t)+ a_k$. 
    Since $w'_k$ is a periodic function as $w_k$, and $a_k$ is a constant, we deduce that $v'_k$ is also $T_k$-periodic. Therefore, by the Riemann-Lebesgue lemma, for a fixed measurable subset $E \subset \RR$ we have
    $$
    v'_k\left(\frac{t}{\e}\right) \rightharpoonup \frac{1}{T_k} \int_0^{T_k} v'_k(t) \dd t, \quad \text{in }\LL^p(E) \text{, as $\e \to 0$}. 
    $$
     Now take $\varphi \in \CC_0(\Omega)$ and $h \in \CC_0(\RR^{2 \times 2})$, with the notation $x=x_1n+x_2n^\perp$ and for each $x_2 \in \RR$, define
     $
     \Omega_{x_2}:=\{ x_1 \in \RR: x=x_1n+x_2n^\perp \in \Omega \}$.
     By Fubini's theorem and the Riemann--Lebesgue lemma in the variable $x_1$, we have
    \begin{align*}
    \lim_{\e\to 0}&\int_\Omega \varphi(x) h ( \nabla u_\e(x)) \dd x \\&=\lim_{\e\to 0} \int_\Omega \varphi(x) h \Big( v'_k\left(\frac{n \cdot x}{\e}\right) \otimes n    \Big) \dd x \\
    &=\lim_{\e\to 0}\int_{\RR}\int_{\Omega_{x_2}}\varphi(x_1,x_2) h \Big( v'_k(\frac{x_1}{\e})\otimes
     n\Big) \dd x_1 \, \dd x_2 \\ 
     &= \int_{\RR}\int_{\Omega_{x_2}}\varphi(x_1,x_2) \Big(\frac{1}{T_k}\int_0^{T_k} h \left( v'_k\left(t\right)\otimes
     n\right)\dd t\Big) \dd x_1 \, \dd x_2 \\&=\int_{\RR}\int_{\Omega_{x_2}}\varphi(x_1,x_2) \Big(\frac{|E^k_1|}{T_k} h(v_{k,1}e_1 \otimes n)+ \frac{|E^k_2|}{T_k} h(v_{k,2}e_2 \otimes n)\Big) \dd x_1 \,\dd x_2\\
     &=\int_{\RR}\int_{\Omega_{x_2}}\varphi(x_1,x_2) \Big(\lambda_k h\Big(\frac{a_{k,1}}{\lambda_k}e_1 \otimes n\Big)+ (1-\lambda_k) h\Big(\frac{a_{k,2}}{1-\lambda_k}e_2 \otimes n\Big)\Big) \dd x_1 \, \dd x_2\\
     &=\int_\Omega \varphi(x)\dd x\, \langle \lambda_k \delta_{B_1^k} +(1-\lambda_k)\delta_{B_2^k},h  \rangle = \int_\Omega \varphi(x)\dd x \, \langle \nu_k, h \rangle,
    \end{align*}
    where 
    $\nu_k=\lambda_k \delta_{B_1^k} +(1-\lambda_k)\delta_{B_2^k}$,
    $B_1^k=\frac{a_{k,1}}{\lambda_k}e_1 \otimes n$, and  $B_2^k= \frac{a_{k,2}}{1-\lambda_k}e_2 \otimes n$.

    Since $\nu_k\overset{*}{\rightharpoonup}\nu$ in $\CC_0(\RR^{2\times2})^*$ as $k\to\infty$, we can take a diagonal sequence $\e(k)\to 0$ as $k \to \infty$ such that
    $\nabla u_{\e(k)} \overset{\textbf{Y}}{\to} \nu $ as $k\to \infty$.
\end{proof}

\begin{corollary}\label{corollary Young measure limit rank one}
In the notation of Proposition {\rm \ref{proposition recovery sequence rank one with YM}}, the generating sequence $\nabla u_\e$ gives a recovery sequence for the linear function $u(x)=Ax$, namely
$$
\lim_{\e\to 0} F^p_\e(u_\e) = \lim_{\e \to 0} \int_\Omega |\nabla u_\e|^p \, \dd x =|\Omega| (|a_1|+|a_2|)^p=|\Omega|f_{\rm hom}^p(A).
$$
\end{corollary}
\begin{proof}
Since $\nabla u_\e$ is $p$--equi--integrable and generates $\nu_x= \lambda \delta_{B_1}+(1-\lambda) \delta_{B_2}$ with $B_1$ and $B_2$ as in Proposition \ref{proposition recovery sequence rank one with YM}, by Proposition \ref{YM LSC, normal integrands} (ii) applied to $F(M)=|M|^p$ we find
\begin{eqnarray*}
\lim_{\e \to 0} \int_\Omega |\nabla u_\e|^p \, \dd x = \int_\Omega \int_{\RR^{2\times 2}} |M|^p \dd\nu_x(M) \, \dd x &=& |\Omega|\left(\lambda \frac{|a_1|^p}{\lambda^p}+(1-\lambda)\frac{|a_2|^p}{(1-\lambda)^p} \right) \\ &=&|\Omega|(|a_1|+|a_2|)^p,
\end{eqnarray*}
and conclude the proof.
\end{proof}

\section{Quasiconvexity properties of  \texorpdfstring{$f^p_{\rm hom}$}{f\_hom}}\label{section bounds on f hom}

In this section we complete the proof of Theorem A, showing that the function $f^p_{\rm hom}$ is not convex. Since it satisfies a growth condition of order $p<2$ from above, this implies also that it is not polyconvex, since polyconvex non-convex functions satisfy grow at least quadratically in some direction. 
The goal of this section is then to find two constants $0<c_p\leq C_p$ such that
$$
(|A_1|+|A_2|)^p +c_p \frac{|\det A|}{|A|^{2-p}}\leq f^p_{\rm hom}(A) \leq (|A_1|+|A_2|)^p +C_p \frac{|\det A|}{|A|^{2-p}},
$$
with  $c_p, C_p \simeq (2-p)^{-1}$ as $p\to 2^-$.    
In particular, the lower bound shows the non-convexity of $f^p_{\rm hom}$, as remarked in the Introduction; see~\eqref{f hom not convex counterexample}.

\begin{lemma}[qualified Jensen estimate]\label{jensen with tail lemma}
Let $\nu \in \mathcal{P}(\RR^{2 \times 2})$ be a probability measure supported on $K
$, defined in \eqref{definition of matrix set K}, with barycenter $\overline{\nu}=A$. Let $\Lambda \geq \kappa_p |A|$, where $\kappa_p:=\sqrt{2}(2\sqrt{2}p)^{1/(p-1)}$. Then
$$
\int_{\RR^{2 \times 2}} |M|^p \dd\nu \geq (|A_1|+|A_2|)^p + \frac{1}{2}\int_{|M|> \Lambda}|M|^p \dd\nu.
$$
\end{lemma}
\begin{proof}
We recall the row-wise matrix norm $N(A)=(|A_1|+|A_2|)$ in \eqref{eq:N}, and the bound $N(A)\leq \sqrt{2}|A|$. We split $\nu$ at level $\Lambda$; that is,
$$
\nu=\nu|_{\{|M|\leq \Lambda\}}+\nu|_{\{|M|>\Lambda\}}=:\nu_b+ \nu_t   \quad \text{(bulk and tail)}.
$$
Let 
$$
T:=\int_{\mathbb R^{2 \times 2}} M \dd\nu_t \in \mathbb R^{2 \times 2} \quad \text{and} \quad E_p:= \int_{\mathbb R^{2 \times 2}} |M|^p \dd\nu_t= \int_{\{|M|>\Lambda\}} |M|^p \dd\nu
$$ be the first and $p$-th moments of the tail. Hence the bulk has first moment $\overline{\nu}_b=A-T$. 
Both $\nu_b$ and $\nu_t$ are supported in $K$; in particular, we can apply to $\nu_b$ the two-branch Jensen computation applied to $\nu$ above.

The first idea is to apply the double Jensen inequality as in \eqref{double jensen in YM section}. Split $\nu_b$ as a combination of two probability measures $$\nu_b= \lambda \nu_{b,1}+ (m_b-\lambda)\nu_{b,2},$$ 
with the mass $m_b:=|\nu_b|\leq 1$ and $\lambda \in [0,m_b]$. 
Then by the moment condition
$A-T=\overline{\nu}_b$ we get  $A-T=\lambda \overline{\nu}_{b,1}+(m_b-\lambda)\overline{\nu}_{b,2}$.
Since $\nu_{b,i}$ are supported in $e_i\otimes \RR^2$ subspaces of $\RR^{2\times 2}$, the moment condition row-wise reads
$$
\overline{\nu}_{b,1}= \frac{1}{\lambda} e_1 \otimes (A_1-T_1) \quad \text{and} \quad \overline{\nu}_{b,2}=\frac{1}{m_b-\lambda} e_2 \otimes (A_2-T_2).
$$
By Jensen's inequality applied to each probability measure $\nu_{b,i}$
\begin{eqnarray*}
 \int_{\{|M|\leq \Lambda \}}|M|^p \dd\nu&= &\int_{\mathbb R^{2 \times 2}}|M|^p \, \dd\nu_b(M)\\
 &=& \lambda \int_{\mathbb R^{2 \times 2}}|M|^p \, \dd\nu_{b,1}(M)+(m_b-\lambda)\int_{\mathbb R^{2 \times 2}}|M|^p \, \dd\nu_{b,2}(M)\\ &\geq& \lambda |\overline{\nu}_{b,1}|^p +(m_b-\lambda)|\overline{\nu}_{b,2}|^p = \frac{|A_1-T_1|^p}{\lambda^{p-1}}+ \frac{|A_2-T_2|^p}{(m_b-\lambda)^{p-1}}\\&\geq& m_b^{1-p} (|A_1-T_1|+|A_2-T_2|)^p=m_b^{1-p}N(A-T)^p \\
 &\geq& N(A-T)^p. 
\end{eqnarray*}
The minimum is given by $\lambda/m_b= |A_1-T_1|/(|A_1-T_1|+|A_2-T_2|)$. Hence, the bulk contributes with at least $N(A-T)^p$.

On the tail $|M| > \Lambda$, since $1 < |M|/\Lambda \leq |M|^p/\Lambda^p$ as $p>1$, it holds
\begin{eqnarray*}
|T|\leq \int_{\{|M|> \Lambda\}} |M| \dd\nu = \Lambda \int_{\{|M|> \Lambda\}} \frac{|M|}{\Lambda} \dd\nu \leq \Lambda \int_{\{|M|> \Lambda\}} \frac{|M|^p}{\Lambda^p} \dd\nu   =\Lambda^{1-p}E_p.
\end{eqnarray*}
We now want to estimate $N(A-T)^p$ with $N(A)^p$ plus other terms, using the convexity of $t\mapsto t^p$. Since $N$ is a norm, it holds $N(A-T)\ge |N(A)-N(T)|$. Applying the convexity inequality $t^p\ge s^p +ps^{p-1}(t-s)$ with $t=|N(A)-N(T)|\geq 0$ and $s=N(A)\geq 0$, we get
$$
N(A-T)^p\geq |N(A)-N(T)|^p \ge N(A)^p+pN(A)^{p-1}(|N(A)-N(T)|-N(A))
$$
Moreover we have $
(|N(A)-N(T)|-N(A)) \geq -N(T) \ge -\sqrt{2}|T| \ge -\sqrt{2}\Lambda^{1-p}E_p,
$ thus
$$
N(A-T)^p\geq  N(A)^p-\sqrt{2}pN(A)^{p-1}\Lambda^{1-p}E_p.
$$
Finally, the choice of $\Lambda \geq \kappa_p |A|$ with $\kappa_p$ as in the statement makes the second term on the right-hand side at most $E_p/2$, since by $N(A)\le \sqrt{2}|A|$ we get $-\sqrt{2}pN(A)^{p-1}\Lambda^{1-p}\geq -\frac{1}{2}$.  All the contributions now read
$$
\int_{\RR^{2\times 2}} |M|^p \dd\nu = \int_{\{|M|\leq \Lambda \}}|M|^p \dd\nu +E_p \geq N(A-T)^p +E_p \ge N(A)^p +\frac{E_p}{2},
$$
as desired, by the definition of $E_p$.
\end{proof}
For $f\in \LL^1_{\rm loc}(\RR^d)$ we denote by
$$
  \mathcal M f(x):=\sup_{r>0}\frac{1}{|B_r(x)|}\int_{B_r(x)}|f(y)|\,\dd y,\qquad x\in\RR^d,
$$
the Hardy--Littlewood maximal function of $f$. The function $\mathcal M f$ is lower semicontinuous, so its superlevel sets are open, and $\mathcal M$ is of weak type $(1,1)$ \cite[Chapter~I, Theorem~1]{Ste70}: there is $c=c(d)$ such that $|\{\mathcal M f>\lambda\}|\le\frac{c}{\lambda}\|f\|_{\LL^1(\RR^d)}$ for all $f\in \LL^1(\RR^d)$ and $\lambda>0$. Since $\mathcal M f\le\mathcal M(f\mathbf 1_{\{|f|>\lambda/2\}})+\lambda/2$, applying this to $f\mathbf 1_{\{|f|>\lambda/2\}}$ gives the refined estimate
\begin{equation}\label{eq:weak11}
  |\{\mathcal M f>\lambda\}|\le\frac{2c}{\lambda}\int_{\{|f|>\lambda/2\}}|f|\,\dd x\qquad\text{for all } f\in \LL^1(\RR^d),\ \lambda>0.
\end{equation}
We shall use the following Lipschitz truncation lemma, which goes back to Acerbi and Fusco \cite{AF,AF88}; see also \cite[Theorem~5.34]{AFP}.

\begin{lemma}[Lipschitz truncation]\label{lem:lipschitz-truncation}
Let $1\le p<\infty$, $w\in \WW^{1,p}(\RR^d,\RR^m)$, $\lambda>0$, and $E_\lambda:=\{\mathcal M(|\nabla w|)>\lambda\}$. Then there exists a Lipschitz map $\tilde w\colon\RR^d\to\RR^m$ such that
$\tilde w=w $ almost everywhere in   $\RR^d\setminus E_\lambda$, and $\mathrm{Lip}(\tilde w)\le c_d\,\lambda$,
where $c_d$ depends only on $d$; in particular, $c_d$ is independent of $p$, $\lambda$ and $w$.
If $w$ has compact support, then $E_\lambda$ is bounded and $\tilde w$ has compact support.
\end{lemma}

Using this result we prove the estimate in the following lemma.

\begin{lemma}[gradient estimates in terms of the Jacobian]\label{jacobians pay gradients} Let $D$ be bounded open and $u-Ax \in \WW^{1,p}_0(D,\RR^2)$, with $\det \nabla u=0$ almost everywhere. Then for all $\lambda \ge 4|A|$,
$$
\int_{|\nabla u| > \lambda/4 } |\nabla u|^p \dd x \geq c \, \frac{p-1}{2-p} \lambda^{p-2}|\det A| |D|,
$$
with $c$ universal constant.
\end{lemma}
\begin{proof}
Set $w:=u-Ax\in \WW^{1,p}_0(D,\RR^2)$, extended by zero to $\RR^2$, and let $E_\lambda:=\{\mathcal M(|\nabla w|)>\lambda\}$. By Lemma~\ref{lem:lipschitz-truncation} there exists a compactly supported Lipschitz map $\tilde w\colon\RR^2\to\RR^2$, in particular $\tilde w\in \WW^{1,\infty}(\RR^2,\RR^2)$, such that $\tilde w=w$ almost everywhere in $\RR^2\setminus E_\lambda$ and $\mathrm{Lip}(\tilde w)\le c_2\lambda$, where $c_2$ is independent of $p$. 

Since $\tilde w$ is Lipschitz with compact support, the null-Lagrangian
property of the determinant (valid for compactly supported maps in $\WW^{1,q}$, $q\ge 2$)
gives
$$
\int_{\RR^2 }(\det (A+ \nabla \tilde w)-\det A)\, \dd x =0.
$$
We split $\RR^2$ into $D\setminus E_\lambda$, $E_\lambda$ and
$\RR^2\setminus(D\cup E_\lambda)$. Since $\tilde w=w$ on
$\RR^2\setminus E_\lambda$, we have $\nabla\tilde w=\nabla w$ almost everywhere there.
Hence on $D\setminus E_\lambda$ the integrand equals
$\det\nabla u-\det A=-\det A$, because $\det\nabla u=0$ almost everywhere; on
$\RR^2\setminus(D\cup E_\lambda)$ we have $\tilde w=w=0$, so the integrand
vanishes. Therefore
$$
|D \setminus E_\lambda |\det A = \int_{E_\lambda} (\det (A+ \nabla \tilde w)-\det A)\, \dd x.
$$
On $E_\lambda$ we use the local Lipschitz estimate on the determinant, which is a polynomial of degree $2$, $ |\det(\xi_1)-\det(\xi_2)|\leq (|\xi_1|+|\xi_2| )|\xi_1-\xi_2|, $
which gives
$$
|\det (A+ \nabla \tilde w)-\det A| \leq (|A+ \nabla \tilde w|+|A|)|\nabla \tilde w|\leq C\lambda^2, 
$$
as $|\nabla \tilde w|\leq c_2 \lambda$ and $|A|\leq \lambda/4$. The null-Lagrangian identity now reads
$$
|D \setminus E_\lambda | |\det A| \leq C\lambda^2 |E_\lambda|.
$$
Writing $|D|=|D \setminus E_\lambda|+ |D \cap E_\lambda| \leq |D \setminus E_\lambda| +|E_\lambda|$, we arrive at the estimate
$$
|D||\det A| \leq |D \setminus E_\lambda | |\det A| + |E_\lambda||\det A| \leq (C\lambda^2+ |\det A|) |E_\lambda| \leq C' \lambda^2 |E_\lambda|,
$$
using $|\det A|\leq \frac{1}{2}|A|^2 \leq \frac{1}{2}\frac{\lambda^2}{16}$. 

The weak $(1,1)$ type estimate \eqref{eq:weak11} gives
$$
|E_\lambda| \leq \frac{C}{\lambda} \int_{\{|\nabla w|> \lambda/2\}} |\nabla w| \dd x,
$$
so that 
$$
|D||\det A| \leq C \lambda \int_{\{|\nabla w|> \lambda/2\}} |\nabla w| \dd x.
$$
Fix $\lambda_0\ge4|A|$. The last inequality holds for every $\lambda\ge\lambda_0$,
since then $\lambda\ge4|A|$. Multiplying it by $\lambda^{p-3}$, integrating over
$[\lambda_0,\infty)$ and using Fubini's theorem, we get
\begin{eqnarray*}
\frac{\lambda_0^{p-2}}{2-p}|D| |\det A|&\leq& C \int_{\lambda_0}^\infty \lambda^{p-2} \int_{\{|\nabla w|> \lambda/2\}} |\nabla w| \dd x \, \dd\lambda\\
&=&C \int_{\{|\nabla w|> \lambda_0/2\}} |\nabla w|  \bigg(\int_{\lambda_0}^{2|\nabla w|} \lambda^{p-2} \, \dd\lambda \bigg)\dd x\\
&\leq& C \int_{\{|\nabla w|> \lambda_0/2\}} |\nabla w|  \bigg(\int_{0}^{2|\nabla w|} \lambda^{p-2} \, \dd\lambda \bigg)\dd x\\
&=& C\frac{2^{p-1}}{p-1}\int_{\{|\nabla w|> \lambda_0/2\}} |\nabla w|^p \, \dd x.
\end{eqnarray*}
To convert the bounds on $u$ rather than on $w$, first note that on $\{|\nabla w|> \lambda_0/2\}$ one has
$$
\frac{\lambda_0}{2} < |\nabla w|=|\nabla u-A|\leq |\nabla u|+|A|\leq |\nabla u| +\frac{\lambda_0}{4},
$$
so that $|\nabla u|> \lambda_0/4$. Moreover, by the same line above, $|\nabla w| \leq |\nabla u|+\frac{\lambda_0}{4}\leq |\nabla u|+\frac{|\nabla w|}{2}$, so that $|\nabla w|\leq 2|\nabla u|$. Therefore, we find
$$
\frac{\lambda_0^{p-2}}{2-p}|D| |\det A| \leq C'\frac{2^{p-1}}{p-1}\int_{|\nabla u| > \lambda_0/4 } |\nabla u|^p \dd x,
$$
which we rewrite as
$$
\int_{|\nabla u| > \lambda_0/4 } |\nabla u|^p \dd x \geq c \, \frac{p-1}{2-p} \lambda_0^{p-2}|\det A| |D|;
$$
that is, the claim with $\lambda=\lambda_0$ and $c=1/C'$ universal constant independent of $p$.
\end{proof}

\begin{proposition}[lower bound on $f^p_{\rm hom}$]\label{lower bound on f hom with N(A)}
    There exists a constant $c_p>0$ such that for every $A$
    $$
    f^p_{\rm hom}(A) \geq N(A)^p +c_p \frac{|\det A|}{|A|^{2-p}},
    $$
    where $N$ is the norm in \eqref{eq:N}.
    The constant $c_p \simeq (2-p)^{-1}$ as $p\to 2^-$.
\end{proposition}

\begin{proof} If $\det A=0$ the inequality is trivial.
In the full-rank case, recall the cell formula \eqref{formula f_hom with A non singular}
$$
f^p_{\rm hom}(A) = \lim_{T\to+\infty}\frac{|\det A|}{T^2}\inf\Bigl\{\int_{TQ_A}|\nabla u|^p\,\dd x: u\in \WW^{1,p}(TQ_A,
\mathcal Z),\
u= u_A \sim \partial TQ_A\Bigr\}.
$$
Fix $T$ an odd integer, so that $\partial(TQ)\subset\mathcal Z$ and hence
$P_0(Ax)=Ax$ on $\partial D$, $D=A^{-1}(TQ)$. Let $u$ be a competitor,
so that $u-Ax\in \WW^{1,p}_0(D,\RR^2)$, and let $\nu_u$ be the gradient distribution,  explicitly
$$
\nu_u:=(\nabla u)_\#\Big(\tfrac1{|D|}\,\mathcal L^2\LLL D\Big),
\qquad
\int_{\mathbb R^{2 \times 2}} g(M)\,\dd\nu_u=\fint_D g(\nabla u)\,\dd x
$$
for every Borel $g\ge0$. Then $\nu_u$ is a probability measure on
$\RR^{2\times2}$ with $\nu_u(\RR^{2\times2}\setminus K)=0$, since
$\nabla u\in K$ almost everywhere by Lemma \ref{nabla u_e has only one row}. Its first moment is finite,
and
$$
\int_{\mathbb R^{2 \times 2}}M\,\dd\nu_u=\fint_D\nabla u\,\dd x=A,
$$
because $\int_D\nabla(u-Ax)\,\dd x=0$ for $u-Ax\in \WW^{1,p}_0(D)$. Finally,
$\det\nabla u=0$ almost everywhere in $D$, since every matrix in $K$ has a zero row.
Note that
$$
\int_{\{|M|>\Lambda\}}|M|^p\,\dd\nu_u
=\fint_D|\nabla u|^p\,\chi_{\{|\nabla u|>\Lambda\}}\,\dd x ,
$$
so Lemma \ref{jensen with tail lemma}, applied to $\nu_u$, and
Lemma \ref{jacobians pay gradients}, applied to $u$, estimate the same quantity.
 We choose
$$
\Lambda:=\kappa_p|A|,\qquad \kappa_p:=\sqrt2\,(2\sqrt2\,p)^{\frac1{p-1}} .
$$
As shown above, $\nu_u$ is a
probability measure concentrated on $K$, with finite first moment and
barycenter $A$, and $\Lambda\ge\kappa_p|A|$. Hence
\begin{eqnarray*}
\fint_D|\nabla u|^p\,\dd x=\int|M|^p\,\dd\nu_u
&\ge& N(A)^p+\frac12\int_{\{|M|>\Lambda\}}|M|^p\,\dd\nu_u\\
&=& N(A)^p+\frac12\fint_D|\nabla u|^p\,\chi_{\{|\nabla u|>\Lambda\}}\,\dd x .    
\end{eqnarray*}
Recall 
$\det\nabla u=0$ almost everywhere in $D$ and $u-Ax\in \WW^{1,p}_0(D,\RR^2)$. We apply 
Lemma \ref{jacobians pay gradients} with $\lambda:=4\Lambda=4\kappa_p|A|$, which is admissible since
$\lambda\ge4|A|$, and obtain
$$
\fint_D|\nabla u|^p\,\chi_{\{|\nabla u|>\Lambda\}}\,\dd x
=\frac1{|D|}\int_{\{|\nabla u|>\lambda/4\}}|\nabla u|^p\,\dd x
\ \ge\ \frac{c\,(p-1)}{2-p}\,\big(4\kappa_p|A|\big)^{p-2}\,|\det A| .
$$
Combining the two estimates,
$$
\fint_D|\nabla u|^p\,\dd x\ \ge\ N(A)^p+c_p\,\frac{|\det A|}{|A|^{2-p}},
\qquad
c_p:=\frac{c\,(p-1)\,(4\kappa_p)^{p-2}}{2\,(2-p)} .
$$
The right-hand side depends neither on $u$ nor on $T$. Taking the infimum
over all competitors and letting $T\to\infty$ along odd integers, we obtain
$$
f^p_{\rm hom}(A)\ \ge\ N(A)^p+c_p\,\frac{|\det A|}{|A|^{2-p}},
$$
and conclude the proof.
\end{proof}

\begin{proposition}[upper bound on $f^p_{\rm hom}$]\label{upper bound of f hom with N(A)}
    There exists a constant $C_p>0$ such that for every $A$
    $$
    f^p_{\rm hom}(A) \leq N(A)^p +C_p \frac{|\det A|}{|A|^{2-p}}.
    $$
    The constant $C_p \simeq (2-p)^{-1}$ as $p\to 2^-$.
\end{proposition}

\begin{proof}
Decompose $A= \sigma_1 \eta_1 \otimes q_1 + \sigma_2 \eta_2 \otimes q_2$ with respect to its singular values $\sigma_i$, with $0 < \sigma_2 \leq \sigma_1$, $\eta_1 \perp \eta_2$ and $q_1 \perp q_2$ orthonormal. The idea of this decomposition into orthogonal rank-one tensors is to use the rank-one convexity of $f^p_{\rm hom}$ to find an upper bound. Write 
$$
R:=\eta_1 \otimes q_1 + \eta_2 \otimes q_2 \in O(2), \quad A= \theta(\sigma_1 R)+ (1-\theta)(\sigma_1 \eta_1\otimes q_1),\qquad \theta:=\frac{\sigma_2}{\sigma_1}\in [0,1].
$$
The two matrices $\sigma_1 R$ and $\sigma_1 \eta_1\otimes q_1$ differ by the rank-one matrix $\sigma_1 \eta_2\otimes q_2$, so $A$ is written as a convex combination along a rank-one direction. Moreover
$$
f^p_{\rm hom}(\sigma_1 R)=\sigma_1^p f^p_{\rm hom}(I)
$$
by the $p$-homogeneity (Lemma \ref{lem:homog}) and the frame indifference under $A\mapsto AQ$, $Q\in O(2)$.
Using the rank-one convexity and $f^p_{\rm hom}=N(\cdot)^p$ on the rank-one cone,
\begin{equation}
    \label{eq:laminate}
f^p_{\rm hom}(A)\leq (1-\theta)(\sigma_1|\eta_1|_{\ell^1})^p+ \theta \sigma_1^p f^p_{\rm hom}(I).
\end{equation} 
Denote for simplicity $A_0:=\sigma_1 \eta_1 \otimes q_1$ so that $N(A_0)=\sigma_1\,|\eta_1|_{\ell^1}$.\\
Since
$\theta\,\sigma_1^p=\sigma_2\,\sigma_1^{p-1}=\sigma_1\sigma_2\,\sigma_1^{p-2}
=|\det A|\,\sigma_1^{p-2}$, the right-hand side of \eqref{eq:laminate} equals
$$
N(A_0)^p+\theta\,\sigma_1^p\big(f^p_{\rm hom}(I)-|\eta_1|_{\ell^1}^p\big)
=N(A_0)^p+\frac{|\det A|}{\sigma_1^{2-p}}\big(f^p_{\rm hom}(I)-|\eta_1|_{\ell^1}^p\big).
$$
The bracket is nonnegative, since
$f^p_{\rm hom}(I)\ge N(I)^p=2^p\ge2^{p/2}\ge|\eta_1|_{\ell^1}^p$, and it is at most
$f^p_{\rm hom}(I)-1$, since $|\eta_1|_{\ell^1}\ge|\eta_1|=1$. We now want to bound $N(A_0)^p$ with $N(A)^p$ plus a term proportional to $|\det A||A|^{p-2}$.  First of all, using that $N$ is a norm and
$A-A_0=\sigma_2\,\eta_2\otimes q_2$, 
$$
|N(A_0)-N(A)|\le N(\sigma_2\,\eta_2\otimes q_2)=\sigma_2\,|\eta_2|_{\ell^1}
\le\sqrt2\,\sigma_2 ,
$$
because $|\eta_2|_{\ell^1}\le\sqrt2\,|\eta_2|=\sqrt2$. For $s,h\ge0$, the mean value theorem gives
$(s+h)^p=s^p+p\,\xi^{p-1}h \le s^p+p\,(s+h)^{p-1}h$. With $s=N(A)$ and $h=\sqrt2\,\sigma_2$,
$$
N(A_0)^p\le\big(N(A)+\sqrt2\,\sigma_2\big)^p
\le N(A)^p+\sqrt2\,p\,\big(N(A)+\sqrt2\,\sigma_2\big)^{p-1}\sigma_2 .
$$
Since $|A|^2=\sigma_1^2+\sigma_2^2$, we have $\sigma_2\le|A|$ and
$\sigma_1\ge|A|/\sqrt2$; moreover $N(A)\le\sqrt2\,|A|$. Hence
$N(A)+\sqrt2\,\sigma_2\le2\sqrt2\,|A|$ and
$\sigma_2=|\det A|/\sigma_1\le\sqrt2\,|\det A|/|A|$, so that
$$
N(A_0)^p\le N(A)^p+2p\,(2\sqrt2)^{p-1}\,\frac{|\det A|}{|A|^{2-p}} .
$$

To conclude, since $p-2<0$ and $\sigma_1\ge|A|/\sqrt2$, we have
$\sigma_1^{p-2}\le2^{\frac{2-p}{2}}|A|^{p-2}$. Therefore
$$
f^p_{\rm hom}(A)\le N(A)^p+C_p\,\frac{|\det A|}{|A|^{2-p}},
\qquad
C_p:=2p\,(2\sqrt2)^{p-1}+2^{\frac{2-p}{2}}\big(f^p_{\rm hom}(I)-1\big).
$$
The value of $f^p_{\rm hom}(I)$ is bounded by $\frac{8}{3(2-p)}$ by \eqref{f on the ID upper bound}. \end{proof}

Finally, we remark that from the explicit form of $c_p$ and $C_p$, a direct computation allows one to rewrite the double estimate as in \eqref{eq:bounds}.

\section{The supercritical case \texorpdfstring{$p\ge 2$}{p ≥ 2}}\label{section p>2}
In this section, we consider the case $p\ge 2$, in which the limit is represented as a simple but degenerate quasiconvex energy. For convenience, we recall the statement of the homogenization result, already given in the Introduction as Theorem B.

\begin{theorem}[supercritical homogenization theorem]\label{thm:super:hom}
Let $p\ge2$. Then $F^p_\e$ given in \eqref{eq1} $\Gamma$-converges, with respect to the strong
$\LL^p$-topology, to
\begin{equation}\label{fompp}
F^p_{\rm hom}(u)=
\begin{cases}
\displaystyle\int_\Omega\bigl(|\nabla u_1(x)|+|\nabla u_2(x)|\bigr)^p\,\dd x
&\hbox{if $u\in \WW^{1,p}(\Omega,\RR^2)$ and $\det\nabla u=0$ a.e.}\cr
+\infty&\hbox{otherwise.}
\end{cases}
\end{equation}\end{theorem}

The lower bound is obtained by noting that actually $F^p_\e= F^p_{\rm hom}$ on $\WW^{1,p}(\Omega,\e\Z)$.

\begin{proposition}[lower bound]\label{lower bound}
The functional defined by
\eqref{fompp} is a lower bound for the $\Gamma$-liminf of $F_\e$.
\end{proposition}

\begin{proof} We first note that the $\Gamma$-limit is finite only on functions with $\det\nabla u(x)=0$ almost everywhere, by the weak continuity of the determinant in $\WW^{1,p}$, since   $\det\nabla u_\e(x)=0$ almost everywhere if $F^p_\e(u_\e)<+\infty$, and using the continuity of the determinant for $p\ge 2$.

\smallskip
Let $u_\e\to u $ be a sequence with $F^p_\e(u_\e)\le C<+\infty$. 
We note that if $u_\e(x)\not\in \frac{\e}{2}+\e\ZZ^2$ then for $|x-y|$ small we have either $(u_\e)_1(x)=(u_\e)_1(y)$ or  $(u_\e)_2(x)=(u_\e)_2(y)$, so that we conclude that
$$
|\nabla u_\e(x)|=|\nabla (u_\e)_1(x)|+|\nabla (u_\e)_2(x)|
$$
for almost all $x\in\Omega$.
By the convexity of the function $g(A)=(|A_1|+|A_2|)^p$, where $A_1$ and $A_2$ are the rows of $A$,
we then have
\begin{eqnarray*}
  \liminf_{\e \to 0} \int_\Omega |\nabla u_\e|^p \,\dd x &=& \liminf_{\e \to 0} \int_\Omega g(\nabla u_\e)\, \dd x \\
  &\geq& \int_\Omega g(\nabla u) \, \dd x= \int_\Omega ||\nabla u_1(x)|+|\nabla u_2(x)||^p\,\dd x,   
\end{eqnarray*}
and conclude the proof.
\end{proof}

\subsection{Recovery sequences for a class of special maps}\label{subsec:treesp}

In this section, we complete the proof of Theorem \ref{thm:super:hom} for functions whose image is a curve. The proof generalizes the zig-zag construction for affine functions in the case $p<2$. In the general case, the image is more complex, but for $p>2$, thanks to the continuity of $u$, it can be described using some notions of analytic topology (next section). These allow us to consider a dense class of functions whose images are ``trees'', on whose edges one can repeat the constructions below (or directly use Lemma \ref{lem:staircase}).

Note that if $\det\nabla u=0$, the two vectors $\nabla u_1(x)$ and $\nabla u_2(x)$ are linearly dependent, so that
\begin{equation}
    ||\nabla u_1(x)|+|\nabla u_2(x)||=\max\Big\{|\nabla u_1(x)+\nabla u_2(x)|,|\nabla u_1(x)-\nabla u_2(x)|
    \Big\}.
\end{equation}

\begin{proposition}[upper bound for special functions--I]\label{upper bound}
Let $u:\Omega\to \RR^2$ be such that $u(\Omega)$ is a simple curve. Then there exists $u_\e\to u$
such that
\begin{equation}\label{recove}
\lim_{\e\to 0} F^p_\e(u_\e)=\int_\Omega ||\nabla u_1(x)|+|\nabla u_2(x)||^p\,\dd x.
\end{equation}
\end{proposition}

\begin{proof}
By approximation, we may consider the case when $u(\Omega)$ is composed 
of line segments. Furthermore, we may suppose that the direction $\tau$ of each such
segment is rational (that is, $\tau\in\RR^2$ with either $\tau_1=0$ or $\tau_2/\tau_1$ 
rational). 

We first consider the case $u(\Omega)=[0,\alpha]=:S$ (where we
have used the notation $[\alpha_1,\alpha_2]$ to
denote the segment between the vectors $\alpha_1$ and $\alpha_2$ in $\RR^2$)
and suppose $\alpha_i=\e N\tau_i$ with 
$N\in\NN$ and $\tau_i\in\ZZ$. This again is not restrictive up to an error of order $\e$.

We may approximate $S$ by the family of line-segment curves $S_\e$ defined by
$$
S_\e:= \bigcup_{k=0}^{N-1} \Bigl(k\e\tau+\Bigl([0,\e\tau_1]\times\{0\}\Bigr)\cup
\Bigl(\{\e\tau_1\}\times[0,\e\tau_2]\Bigr)\Bigr)\,.
$$
Note that the extreme points of $S_\e$ are the endpoints of $S$.

If $\tau_1\tau_2\ge0$ (so that we may suppose that both of them are non-negative)
we consider the bijection $\Phi_\e:S\to S_\e$ defined by
$$
\Phi_\e(\xi)=\begin{cases}k\e\tau+(\xi_1+\xi_2-k\e(\tau_1+\tau_2), 0) &\hbox{ if $\xi\in k\e\tau+ \Bigl[0,\frac{\tau_1}{\tau_1+\tau_2}\e\tau\Bigr]$}\cr
k\e\tau+(\e\tau_1, \xi_1+\xi_2-\e((k+1)\tau_1+k\tau_2) &\hbox{ if $\xi\in k\e\tau+ \Bigl[\frac{\tau_1}{\tau_1+\tau_2}\e\tau, \e\tau\Bigr]$,} 
\end{cases}
$$
and $k\in\{0,\ldots,N-1\}$.

We then define 
$u_\e(x)=\Phi_\e(u(x))$,
and note that 
$$
\Big|\frac{\partial u_\e}{\partial x_i}\Big|=\Big|\frac{\partial u_1}{\partial x_i}+\frac{\partial u_2}{\partial x_i}\Big|=\Big|\frac{\partial u_1}{\partial x_i}\Big|+\Big|\frac{\partial u_2}{\partial x_i}\Big|,
$$
the last equality due to condition $\tau_2\frac{\partial u_1}{\partial x_i}=\tau_1\frac{\partial u_2}{\partial x_i}$ which is valid almost everywhere.
This implies 
\begin{eqnarray*}
|\nabla u_\e(x)|^2&=&
\Big|\frac{\partial u_\e}{\partial x_1}\Big|^2+\Big|\frac{\partial u_\e}{\partial x_2}\Big|^2
=
\Big|\frac{\partial u_1}{\partial x_1}+\frac{\partial u_2}{\partial x_1}\Big|^2+\Big|\frac{\partial u_1}{\partial x_2}+\frac{\partial u_2}{\partial x_2}\Big|^2
\\
&=&
\Big(\frac{\tau_1+\tau_2}{\tau_1}\Big)^2
\Big(\Big|\frac{\partial u_1}{\partial x_1}\Big|^2+\Big|\frac{\partial u_1}{\partial x_2}\Big|^2\Big)\ =\ 
\Big(\frac{\tau_1+\tau_2}{\tau_1}\Big)^2|\nabla u_1|^2.   
\end{eqnarray*}

Hence, reasoning symmetrically for $u_2$,
$$
|\nabla u_\e|=\frac{\tau_1+\tau_2}{\tau_1}|\nabla u_1|,\qquad
|\nabla u_\e|=\frac{\tau_1+\tau_2}{\tau_2}|\nabla u_2|.
$$
We then have
$$
|\nabla u_\e|=\frac{\tau_1}{\tau_1+\tau_2}|\nabla u_\e|+\frac{\tau_2}{\tau_1+\tau_2} |\nabla u_\e|= |\nabla u_1|+ |\nabla u_2|.
$$

If $\tau_1\tau_2<0$ then the bijection is defined analogously with $\xi_1-\xi_2$ in 
place of $\xi_1+\xi_2$, and the same argument gives
$$
\Big|\frac{\partial u_\e}{\partial x_i}\Big|=\Big|\frac{\partial u_1}{\partial x_i}-\frac{\partial u_2}{\partial x_i}\Big|=\Big|\frac{\partial u_1}{\partial x_i}\Big|+\Big|\frac{\partial u_2}{\partial x_i}\Big|,
$$
from which  in all cases we derive
$
|\nabla u_\e(x)|^p=||\nabla u_1(x)|+|\nabla u_2(x)||^p$,
and we have the desired equality in the limit.

If $u(\Omega)$ is composed of more line segments, we can repeat the construction separately on each line segment.    
\end{proof}

\begin{proposition}[upper bound for special functions--II]\label{upper bound2}
Let $u:\Omega\to \RR^2$ be a piece\-wise-affine function with {\rm det}$\nabla u=0$ almost everywhere. 
Then there exist $u_\e\to u$ such that {\rm(\ref{recove})} holds.
\end{proposition}

\begin{proof}Up to an approximation argument, we may suppose that $u$ takes values 
in $\e\Z$ at the vertices of each triangle $T$ of its underlying triangulation. The condition 
{\rm det}$\nabla u=0$ implies that $u(T)$ is a segment for all such $T$.
Again by an approximation argument, we may also suppose that the direction
of each such segment is rational and that we may 
separately apply the construction in the first part of the proof of Proposition \ref{upper bound} 
to the restriction of $u$ to (the interior of) each $T$. 
We note now that if $T$ and $T'$ have as common side a segment 
with endpoints $P$ and $P'$, then we have two cases:
$u({P})\neq u(P')$, so that the image of $T$ and that of $T'$ lie in the same line.
or $u({P})= u(P')$. In both cases, the constructions of Proposition \ref{upper bound} 
on $T$ and $T'$ agree on the segment $[P,P']$; in the second case, simply because
$u_\e= u({P})= u(P')$ on the whole segment. 
\end{proof}

\subsection{Recovery sequences for general maps}\label{subsec:trees}

 Throughout this subsection $\Omega\subset\RR^2$ is a bounded convex open set, $p>2$, and $u\in \WW^{1,p}(\Omega,\RR^2)$ satisfies
$\det\nabla u=0$ almost everywhere; by Morrey's embedding we identify $u$ with its
continuous extension to $\overline{\Omega}$. We write $N(A)=|A_1|+|A_2|$, and for a continuous
curve $\gamma$ we let ${\rm Var}_1(\gamma)$ denote its length with respect to
$|\cdot|_1$.

The following result, together with 
Proposition~\ref{lower bound},
 proves Theorem \ref{thm:super:hom} in the case $p>2$ and $\Omega$ convex.
 
\begin{theorem}\label{thm:limsup-p>2} Let $p>2$, and let $u\in \WW^{1,p}(\Omega,\RR^2)$ satisfy
$\det\nabla u=0$ almost everywhere. Then 
$$\Gamma\hbox{-}\limsup_{\e\to0}F^p_\e(u)\le\int_\Omega N(\nabla u)^p\,\dd x,
$$
with respect to the strong $\LL^p$-topology.
\end{theorem}

\begin{lemma}\label{lem:null-image} Let $p>2$, and let $u\in \WW^{1,p}(\Omega,\RR^2)$ satisfy
$\det\nabla u=0$ almost everywhere. Then 
$|u(\overline \Omega)|=0$.
\end{lemma}

\begin{proof}
Let $\tilde u\in \WW^{1,p}(\RR^2,\RR^2)$ be an extension of $u$. Since $p>2$,
$\tilde u$ is continuous and satisfies Lusin's condition (N) \cite{MM}, so
that the area formula holds. Hence
$$|u(\Omega)|\le
\int_{\RR^2}\#\bigl(u^{-1}(y)\cap\Omega\bigr)\,\dd y=\int_\Omega|\det\nabla u|\,\dd x=0.
$$
Moreover, $|\partial\Omega|=0$, so that
$|u(\partial\Omega)|=|\tilde u(\partial\Omega)|=0$ by (N).
\end{proof}

\begin{lemma}[factorization through a dendrite]\label{lem:dendrite}
Let $X$ be the set of connected components of the fibers $u^{-1}(y)$,
$y\in u(\overline{\Omega})$, endowed with the quotient topology; let $q:\overline{\Omega}\to X$ map $x$ to the
component containing it and set $\psi(q(x)):=u(x)$. Then $u=\psi\circ q$, and:

\smallskip{\rm (i)} $q$ is continuous, closed, surjective, and monotone; that is, its
fibers are connected;

\smallskip{\rm (ii)} $\psi$ is continuous and light; that is, its fibers are totally
disconnected;

\smallskip{\rm (iii)} $X$ is a {\em dendrite}; that is, a locally connected continuum containing
no simple closed curve.
\end{lemma}

\begin{proof}
Items (i) and (ii) are standard, and can be directly derived from the monotone-light factorization \cite[Theorem 4.2]{Whyburn}. 

It remains to prove item (iii). Since $\overline{\Omega}$
is a locally connected continuum, so is its continuous image $X$. By Lemma~\ref{lem:null-image} the compact set $u(\overline{\Omega})$ has
an empty interior, hence $\dim u(\overline{\Omega})\le1$. Since $\psi:X\to u(\overline{\Omega})$ is a closed
surjection with zero-dimensional fibers, the Hurewicz formula for closed maps
gives $\dim X\le\dim u(\overline{\Omega})\le1$ \cite[Theorem VI.7]{HW}.

 We are left to prove that $X$ contains no simple closed curve. Let $g:X\to S^1$ be continuous. Since $\overline{\Omega}$ is convex,
$g\circ q=e^{ih}$ for some continuous $h:\overline{\Omega}\to\RR$. On each fiber $q^{-1}(a)$,
which is connected, $g\circ q$ is constant; hence $h$ takes values in a
discrete set there, and then it is constant on the fiber. Thus, $h=\bar h\circ q$, with $\bar h$
continuous because $q$ is a quotient map and $g=e^{i\bar h}$ is
null-homotopic.
If $S\subset X$ were a simple closed curve,
then, since $\dim X\le1$, a homeomorphism $S\to S^1$ would extend to a
continuous map $X\to S^1$ \cite{HW}. Since this extension, and
hence the homeomorphism, would be null-homotopic, this is a contradiction.
\end{proof}

The following result will be used to prove the H\"older continuity of $d_1$.
\begin{lemma}
\label{lem:rays}
Let $p>2$ and $v\in \WW^{1,p}(\RR^2,\RR^2)$, identified with its continuous
representative. Then for every $x\in\RR^2$ and every $R>0$ the following holds
for almost every $\omega\in S^1$: the map $\rho\mapsto v(x+\rho\,\omega)$ is
absolutely continuous on $[0,R]$, and
$$
v(x+\rho\,\omega)-v(x)=\int_0^\rho\nabla v(x+s\,\omega)\,\omega\,\dd s
\qquad\hbox{for every }\rho\in[0,R].
$$
\end{lemma}

\begin{proof}
In polar coordinates centered at $x$, for every $f\in \LL^p(B_R(x))$,
$$
\int_{S^1}\int_0^R|f(x+s\,\omega)|\,\dd s\,\dd \omega
=\int_{B_R(x)}\frac{|f(w)|}{|w-x|}\,\dd w
\le\|f\|_{\LL^p(B_R(x))}\,\big\||\cdot-x|^{-1}\big\|_{\LL^{p'}(B_R(x))},
$$
and the last norm is finite because $p'=p/(p-1)<2$. Let $v_k\in \CC^\infty$ with
$v_k\to v$ in $\WW^{1,p}(B_R(x))$. By Morrey's inequality, $v_k\to v$ uniformly
on $B_R(x)$. Applying the estimate to $f=\nabla v$ and 
$f=\nabla v_k-\nabla v$, and passing to a subsequence, we find that for almost
every $\omega\in S^1$
$$
\int_0^R|\nabla v(x+s\,\omega)|\,\dd s<\infty
\qquad\hbox{and}\qquad
\int_0^R|\nabla v_k-\nabla v|(x+s\,\omega)\,\dd s\to0 .
$$
For such $\omega$ and all $\rho\in[0,R]$ we pass to the limit in
$$
v_k(x+\rho\,\omega)-v_k(x)=\int_0^\rho\nabla v_k(x+s\,\omega)\,\omega\,\dd s ,
$$
using the uniform convergence on the left-hand side.
\end{proof}
By Lemma~\ref{lem:dendrite}, any two points $a,b\in X$ are joined by a unique
arc $[a,b]\subset X$ \cite{Nadler}. We set
$$
d_1(a,b):={\rm Var}_1\bigl(\psi;[a,b]\bigr),
$$
the $|\cdot|_1$-length of $\psi$ along $[a,b]$ and $d_1(a,a):=0$.

Recall that a metric space $(Y,d)$ is an \emph{$\RR$-tree} if any two points
$a,b\in Y$ are joined by a unique arc $[a,b]\subset Y$, and this arc is a
geodesic segment: that is, it is isometric to the interval $[0,d(a,b)]$. In
particular, $d(a,b)=d(a,c)+d(c,b)$ for every $c\in[a,b]$.
Lemma~\ref{lem:d1}(iii) below shows that the arcs $[a,b]$ of $X$ are geodesic
segments for $d_1$, so that $(X,d_1)$ is an $\RR$-tree.
\begin{lemma}[the taxicab tree]\label{lem:d1} We have the following.

\smallskip{\rm (i)} For every continuous curve $\gamma$ in $\overline{\Omega}$ from $x$ to $y$,
$d_1(q(x),q(y))\le{\rm Var}_1(u\circ\gamma)$.

\smallskip
{\rm (ii)} There is $C=C(\Omega,p)$ such that
$d_1(q(x),q(y))\le C\,\|\nabla u\|_{\LL^p}\,|x-y|^{1-2/p}$ for all $x,y\in \overline{\Omega}$.

\smallskip{\rm (iii)} $(X,d_1)$ is a compact $\RR$-tree with the topology of $X$, and
$\psi:(X,d_1)\to(\RR^2,|\cdot|_1)$ is $1$-Lipschitz.
\end{lemma}

\begin{proof}
(i) If $z$ is an interior point of the arc $[a,b]$, every path from $a$ to $b$
in $X$ passes through $z$. Otherwise, $a$ and $b$ would be in the same
component of the open set $X\setminus\{z\}$; this component is arcwise
connected, and would contain a second arc from $a$ to $b$, contradicting the uniqueness of $[a,b]$. Let
$a=q(x)$, $b=q(y)$, and let $a=\alpha_0,\dots,\alpha_m=b$ be points in order
along $[a,b]$. Choosing successive first-hitting times, we find
$0=t_0\le t_1\le\dots\le t_m=1$ with $q(\gamma(t_k))=\alpha_k$, so that
$$
\sum_k|\psi(\alpha_k)-\psi(\alpha_{k-1})|_1
=\sum_k|u(\gamma(t_k))-u(\gamma(t_{k-1}))|_1\le{\rm Var}_1(u\circ\gamma).
$$
Taking the supremum over the points $\alpha_k$ gives (i).

(ii) Let $\tilde u\in \WW^{1,p}(\RR^2,\RR^2)$ be an extension of $u$. By
Lemma~\ref{lem:rays} applied to $\tilde u$, for every $x\in \overline{\Omega}$ and almost
every $z\in\Omega$ (the points lying on exceptional rays from $x$ form a null
set), $u$ is absolutely continuous on the segment $[x,z]\subset \overline{\Omega}$, and
$$
{\rm Var}_1(u;[x,z])=\int_0^{|z-x|}\big|\nabla u(x+s\,\omega)\,\omega\big|_1\,\dd s,
\qquad \omega=\frac{z-x}{|z-x|} .
$$
Let $r=|x-y|$ and $S=\Omega\cap B_r(x)$. Since $\Omega$ is convex,
$|S|\ge c_\Omega r^2$. Apply (i) to the curves $x\to z\to y$, $z\in S$, and
average over $z$ as in the proof of Morrey's inequality
\cite[Lemma~7.16]{GT}. This gives
$$
d_1(q(x),q(y))\le C\int_{\Omega\cap B_{2r}(x)}\frac{|\nabla u(w)|}{|w-x|}\,\dd w
+C\int_{\Omega\cap B_{2r}(y)}\frac{|\nabla u(w)|}{|w-y|}\,\dd w ,
$$
and H\"older's inequality gives (ii). 

(iii) By (ii), $d_1$ is finite and $q:\overline{\Omega}\to(X,d_1)$ is continuous. Moreover
$d_1(a,b)>0$ for $a\neq b$, since the light map $\psi$ is not constant on
$[a,b]$, and $d_1$ is additive along the arcs, so $d_1$ is a metric. Since
$q$ is a quotient map, the identity from $X$ to $(X,d_1)$ is a continuous
bijection from a compact space onto a Hausdorff space, hence a homeomorphism.
In each arc $[a,b]$ the function $t\mapsto d_1(a,t)$ is continuous and strictly
increasing, so $[a,b]$ is isometric to the segment $[0,d_1(a,b)]$; therefore
$(X,d_1)$ is a tree $\RR$. Finally,
$|\psi(a)-\psi(b)|_1\le{\rm Var}_1(\psi;[a,b])=d_1(a,b)$.
\end{proof}

We now recall the notions related to $\RR$-trees used in the construction of the
recovery sequences, stated for the $\RR$-tree $(X,d_1)$; see
\cite[Chapters~II.1 and~II.2]{BH}.
A \emph{subtree} of $(X,d_1)$ is a subset $T$ with $[a,b]\subset T$ for all
$a,b\in T$; it is \emph{finite} if it is spanned by finitely many points: that is,
$T=\bigcup_{i=2}^m[a_1,a_i]$ for some $a_1,\dots,a_m\in X$. A finite subtree is
compact, being a finite union of arcs. For $t\in T$, the \emph{order} of $t$ in
$T$ is the number of connected components of $T\setminus\{t\}$. If $T$ is a
finite subtree, the set $V$ of its points of order different from $2$
(endpoints and branch points) is finite, and the closures of the connected
components of $T\setminus V$ are finitely many arcs $[v,v']$, with $v,v'\in V$,
called the \emph{edges} of $T$. Each edge is a geodesic segment of length
$d_1(v,v')>0$. Note that in general the edges are not the arcs $[a_1,a_i]$
that span $T$.

Finally, $\RR$-trees are CAT$(0)$ spaces, and every subtree $T$ is convex, since
the geodesic segments joining its points are the arcs $[a,b]$. Hence, if $T$ is
a finite subtree, every $t\in X$ has a unique nearest point $r_T(t)\in T$, the
\emph{nearest-point retraction} $r_T\colon X\to T$ is $1$-Lipschitz
\cite[Proposition II.2.4]{BH}, and $d_1(r_T(t),t)=d_1(t,T)$.
\begin{lemma}\label{lem:tree-net}
Let $T\subset X$ be a finite subtree, and let $\lambda>0$ be the minimal
$d_1$-length of its edges. For $0<\e\le h\le\lambda$ there is a continuous map
$\psi_{\e,h}:T\to\e\Z$ such that
$$
|\psi_{\e,h}(a)-\psi_{\e,h}(b)|\le\Big(1+\frac{2\e}{h}\Big)d_1(a,b),
\qquad
|\psi_{\e,h}(a)-\psi(a)|_1\le h+3\e 
$$
for all $a,b\in T$.
\end{lemma}

\begin{proof}
Subdivide each edge of $T$, of length $L\ge h$, into $n$ arcs of equal
$d_1$-length $\ell:=L/n$, where $n$ is the integer part of $L/h$, so that
$\ell\in[h,2h)$. To each endpoint $z$ of these arcs, associate a vertex
$\beta_z$ of $\e\Z$ with $|\beta_z-\psi(z)|_1\le\e$. For each arc $[z,z']$, let
$v_{zz'}\colon[0,\ell]\to\e\Z$ be the path given by Lemma~\ref{lem:staircase}
from $\beta_z$ to $\beta_{z'}$ on $[0,\ell]$, and set
$$
\psi_{\e,h}(a):=v_{zz'}\big(d_1(z,a)\big)\qquad\hbox{for }a\in[z,z'].
$$
These definitions agree at the common endpoints of the arcs, so
$\psi_{\e,h}$ is continuous. By Lemma~\ref{lem:d1}(iii),
$|\beta_{z'}-\beta_z|_1\le|\psi(z')-\psi(z)|_1+2\e\le\ell+2\e$; hence, by
Lemma~\ref{lem:staircase}(i), for $a,b\in[z,z']$
$$
|\psi_{\e,h}(a)-\psi_{\e,h}(b)|\le\frac{\ell+2\e}{\ell}\,\bigl|d_1(z,a)-d_1(z,b)\bigr|
\le\Big(1+\frac{2\e}h\Big)d_1(a,b).
$$
For general $a,b\in T$, the geodesic $[a,b]$ is the union of finitely many
subarcs of the arcs $[z,z']$, and the first estimate follows from the
additivity of $d_1$ along $[a,b]$.

For the second estimate, let $a\in[z,z']$, $s:=d_1(z,a)$, and $\theta:=s/\ell$.
By Lemma~\ref{lem:staircase}(ii) and $|\cdot|_1\le\sqrt2\,|\cdot|$,
$$
\big|\psi_{\e,h}(a)-(1-\theta)\beta_z-\theta\beta_{z'}\big|_1\le2\e,
\qquad
\big|(1-\theta)(\beta_z-\psi(z))+\theta(\beta_{z'}-\psi(z'))\big|_1\le\e .
$$
Moreover, $\psi$ is $1$-Lipschitz by Lemma~\ref{lem:d1}(iii), and
$d_1(z,a)=s$, $d_1(a,z')=\ell-s$, so that
$$
\big|(1-\theta)\psi(z)+\theta\psi(z')-\psi(a)\big|_1\le(1-\theta)s+\theta(\ell-s)
=2\ell\theta(1-\theta)\le\frac\ell2<h .
$$
Summing up, $|\psi_{\e,h}(a)-\psi(a)|_1\le h+3\e$.
\end{proof}
\begin{proof}[Proof of Theorem~{\rm \ref{thm:limsup-p>2}}]
The idea is to repeat the construction in Lemma \ref{lem:tree-net} for suitable approximating trees $T_\delta$. This allows one to construct recovery sequences for the corresponding $u_\delta$, which are suitable approximations of $u$.

\smallskip

For $\delta>0$ let $\{a_1,\dots,a_m\}$ be a
$\delta$-net of the compact space $(X,d_1)$, and let $T_\delta$ be the union of
the arcs $[a_1,a_i]$, $i=1,\dots,m$. Then $T_\delta$ is a finite subtree, and
$d_1(t,T_\delta)\le\delta$ for every $t\in X$. Let $r_\delta:=r_{T_\delta}$ be the nearest-point retraction; it is $1$-Lipschitz, and $d_1(r_\delta(t),t)\le\delta$ for all $t\in X$.

We now construct recovery maps for the target function $$u_\delta:=\psi\circ r_\delta\circ q,$$ related to $u$ and $T_\delta$. Let $\lambda_\delta$ be the minimal edge length of
$T_\delta$. For $\e$ small, choose $h=h(\e)$ with $\e\le h\le\lambda_\delta$,
$h\to0$ and $\e/h\to0$, and set
$$
v_\e:=\psi_{\e,h}\circ r_\delta\circ q:\overline{\Omega}\to\e\Z .
$$
By Lemmas~\ref{lem:d1} and~\ref{lem:tree-net}, $v_\e$ is continuous. Fix
$\omega\in S^1$. For almost every line $L$ parallel to $\omega$, $u$ is absolutely
continuous on the segment $L\cap\Omega$. For $x,y\in L\cap\Omega$,
Lemma~\ref{lem:tree-net}, the $1$-Lipschitz continuity of $r_\delta$ and
Lemma~\ref{lem:d1}(i) give
$$
|v_\e(x)-v_\e(y)|\le\Big(1+\frac{2\e}h\Big)\,d_1(q(x),q(y))
\le\Big(1+\frac{2\e}h\Big)\int_{[x,y]}|\partial_\omega u|_1\,\dd\mathcal H^1 .
$$
Hence, $v_\e\in \WW^{1,p}(\Omega,\e\Z)$, and for every $\omega$ in a countable dense
subset of $S^1$
$$
|\nabla v_\e\,\omega|\le (1+\tfrac{2\e}h)|\nabla u\,\omega|_1\qquad\hbox{almost everywhere in }\Omega .
$$
By Lemma~\ref{nabla u_e has only one row}, $\nabla v_\e=e_i\otimes w$ almost everywhere,
so that $|\nabla v_\e\,\omega|=|w\cdot\omega|$, while
$$|\nabla u\,\omega|_1=|\nabla u_1\cdot\omega|+|\nabla u_2\cdot\omega|\le N(\nabla u).$$ Both
sides are continuous in $\omega$; taking $\omega=w/|w|$ we obtain
$$
|\nabla v_\e|\le\Big(1+\frac{2\e}h\Big)N(\nabla u)$$
almost everywhere in $\Omega $.
By Lemma~\ref{lem:tree-net}, applied at the point $r_\delta(q(x))\in T_\delta$,
we have $|v_\e-u_\delta|_1\le h+3\e$, so that $v_\e\to u_\delta$ uniformly as
$\e\to0$, and
$$
\Gamma\hbox{-}\limsup_{\e\to0}F^p_\e(u_\delta)
\le\limsup_{\e\to0}\int_\Omega|\nabla v_\e|^p\,\dd x\le\int_\Omega N(\nabla u)^p\,\dd x .
$$
Finally, since $\psi$ is $1$-Lipschitz by Lemma~\ref{lem:d1}(iii),
$|u_\delta-u|_1\le d_1(r_\delta\circ q,q)\le\delta$, so that $u_\delta\to u$
uniformly as $\delta\to0$. Since the $\Gamma$-limsup is lower semicontinuous with
respect to the $\LL^p$-convergence, the claim follows.
\end{proof}

\begin{remark}\rm
(a) Convexity is used to lift maps $\overline{\Omega}\to S^1$, and in the Morrey estimate of
Lemma~\ref{lem:d1}(ii). Lifting only requires that every continuous map
$\overline{\Omega}\to S^1$ be null-homotopic. By Borsuk's separation theorem, this holds if
and only if $\RR^2\setminus \overline{\Omega}$ is connected. The Morrey estimate holds for
Lipschitz domains, with polygonal paths. The theorem therefore extends to
bounded Lipschitz domains such that $\RR^2\setminus\overline\Omega$ is
connected. If $\overline\Omega$ separates the plane, $X$ may contain circles:
for $u(x)=x/|x|$ on an annulus, $X$ is a circle.

(b) For $p=2$ the same proof applies to every continuous $u\in \WW^{1,2}$ with
$|u(\overline{\Omega})|=0$. Lemma~\ref{lem:d1}(ii) is replaced by the finiteness of $d_1$ on a
dense subset of $X$, which follows from Fubini's theorem. The nearest-point
retractions are replaced by first-point retractions onto an increasing
sequence of finite subtrees whose union is dense; these converge uniformly to
the identity. However, $\WW^{1,2}$ maps need not be continuous, and continuous
$\WW^{1,2}$ maps need not satisfy Lusin's condition (N). The general case $p=2$
is open.
\end{remark}

\section{Concluding remarks: the metric point of view}\label{metric remarks}
We finally note that the convergence of families $F^p_\e$ can be related to a convergence on varying metric spaces.
Indeed, for every $\e>0$ we consider the distance $d_\e$ on $\e\Z$ defined by
\begin{equation}\label{eq:deps}
d_\e(u,v)=\min
\Bigl\{\int_0^1|\phi'(t)|\,\dd t: \phi(t)\in \e\Z \hbox{ for all }t\in[0,1], \phi(0)=u,\ \phi(1)=v\Bigr\}.
\end{equation}

\begin{proposition}\label{dis-con}
The distances $d_\e$ converge with uniform continuity to the distance
$$
d(u,v)=|u-v|_1=|u_1-v_1|+|u_2-v_2|;
$$
more precisely, if $(u_\e,v_\e)\to (u,v)$ then $|d_\e(u_\e,v_\e)-d(u,v)|\le |u_\e-u|+|v_\e-v|$
\end{proposition}
\begin{proof}For $u_\e,v_\e \in \e \mathcal{Z}$, there exists a path on $\e \mathcal{Z}$ connecting $u_\e$ and $v_\e$ and so that
$$
d_\e(u_\e,v_\e)= d_1(u_\e,v_\e)
$$
(e.g. by increasing or decreasing monotonically each component when going from $u_\e$ to $v_\e$, or vice versa, the total sum of the horizontal and vertical lengths will be the same), thus for $u,v \in \mathbb{R}^2$
$$
d_\e(u_\e,v_\e)= d_1(u_\e,v_\e) \leq d_1(u_\e,u) + d_1(u,v)+d_1(v,v_\e),
$$
and this gives one inequality. The other is found in the same way:
$$
d_1(u,v)\leq d_1(u,u_\e) + d_1(u_\e,v_\e)+d_1(v_\e,v) = d_1(u,u_\e) + d_\e(u_\e,v_\e)+d_1(v_\e,v),
$$
    and the claim is proved.
\end{proof}

Let $u_\e\in \WW^{1,p}(\Omega , \e\Z)$ with $p\ge 2$ be converging to $u$ and let $z_\e\in \e \mathbb Z^2$ be converging to $z_0$.
We define the scalar function $v_\e\in \WW^{1,p}(\Omega)$ by
\begin{equation}\label{veze}
  v_\e(x)= d_\e(u_\e(x),z_\e).
\end{equation}
Then the convergence of $u_\e$ can be related to that of these scalar functions $v_\e$, whose limit is $v(x)= d(u(x),z_0)$.
Note that $|\nabla v_\e(x)|=|\nabla u_\e(x)|$, so that $F_\e(u_\e)=\int_\Omega|\nabla v_\e|^p\dd x$.
The liminf inequality for $F^p_\e$ in Proposition \ref{lower bound} can be recovered from the lower-semicontinuity inequality along the sequence $v_\e$ by locally optimizing the choice of $z_0$.
This corresponds to the determination of the full Reshetnyak density, which is obtained by
taking the least $g\in \LL^p$ dominating $|\nabla d_\e(u_\e(\cdot),z)|$ for $z$
in a countable dense subset of $\e\Z$. 

The limit in Theorem 
\ref{thm:super:hom} has a natural metric interpretation. For $\e>0$ and
$v,w\in\e\Z$ let $d_\e$
be the length distance of the network defined in \eqref{eq:deps}. By Proposition \ref{dis-con}, $d_\e$
coincides with the restriction to $\e\Z$ of the $\ell^1$ distance
$|v-w|_1$, and $d_\e\to|\cdot|_1$ with uniform continuity. Hence
$(\e\Z,d_\e)$ embeds isometrically in $(\R2 ,|\cdot|_1)$ as an
$\e/2$-dense subset; in particular, the pointed Gromov--Hausdorff distance
between $(\e\Z,d_\e)$ and $(\R2 ,|\cdot|_1)$ is $O(\e)$, and the networks
converge, as metric spaces, to the normed plane $(\R2 ,|\cdot|_1)$. Since $|w|\le|w|_1\le\sqrt2\,|w|$, the class
$\WW^{1,p}(\Omega,\e\Z)$ coincides, with energies comparable uniformly in $\e$,
with the Sobolev space $\WW^{1,p}(\Omega , (\e\Z,d_\e))$ of maps into a metric
space, which can be defined equivalently after Korevaar and Schoen \cite{KS},
Reshetnyak \cite{Reshetnyak} or Ambrosio \cite{Ambrosio}; the equivalence of
these notions for Euclidean domains and complete separable targets is due to
Chiron \cite{Chiron}, and an explicit formula relating the Korevaar--Schoen
$p$-energy to the Lipschitz-composition description is given by Logaritsch
and Spadaro \cite{LS}.

All these energies are built on the approximate metric derivative: for a
Sobolev map $u$ with values in a metric space $(X,d)$, at almost every $x$
there is a seminorm $\apmd u_x$ on $\R2$ such that
$d(u(y),u(x))-\apmd u_x(y-x)=o(|y-x|)$ as $y\to x$, in the approximate sense.
Its existence is due to Kirchheim \cite{Kirchheim-rect} for Lipschitz maps and
to Karmanova \cite{Karmanova} for Sobolev maps; we follow the notation of
Lytchak and Wenger \cite[Section~4]{LW1}. Reshetnyak's energy density is the
maximal stretch $\max_{|v|=1}\apmd u_x(v)$, while the Korevaar--Schoen density
is an average of $(\apmd u_x(v))^p$ over directions; for general targets,
already for $(\R2 ,|\cdot|_1)$, the two energies are genuinely different, a
phenomenon related to the competing notions of area in normed spaces
\cite{LW1,LW2}. In our problem the ambiguity disappears. If
$\nabla u(x)=a\otimes\nu$ with $|\nu|=1$, then
$\apmd u_x(v)=|a|_1\,|\nu\cdot v|$ for the $\ell^1$ distance, a rank-one
seminorm; moreover, every $u\in \WW^{1,p}(\Omega , \e\Z)$ has $\nabla u(x)\in K$
almost everywhere, and a matrix with at most one non-zero row acts
identically in the Euclidean and in the $\ell^1$ norm. We obtain the
following normalization statement.

\begin{remark}[all energies agree on the constrained class]
\label{rem:energies-agree}\rm
Let $u\in \WW^{1,p}(\Omega , \R2 )$ with $\det\nabla u=0$ almost everywhere, regarded as a
map into $(\R2 ,|\cdot|_1)$. Then its Reshetnyak $p$-energy equals
$\int_\Omega \! (|\nabla u_1|+|\nabla u_2|)^p\, \dd x$, and its Korevaar--Schoen
$p$-energy equals $\gamma_p$ times this quantity, where $\gamma_p>0$ depends
only on $p$ and on the normalization chosen for the Korevaar--Schoen energy
(with the averaged normalization,
$\gamma_p=\fint_{S^{1}} \! |v_1|^p\,\dd\mathcal H^{1}(v)$). If moreover
$u\in \WW^{1,p}(\Omega,\e\Z)$, then $|\nabla u_1|+|\nabla u_2|=|\nabla u|$
almost everywhere, so that the Reshetnyak and the Korevaar--Schoen $p$-energies of
$u\colon\Omega\to(\e\Z,d_\e)$ equal $F^p_\e(u)$ and $\gamma_pF^p_\e(u)$,
respectively.
\end{remark}

Consequently, $F^p_{\rm hom}$ is the Reshetnyak $p$-energy of $u$ regarded as a map
into $(\R2 ,|\cdot|_1)$, the metric limit of the networks, and, the
normalizing factor being the same for every competitor, none of the
$\Gamma$-convergence statements depends on the chosen definition of energy.
We also note that $(\R2 ,|\cdot|_1)$ does not satisfy the Euclidean-tangent
condition of \cite{LW1}, yet every map in the domain of $F^p_{\rm hom}$ satisfies it
pointwise, its approximate metric derivative being induced by a degenerate
scalar product: it is the constraint, and not the geometry of the target
alone, that governs the limit.

\bigskip\noindent
{\bf Acknowledgements.} 
Andrea Braides gratefully acknowledges the hospitality of the Mathematical Instute at Oxford, where the problem in this paper
was first discussed a long, long time ago. He also acknowledges funding by the Excellence Project MatMod@TOV (Grant Agreement No.~CUP\_E83C23000330006) of
the Department of Mathematics of the University of Rome Tor Vergata. 

\bigskip\noindent{\bf Data availability.} Not applicable.

\bigskip\noindent{\bf Declarations.}

\bigskip\noindent{\bf Conflict of interests.} The authors declare that they have no conflict of interest.

\bigskip\noindent
{\bf Use of AI tools.} 
A large language model was used in the preparation of this article.
Beyond general assistance with the exposition, the model provided valuable bibliographic suggestions that helped highlight the significance of the results,
as well as the idea of employing tools from graph theory to construct recovery sequences for the case $p > 2$. Apart from this, the overall approach,
the ideas, and the mathematical framework of the article were developed prior to the use of artificial intelligence tools.
All other results, arguments, and techniques presented here were developed exclusively by the authors.
The authors assume full responsibility for the content of the article.

\end{document}